\documentclass[aps,nofootinbib,preprintnumbers]{report}
\usepackage{amsfonts}
\usepackage[T1]{fontenc}
\usepackage{color}
\usepackage{amssymb}
\usepackage{amsmath}
\def\co#1{\smash{\mathop{\longrightarrow}\limits^{#1}}}
\def\eq#1{\smash{\mathop{=}\limits^{#1}}}
\def\N{{\mathbb N}}
\def\Z{{\mathbb Z}}
\def\R{{\mathbb R}}
\def\C{{\mathbb C}}
\def\p{{\bf P}}
\def\E{{\bf E}}
\def\var{{\bf Var}}
\def\re{\Re{\mathbb e}\thinspace}
\def\e{{\rm e}}
\def\d{{\rm d}}
\def\i{{\rm i}}

\newtheorem{theorem}{Theorem}[section]
\newtheorem{proposition}{Proposition}[section]
\newtheorem{corollary}{Corollary}[section]
\newtheorem{remark}{Remark}[section]
\newtheorem{lemma}{Lemma}[section]

\begin{document}

\title {Approximation Theorems Related to the Coupon Collector's Problem\\[18pt]
Ph.D.~Thesis}
\author{by\\[18pt]
Anna P\'osfai\\[30pt]
\emph{Supervisors:} Prof.~S\'andor Cs\"org\H{o} and Prof.~Andrew D.~Barbour\\
[80pt]
Doctoral School in Mathematics and Computer Science\\
University of Szeged\\
Bolyai Institute\\[60pt]}
\date{2010}
\maketitle

\tableofcontents

\chapter{Introduction}

\section{The coupon collector's problem}

The coupon collector's problem is one of the classical problems of probability theory. The simplest and probably original version of the problem is the following: Suppose that there are $n$ coupons, from which coupons are being collected with replacement. What is the probability that more than $t$ sample trials are needed to collect all $n$ coupons? One of the first discussions of the problem is due to P\'olya \cite{Polya}. It is brought up 7 times in Feller \cite{F}. The problem has numerous variants and generalizations. It is related to urn problems and the study of waiting times of various random phenomena (e.g.~\cite{Holst}, \cite{GH}, \cite{BD}), etc.

We shall be interested in the following version of the problem. A coupon collector samples with replacement a set of $n\ge 2$
distinct coupons so that at each time any one of the $n$ coupons is
drawn with the same probability $1/n$. For a fixed integer
$m\in\{0,1,\ldots,n-1\}$, this is repeated until $n-m$ distinct
coupons are collected for the first time. Let $W_{n,m}$ denote the
number of necessary repetitions to achieve this. Thus the random
variable $W_{n,m}$, called the coupon collector's waiting time, can
take on the values $n-m, n-m+1, n-m+2,\ldots$, and gives the number
of draws necessary to have a collection, for the first time, with
only $m$ coupons missing. In particular, $W_{n,0}$ is the waiting
time to acquire, for the first time, a complete collection.

The starting point in the study of the behavior of the distribution of the coupon collector's waiting
time is the well-known equality in distribution (\cite{F}, p.~225)
\begin{equation}\label{Wgeo}
W_{n,m}\;\eq{\cal D}\;X_{n/n}+X_{(n-1)/n}+\cdots+X_{(m+1)/n},
\end{equation}
where $X_{n/n}, X_{(n-1)/n},\ldots,X_{(m+1)/n}$ are independent
random variables with geometric distributions pertaining to the
success probabilities $n/n, (n-1)/n,\ldots$, $(m+1)/n$,
respectively, so that $\p\{X_{k/n} = j\} = \big(1 - \frac{k}{n}\big)^{j-1}\frac{k}{n}$, $j\in\N := \{1,2,\ldots\}$, for every
$k\in\{m+1,\ldots, n\}$.

Since the mean and variance of a geometric random variable with parameter $p$ are $1/p$ and $(1-p)/p^2$ respectively, the mean and variance of the waiting time are
\begin{equation}\label{mu}
\mu_n = \mu_n(m) := \E(W_{n,m})=n\sum_{k=m+1}^{n}{1\over k},
\end{equation}
and
\begin{equation}\label{0var}
\sigma_n^2 = \sigma_n^2(m) :=\var(W_{n,m})=n\sum_{k=m+1}^{n}{{n-k}\over k^2}=n\sum_{k=m+1}^{n-1}{{n-k}\over k^2}.
\end{equation}

\section{Limit theorems in the coupon collector's problem}

Different limit theorems have been proved for the asymptotic distribution of $W_{n,m}$, depending on how $m$ behaves as $n\to\infty$. From now on all asymptotic relations throughout are meant as $n\to\infty$ unless otherwise specified.

The first result was proved by Erd\H{o}s and R\'{e}nyi \cite{ER} for complete collections when $m=0$ for all $n\in\N$, obtaining a limiting Gumbel extreme value distribution:
$$\frac{W_{n,0}-\mu_n}{n}\,\co{\cal D}\,\textrm{Gumbel}(0),$$
where the probability measure $\textrm{Gumbel}(0)$ is defined to be the Gumbel distribution shifted by Euler's constant:
\begin{equation*}
\textrm{Gumbel}(0)\{(-\infty,x]\}= \e^{-\e^{-(x+\gamma)}},\quad x\in\R,
\end{equation*}
where $\gamma=\lim_{n\to\infty}\left(\sum_{k=1}^{n}{{1}\over{k}}-\log n\right) = 0,577215\ldots$.

This result was extended by Baum and Billingsley \cite{BB}, who examined all relevant sequences of $m$. They determined four different limiting distributions:
\begin{enumerate}
  \item[\bf 1.] {\bf Degenerate distribution at 0}\\
    $$\textrm{If } \frac{n-m}{\sqrt{n}}\to 0,\textrm{ then }W_{n,m}-(n-m)\,\co{\cal D}\, 0,$$
    that is the limiting probability measure is concentrated on 0.
  \item[\bf 2.] {\bf Poisson distribution}\\
    $$\textrm{If } \frac{n-m}{\sqrt{n}}\to\sqrt{2\lambda},\textrm{ then }W_{n,m}-(n-m)\,\co{\cal D}\,\mathrm{Po}(\lambda),$$
    where $\mathrm{Po}(\lambda)$ is the Poisson distribution with parameter $\lambda$ defined by $\mathrm{Po}(\lambda)\{k\}=\frac{\lambda^k}{k!}\e^{-\lambda}$, $k=0,1,2,\ldots$.
  \item[\bf 3.] {\bf Normal distribution}\\
    $$\textrm{If }\frac{n-m}{\sqrt{n}}\to\infty \textrm{ and }m\to\infty, \textrm{ then }\frac{W_{n,m}-\mu_n}{\sigma_n}\,\co{\cal D}\, \mathrm{N}(0,1),$$
    where $\mathrm{N}(0,1)$ denotes the standard normal distribution, whose probability density function with respect to the Lebesgue measure is  $\frac{1}{\sqrt{2\pi}}\e^{-x^2/2}$, $x\in{\mathbb R}$.
  \item[\bf 4.] {\bf Gumbel-like distribution}\\
    $$\textrm{If } m\equiv m, \textrm{ then }\frac{W_{n,m}-\mu_n}{n}\,\co{\cal D}\, \textrm{Gumbel}(m),$$
    where we call $\textrm{Gumbel}(m)$ the Gumbel-like distribution with parameter $m$, and define it to be the probability measure with probability density function with respect to the Lebesgue measure
    \begin{equation*}
    {1\over m!}\e^{-(m+1)\left(x+\gamma-\sum_{k=1}^{m}{{1}\over{k}}\right)}\,\e^{-\e^{-\left(x+\gamma-\sum_{k=1}^{m}{{1}\over{k}}\right)}}, \quad x\in\R.
    \end{equation*}
\end{enumerate}

\section{Aims of the thesis}

One of the aims of this thesis is to refine the limit theorems of the previous section. Our basic goal is to approximate the distribution of the coupon collector's appropriately centered and normalized waiting time with well-known measures with high accuracy, and in many cases prove asymptotic expansions for the related probability distribution functions and mass functions. The approximating measures shall be chosen from five different measure families. Three of them -- the Poisson distributions, the normal distributions and the Gumbel-like distributions -- shall be probability measure families whose members occur as limiting laws in the limit theorems of Baum and Billingsley.

The fourth set of measures considered shall be a certain $\{\pi_{\mu,a}: \mu>0, a>0\}$ family of compound Poisson measures which we now define. For each $\mu>0$ and $a>0$ let $\pi_{\mu,a}$ denote the probability distribution of $Z_1+2Z_2$, where $Z_1$ and $Z_2$ are independent random variables defined on a common probability space, $Z_1\sim\mathrm{Po}(\mu)$ and $Z_2\sim\mathrm{Po}(a/2)$. Since $\mathrm{Po}(\lambda)$, $\lambda>0$, has probability generating function $\exp\left\{\lambda(z-1)\right\}$, the probability generating function of $Z_1+2Z_2$ is
\begin{align*}
g(z)&:=\E\left(z^{Z_1+2Z_2}\right)
=\E\left(z^{Z_1}\right)\E\left((z^2)^{Z_2}\right)
=\e^{\mu(z-1)}\e^{\frac{a}{2}(z^2-1)}\\
&=\exp\left\{\left(\mu+\frac{a}{2}\right)\left(\frac{\frac{a}{2}}{\mu+\frac{a}{2}}z^2+\frac{\mu}{\mu+\frac{a}{2}}z-1\right)\right\}.
\end{align*}
By the basic properties of probability generating functions, we see that $Z_1+2Z_2$ does have a compound Poisson distribution, that is, it equals in distribution a random variable of the form $\sum_{k=1}^{N}X_k$, where $N,X_1,X_2,\ldots$ are independent random variables given on a common probability space such that $N$ has Poisson distribution and $X_1,X_2,\ldots$ are identically distributed, namely $N\sim\mathrm{Po}\left(\mu+\frac{a}{2}\right)$ and each $X_k$, $k=1,2,\ldots$, takes on the values 1 and 2 in the proportion $\mu:\frac{a}{2}$.

The fifth set of approximating measures we consider shall be the family of Poisson--Charlier signed measures. For any positive real numbers $\lambda$, $\widetilde{a}^{(1)},\ldots,\widetilde{a}^{(S)}$ and $S\in\N$, the Poisson--Charlier signed measure $\nu=\nu(\lambda,\widetilde{a}^{(1)},\ldots,\widetilde{a}^{(S)})$ is a signed measure concentrated on the nonnegative integers defined by
\begin{equation}\label{Poisson--Charlier_measure}
\nu\{j\}=\mathrm{Po}\{j\}(\lambda)\left(\sum_{r=1}^{S}(-1)^{r}\widetilde{a}^{(r)}C_r(j,\lambda)\right),\quad j\in\N,
\end{equation}
where
\begin{equation}\label{Charlier_polynomial}
C_r(j,\lambda):=\sum_{k=0}^{r}{r\choose k}{j\choose k}k!\lambda^{-2k}
\end{equation}
is the $r$-th Charlier polynomial (\cite{Ch} p.~170).

In the next chapter we explain the basic underlying ideas of the methods used in the proofs of the thesis, and list some of the important results we shall use. Then, each of the following five chapters is dedicated to the approximation of the coupon collector's waiting time with members of one of the five chosen measure families.

The results of Chapter 3 were published in \cite{PCS}, those of Chapter 4 were published in \cite{P_n}. The results of the first three sections of Chapter 5 can be found in \cite{P_p}, some details are contained in \cite{ParXiv}. The results of Chapter 6 were published in \cite{P_cp}.


\chapter{Methods to measure the closeness of probability distributions}

\section{Probability metrics}

There are several ways of defining the distance of two probability distributions. (See e.g.~\cite{GS} and \cite{R}.) Throughout this section let $\mu$ and $\nu$ be two probability measures on the measurable space $({\mathbb R},{\cal B})$, where ${\cal B}$ denotes the $\sigma$-algebra of the Borel sets of the real line. Let $X$ be a real random variable with distribution $\mu$ and distribution function $F$, and let $Y$ be a real random variable with distribution $\nu$ and distribution function $G$. For an arbitrary family ${\cal H}$ of bounded real measurable functions on the real line we define
\begin{equation}\label{d_H}
d_{\cal H}(\mu,\nu)=\sup_{h\in{\cal H}}\left|\int_{-\infty}^{\infty}\!h\,\d\mu-\int_{-\infty}^{\infty}\!h\,\d\nu\right|
=\sup_{h\in{\cal H}}\left|\E(h(X))-\E(h(Y))\right|,
\end{equation}
which we call the probability metric associated with the family of test functions ${\cal H}$, if ${\cal H}$ is coarse enough to assure $d_{\cal H}(\mu,\nu)=0 \Rightarrow \mu=\nu$.

In this thesis we shall be interested in the probability metrics resulting from ${\cal H}=\{\textrm{indicator functions of }(-\infty,x], x\in{\mathbb R}\}$ and ${\cal H}=\{\textrm{indicator functions of all Borel sets}\}$, which are known as Kolmogorov distance and total variation distance respectively.

\bigskip

\noindent {\bf Kolmogorov distance}

\bigskip

The Kolmogorov distance between $\mu$ and $\nu$ is defined to be
\begin{equation}\label{d_K_def}
d_{\mathrm{K}}(\mu,\nu)=\sup_{x\in{\mathbb R}}\big|\mu((-\infty,x])-\nu((-\infty,x])\big|.
\end{equation}
Clearly, this is exactly the supremum distance of the corresponding distribution functions: $d_{\mathrm{K}}(\mu,\nu)=\sup_{x\in{\mathbb R}}|F(x)-G(x)|$.

Obviously, $0\leq d_{\mathrm{K}}(\mu,\nu)\leq 1$, $d_{\mathrm{K}}(\mu,\nu)=0$ iff $\mu=\nu$, and $d_{\mathrm{K}}(\mu,\nu)=1$ iff $\sup\{x\in{\mathbb R}: F(x)<1\}\leq\inf\{x\in{\mathbb R}: G(x)>0\}$ or $\sup\{x\in{\mathbb R}: G(x)<1\}\leq\inf\{x\in{\mathbb R}: F(x)>0\}$. Since $d_{\mathrm{K}}(\cdot,\cdot)$ as a function of two variables is also symmetric and satisfies the triangle inequality, $d_{\mathrm{K}}(\cdot,\cdot)$ is indeed a metric on the space of probability distributions on $({\mathbb R},{\cal B})$.

Convergence in Kolmogorov metric is stronger than convergence in distribution, that is if $\mu_n$, $n\in{\mathbb N}$, is a sequence of probability measures on $({\mathbb R},{\cal B})$ with corresponding distribution functions $F_n$, $n\in{\mathbb N}$, such that $d_{\mathrm{K}}(\mu_n,\mu)\to0$, then $\mu_n$ converges weakly to $\mu$, meaning that $F_n(x)\to F(x)$ for each $x\in{\mathbb R}$ continuity point of $F$. The converse is not true in general. (See \cite{GS} p.~14 Theorem 6) One possible metric that metrizes weak convergence of probability measures on $({\mathbb R},{\cal B})$ is the Levy metric defined by
\begin{equation*}
d_{\mathrm{L}}(\mu,\nu)=\inf\{\varepsilon>0: F(x-\varepsilon)-\varepsilon\leq G(x) \leq F(x+\varepsilon)+\varepsilon, \forall x\in{\mathbb R}\}.
\end{equation*}
We have
\begin{equation*}
d_{\mathrm{L}}(\mu,\nu)\leq d_{\mathrm{K}}(\mu,\nu)\leq \left(1+\sup_{x\in{\mathbb R}}|F'(x)|\right)d_{\mathrm{L}}(\mu,\nu),
\end{equation*}
where the first inequality is true for any choices of $\mu$ and $\nu$ (see \cite{H} p.~34), while the second one holds true only if $\mu$ is absolutely continuous with respect to the Lebesgue measure (see \cite{P_limit} p.~43). This implies that the weak convergence $\mu_n\rightrightarrows\mu$ is equivalent to $d_{\mathrm{K}}(\mu_n,\mu)\to0$ in the case when the limiting probability distribution $\mu$ is absolutely continuous and has a bounded density function.

\bigskip

\noindent {\bf Total variation distance}

\bigskip

The total variation distance between $\mu$ and $\nu$ is defined to be
\begin{equation}
d_{\mathrm{TV}}(\mu,\nu)=\sup_{B\in{\cal B}}\big|\mu(B)-\nu(B)\big|.
\end{equation}\label{d_TV_def}
The definition above may be given in other equivalent forms. By \cite{BHJ} p.~253, if $\mu$ and $\nu$ are both absolutely continuous with respect to a $\sigma$-finite measure $\lambda$ (for example $\lambda=\mu+\nu$), and $f$ and $g$ are the densities of $\mu$ and $\nu$ with respect to $\lambda$, then
\begin{align*}
d_{\mathrm{TV}}(\mu,\nu)&=|\mu(B_0)-\nu(B_0)|,\quad\textrm{where }B_0=\{x\in{\mathbb R}:f(x)>g(x)\}\\
&=\frac{1}{2}\int_{-\infty}^{\infty}\!|f-g|\,\d\lambda\\
&=1-\int_{-\infty}^{\infty}\!\min\{f,g\}\,\d\lambda.
\end{align*}

Later we shall be interested in the case when $\mu$ and $\nu$ are the distributions of certain integer valued random variables $X$ and $Y$ defined on a common probability space $(\Omega,{\cal A},\p)$. If we choose $\lambda$ to be the $\sigma$-finite measure that puts unit mass on each of the integers, the formulas above yield
\begin{align*}
d_{\mathrm{TV}}(\mu,\nu)&=|\p(X\in B_0)-\p(Y\in B_0)|,\quad\textrm{where }B_0=\{k\in{\mathbb Z}:\p(X=k)>\p(Y=k)\}\\
&=\frac{1}{2}\sum_{k\in{\mathbb Z}}|\p(X=k)-\p(Y=k)|\\
&=1-\sum_{k\in{\mathbb Z}}\min\{\p(X=k),\p(Y=k)\}.
\end{align*}

As in the case of the Kolmogorov metric, it is easy to see that $0\leq d_{\mathrm{TV}}(\mu,\nu)\leq 1$, $d_{\mathrm{TV}}(\mu,\nu)=0$ iff $\mu=\nu$, and $d_{\mathrm{TV}}(\mu,\nu)=1$ iff $\mu$ and $\nu$ are mutually singular. Since $d_{\mathrm{TV}}(\cdot,\cdot)$ as a function of two variables is obviously symmetric and satisfies the triangle inequality, $d_{\mathrm{TV}}(\cdot,\cdot)$ is indeed another metric on the space of probability distributions on $({\mathbb R},{\cal B})$.

We note that $d_{\mathrm{K}}(\mu,\nu)\leq d_{\mathrm{TV}}(\mu,\nu)$. It follows by our remarks concerning convergence in Kolmogorov distance that convergence in total variation distance is stronger than convergence in distribution, that is if $\mu_n$, $n\in{\mathbb N}$, is a sequence of probability measures on $({\mathbb R},{\cal B})$ such that $d_{\mathrm{TV}}(\mu_n,\mu)\to0$, then $\mu_n$ converges weakly to $\mu$. It is easy to give an example showing that the converse of this statement is not true in general. We may take any sequence of discrete real valued random variables for which the central limit theorem holds true. In this case although the induced probability measures converge weakly to the standard normal distribution, the corresponding total variation distances all equal 1, simply because an absolutely continuous and a discrete probability measure are always mutually singular. However, if the probability distributions $\mu_n$, $n\in{\mathbb N}$, and $\mu$ are concentrated on a countable subset of ${\mathbb R}$, then $\mu_n\rightrightarrows\mu$ implies $d_{\mathrm{TV}}(\mu_n,\mu)\to0$ (see \cite{GS} p.~14 Theorem 6).


\section{The method of characteristic functions}

Let $X$ be a real random variable defined on a probability space $(\Omega,{\cal A},\p)$ with distribution $\mu=\p\circ X^{-1}$ and distribution function $F(x)=\mu((-\infty,x])$, $x\in\R$. The characteristic function of $X$ is the complex valued function
\begin{equation*}
\varphi_X(t):=\E\left(\e^{\i tX}\right)=\int_{-\infty}^{\infty}\e^{\i tx}\mu(\d x)=\int_{-\infty}^{\infty}\cos(tx)\mu(\d x)+\i\int_{-\infty}^{\infty}\sin(tx)\mu(\d x),
\end{equation*}
well-defined for all $t\in\R$. The characteristic function of any real random variable completely defines its probability distribution, that is there is a one-to-one correspondence between probability measures on $(\R,{\cal B})$ and characteristic functions (\cite{F}, Volume II.~p.~508).

One of the most important applications of characteristic functions is the study of convergence in distribution. The Continuity Theorem (\cite{F}, Volume II.~p.~508) states that in order that a sequence $\{\mu_n\}_{n\in\N}$ of probability distributions converges weakly to a probability distribution $\mu$ it is necessary and sufficient that the sequence $\{\varphi_n\}_{n\in\N}$ of their characteristic functions converges pointwise to a limit $\varphi$, and that $\varphi$ is continuous at the origin. In this case $\varphi$ is the characteristic function of $\mu$. It follows that if $X_n$ is a real random variable with distribution $\mu_n$, distribution function $F_n$ and characteristic function $\varphi_n$, $n\in\N$, and the same goes for $X$ with $\mu$, $F$ and $\varphi$, then
\begin{equation*}
X_n\to X\quad\Leftrightarrow\quad
\mu_n\rightrightarrows\mu\quad\Leftrightarrow\quad
F_{n}(x)\to F(x), x\in C_F\quad\Leftrightarrow\quad
\varphi_n(t)\to\varphi(t),t\in\R,
\end{equation*}
where $C_F$ denotes the set of continuity points of $F$.

From this it is clear that characteristic functions are an important tool for proving limit theorems. There are also classical results on characteristic functions which provide methods to refine limit theorems. One of these is Esseen's smoothing inequality (\cite{P_sums} p.~109), which applied on two probability distribution functions, gives an upper bound on their Kolmogorov distance with the help of the difference of the corresponding characteristic functions.

\bigskip

\noindent{\bf Esseen's smoothing inequality.}
If $F$ is a nondecreasing function, $G$ is a differentiable function of bounded variation and bounded derivative $g$, $\lim_{x\to-\infty}F(x)=\lim_{x\to-\infty}G(x)$ and $\lim_{x\to\infty}F(x)=\lim_{x\to\infty}G(x)$, and $\varphi$ and $\psi$ are the Fourier-Stieltjes transforms of $F$ and $G$ respectively, that is
$$\varphi(t)=\int_{-\infty}^{\infty}\e^{\i tx}\d F(x)\quad\textrm{and}\quad\psi(t)=\int_{-\infty}^{\infty}\e^{\i tx}\d G(x),\quad t\in\R,$$
then for any $T>0$ we have
\begin{equation*}
\sup_{x\in\R}|F(x)-G(x)|\leq\frac{b}{2\pi}\int_{-T}^{T}\left|\frac{\varphi(t)-\psi(t)}{t}\right|\d t + c_b\frac{\sup_{x\in\R}|g(x)|}{T},
\end{equation*}
where $b>1$ is arbitrary and $c_b>0$ is a constant depending only on $b$.

\bigskip

An analogous result for comparing discrete distributions with the help of their characteristic functions is given in a recent paper of Barbour, Kowalski and Nikeghbali. We gather the results of Proposition 2.2., Corollary 2.3. and the formulas (3.15)--(3.17) on p.~11 in \cite{BKN} in the following theorem.

\bigskip

\begin{theorem}{\em{\bf(Barbour, Kowalski and Nikeghbali)}}\label{t_BKN}
Let $\mu$ and $\nu$ be finite signed measures on $\Z$, with Fourier-Stieltjes transforms $\phi$ and $\psi$ respectively, that is
$$
\phi(t)=\sum_{k\in\Z}\e^{\i tk}\mu\{k\}\quad\textrm{and}\quad\psi(t)=\sum_{k\in\Z}\e^{\i tk}\nu\{k\},\quad t\in\R.
$$
Suppose that $\phi=\widetilde{\phi}\chi$ and $\psi=\widetilde{\psi}\chi$ for some functions $\widetilde{\phi}$, $\widetilde{\psi}$, $\chi:\R\to\C$, and that for some constants $t_0,\gamma_0,\gamma,\rho,\eta>0$ and $\gamma_r,\theta_r>0$, $r=1,2,\ldots,S$, $S\in\N$,
\begin{align}\label{uj_tetel1}
&|\widetilde{\phi}(t)-\widetilde{\psi}(t)|\leq\sum_{r=1}^{S}\gamma_r|t|^{\theta_r}+\gamma_0\quad\textrm{and}\quad
|\chi(t)|\leq\gamma\e^{-\rho t^2},\quad0\leq|t|\leq t_0\\ \label{uj_tetel2}
&|\phi(t)-\psi(t)|\leq\eta,\quad t_0<|t|\leq\pi.
\end{align}
Then
\begin{equation}\label{uj_BKN_loc}
\sup_{k\in\Z}|\mu\{k\}-\nu\{k\}|
\leq\sum_{r=1}^{S}\alpha_{\theta_r}\gamma\gamma_r(\rho\vee1)^{-\frac{\theta_r+1}{2}}+\alpha_1\gamma\gamma_0+\alpha_2\eta,
\end{equation}
furthermore, if $\mu$ is a probability measure, then also
\begin{equation}\label{uj_BKN_Kol}
d_{\mathrm{K}}(\mu,\nu)\leq\inf_{a\leq b}\bigg(2|\nu|\{(-\infty,a)\cup(b,\infty)\}+2\varepsilon_{ab}\bigg)
\end{equation}
and
\begin{equation}\label{uj_BKN_d_TV}
2\,\d_{\mathrm{TV}}(\mu,\nu)\leq\inf_{a\leq b}\bigg((b-a+1)\sup_{k\in\Z}|\mu\{k\}-\nu\{k\}|+6|\nu|\{(-\infty,a)\cup(b,\infty)\}+4\varepsilon_{ab}\bigg),
\end{equation}
where $\alpha_{\theta_r}$ are positive constants depending on $\theta_r$, $r=1,2,\ldots,R$, $\alpha_1=\frac{t_0}{\pi}\wedge \frac{1}{2\sqrt{\pi\rho}}$, $\alpha_2=1-\frac{t_0}{\pi}$ and
\begin{equation}\label{uj_BKN_epsilon}
\varepsilon_{ab}:=
\sum_{r=1}^{S}\alpha_{\theta_r}\gamma\gamma_r(\rho\vee1)^{-\frac{\theta_r}{2}}+(b-a+1)(\alpha_1\gamma\gamma_0+\alpha_2\eta).
\end{equation}
\end{theorem}


\section{Stein's method}

Stein's method is a way of deriving explicit estimates for the closeness of two probability distribution. It was introduced by Charles Stein for normal approximation in \cite{S} in 1972. We shall now review the basic idea of the method (see \cite{BChen}). Let $\mu_0$ be a fixed probability measure on $({\mathbb R},{\cal B})$, which we shall approximate with another probability measure $\mu$ on $({\mathbb R},{\cal B})$. The error of the approximation will be measured in the probability metric $d_{\cal H}(\cdot,\cdot)$ defined in (\ref{d_H}), where ${\cal H}$ is a well-chosen fixed family of test functions. Stein's method consists of the following three steps:
\begin{enumerate}
  \item The Stein characterization of $\mu_0$

    \smallskip
    One needs to find a set of functions ${\cal F}_0\subset{\cal F}:=\{f:{\mathbb R}\to{\mathbb R}\textrm{ measurable}\}$ and a mapping $T:{\cal F}_0\to{\cal F}$ called the Stein operator for $\mu_0$ such that
    \begin{equation}\label{stein_char}
    \mu=\mu_0\quad\textrm{if and only if}\quad\int_{-\infty}^{\infty}\!Tf\,\d\mu=0\,\,\textrm{ for all }f\in{\cal F}_0,
    \end{equation}
    or equivalently for any real random variable $X$
    \begin{equation*}
    X\sim\mu_0\quad\textrm{if and only if}\quad \E(Tf(X))=0\,\,\textrm{ for all }f\in{\cal F}_0.
    \end{equation*}
  \item Solving the Stein equation

    \smallskip
    For each test function $h\in{\cal H}$ one needs to find a solution $f=f_h\in{\cal F}_0$ of the
    \begin{equation}\label{stein_eq}
    h(x)-\int_{-\infty}^{\infty}\!h\,\d\mu_0=Tf(x),\quad x\in\R,
    \end{equation}
    Stein equation. If such an $f_h$ exists for each test function $h\in{\cal H}$, then writing the solution in the above Stein equation, integrating both sides of the equation with respect to $\mu$ and taking the supremum of the absolute values of both sides over all test functions yields
    \begin{equation*}
    d_{\cal H}(\mu,\mu_0)
    =\sup_{h\in{\cal H}}\left|\int_{-\infty}^{\infty}\!h\,\d\mu-\int_{-\infty}^{\infty}\!h\,\d\mu_0\right|
    =\sup_{h\in{\cal H}}\left|\int_{-\infty}^{\infty}\!Tf_h\,\d\mu\right|,
    \end{equation*}
    that is for any random variable $X\sim\mu$ we obtain
    \begin{equation}\label{stein_dist}
    d_{\cal H}(\mu,\mu_0)=\sup_{h\in{\cal H}}\left|\E(Tf_h(X))\right|.
    \end{equation}
  \item Bounding $\sup_{h\in{\cal H}}\left|\E(Tf_h(X))\right|$

    \smallskip
    We have obtained a formula in (\ref{stein_dist}) that expresses $d_{\cal H}(\mu,\mu_0)$ as a supremum of certain expectations. In the formula the distribution $\mu_0$ is only present implicitly through the choices of the functions $T\!f_h$. To bound $d_{\cal H}(\mu,\mu_0)$ we need to give estimates for the expectations $\E(Tf_h(X))$, which are surprisingly easier to bound than the original defining formula of the distance, if the Stein operator $T$ was chosen in a clever way. We note that being able to give good approximations to the expectations $\E(Tf_h(X))$ depends heavily on the properties of the solutions $f_h$.
\end{enumerate}

We conclude that the key point in the procedure proposed above is to chose a good Stein operator for $\mu_0$: not only does $T$ need to characterize $\mu_0$ as given in (\ref{stein_char}), it also has to yield a Stein equation (\ref{stein_eq}) that has a solution $f_h$ for each test function $h\in{\cal H}$, moreover these solutions need to have nice properties.

In this thesis we shall apply results proved by Stein's method to approximate the appropriate function of the coupon collector's waiting time with a compound Poisson random variable, and we shall apply the method directly to obtain a Poisson approximation error estimate for the same waiting time. For later use, we now collect the basic results on Poisson approximation in total variation distance with Stein's method.

\bigskip

\noindent {\bf Poisson approximation with Stein's method}

\bigskip

Stein's method was first extended to Poisson approximation by Chen in \cite{C}. The theory was further developed by Barbour, Holst, Janson and others (see \cite{BChen} and \cite{BHJ}).

The Stein operator for $\mathrm{Po}(\lambda)$ is
$$
T:{\cal F}_0:=\{f:\Z_+\to\R \textrm{ bounded}\}\to{\cal F}:=\{f:\Z_+\to\R\},\quad (Tf)(k)=\lambda f(k+1)-kf(k).
$$
It can be proved (\cite{BChen} p.~65) that this operator $T$ characterizes $\mathrm{Po}(\lambda)$ in the required way, that is for any probability measure $\mu$ on $\Z_+$
$$
\mu=\mathrm{Po}(\lambda)\quad\Leftrightarrow\quad\int_{\Z_+}Tf\d \mu=0\,\textrm{ for all }\,f\in{\cal F}_0.
$$
It can also be proved (\cite{BChen} p.~66) that for each $h_A$ indicator function of $A\subset\Z_+$, the Stein equation
\begin{equation}\label{stein_egyenlet}
\lambda f(k+1)-kf(k)=h(k)-\int_{\Z_+}Tf\d \mathrm{Po}(\lambda)
\end{equation}
has a solution $f_h=f_A$, and
\begin{equation}\label{stein_mokorlat}
\sup_{k\in\Z_+}|f_h(k)|\leq \min\left\{1,\sqrt{\frac{2}{\e\lambda}}\right\}
\end{equation}
and
\begin{equation*}
\sup_{k\in\Z_+}|f_h(k+1)-f_h(k)|\leq \frac{1-\e^{-\lambda}}{\lambda}.
\end{equation*}
The method yields the formula
\begin{equation}\label{stein_d_TV_formula}
d_{\mathrm{TV}}({\cal D}(X),\mathrm{Po}(\lambda))=\sup_{A\subset\Z_+}|\E\{\lambda f_A(X)-Xf_A(X)\}|.
\end{equation}






\section{Couplings}

Let $X$ and $Y$ be random variables defined on the probability spaces $(\Omega_1,{\cal A}_1,\p_1)$ and $(\Omega_2,{\cal A}_2,\p_2)$ respectively. A coupling of $X$ and $Y$ is a pair of random variables $X'$ and $Y'$ that are defined on the same probability space $(\Omega,{\cal A},\p)$, and such that $X\eq{\cal D}X'$ and $Y\eq{\cal D}Y'$. Clearly the coupling of $X$ and $Y$ only depends on the distribution of these random variables. It will be useful for us to think of a coupling the following way: given two probability distributions $\mu_1$ and $\mu_2$ on $({\mathbb R},{\cal B})$, a coupling of these probability measures means the construction of a probability space $(\Omega,{\cal A},\p)$ and a random vector $(X',Y')$ on this probability space whose $\mu$ distribution on $({\mathbb R}^2,{\cal B}^2)$ has marginals $\mu_1$ and $\mu_2$.

Couplings are used in a vast variety of proofs (see e.g.~\cite{L} and \cite{T}). In each of them the basic underlying idea is to construct a suitable coupling $(X',Y')$ such that $X'$ and $Y'$ have the dependence structure most adequate for handling the problem considered. In this thesis we shall use the coupling method to give estimates for the total variation distance of certain distributions. We now present the basic relation between couplings and total variation distance.

\bigskip

\noindent{\bf The coupling inequality}

\bigskip
If $(X',Y')$ defined on $(\Omega,{\cal A},\p)$ is a coupling of the random variables $X$ and $Y$, then
\begin{equation}\label{coupl_ineq}
d_{\mathrm{TV}}(X,Y)\leq \p(X'\neq Y').
\end{equation}

Since the proof of the coupling inequality is quite simple and short, we include it here:
\begin{align*}
d_{\mathrm{TV}}(X,Y)&=\sup_{B\in{\cal B}}|\p(X'\in B)-\p(Y'\in B)|\\
&\leq \sup_{B\in{\cal B}}|\p(X'\in B,X'=Y')-\p(Y'\in B,X'=Y')|+\\&\quad\quad\quad\quad+\sup_{B\in{\cal B}}|\p(X'\in B,X'\neq Y')-\p(Y'\in B,X'\neq Y')|\\
&\leq \p(X'\neq Y').
\end{align*}
It can be proved that there always exists a coupling for which there is equality in (\ref{coupl_ineq}) (\cite{L} p.~19).

\bigskip

\noindent {\bf Couplings to bound ${\bf d_{ \mathrm{{\bf TV}}}(W,W+1)}$, where ${\bf W}$ is a sum of independent integer valued random variables}

\bigskip

Let $X_1,X_2,\ldots$ be independent integer valued random variables, and let $W_n=\sum_{j=1}^{n}X_j$, $n=1,2,\ldots$. We are interested in estimating the total variation distance $d_{\mathrm{TV}}(W_n,W_n+1)$, and there is a technique involving couplings for this purpose (\cite{L} Chapter 3).

Assume there is a probability space $(\Omega,{\cal A},\p)$ on which random variables $X'_1,X'_2,\ldots$ and $X''_1,X''_2,\ldots$ are defined in such a way that both of the sequences consist of independent random variables and $X_j\,\eq{\cal D}\,X_j'\,\eq{\cal D}\,X_j''$, $j=1,2,\ldots$. We identify each of these sequences with a random walk on the integers: let $W'=(W'_0,W'_1,\ldots)$ be a random walk that starts form 0 at the initial moment and has consecutive step sizes $X'_1,X'_2,\ldots$, that is
$$W_0'=0 \textrm{ and } W_n'=\sum_{j=1}^{n}X_j',\, n=1,2,\ldots,$$
and let $W''=(W''_0,W''_1,\ldots)$ be a random walk that starts form 1 at the initial moment and has consecutive step sizes $X''_1,X''_2,\ldots$, that is
$$W_0''=1 \textrm{ and } W_n''=1+\sum_{j=1}^{n}X_j'',\, n=1,2,\ldots.$$
Consider the random time
\begin{equation*}
T=\inf\{k: W_k'=W_k''\}
\end{equation*}
when the random walks $W'$ and $W''$ first meet. (We use the convention that the infimum of the empty set is infinity.) Put $\widetilde{W}''=(\widetilde{W}''_0,\widetilde{W}''_1,\ldots)$, where
\begin{equation*}
\widetilde{W}_k''=\left\{
    \begin{array}{ll}
      W''_k, & \hbox{$k\leq T$,} \\
      W'_k, & \hbox{$k>T$,}
    \end{array}
  \right.
\end{equation*}
for each $k\in{\mathbb N}$.
It is evident that $W''\,\eq{\cal D}\,\widetilde{W}''$, particularly $W_n''\,\eq{\cal D}\,\widetilde{W}''_n$. Since $(W_n',\widetilde{W}''_n)$ is a coupling of $(W_n,W_n+1)$, by the coupling inequality we have
\begin{equation}\label{d_TV_Tn}
d_{\mathrm{TV}}(W_n,W_n+1)\leq \p(W_n'\neq \widetilde{W}_n'')=\p(T>n).
\end{equation}
We see that if $\p(T<\infty)=1$, then $d_{\mathrm{TV}}(W_n,W_n+1)\to0$ as $n\to\infty$, and calculating $\p(T>n)$ yields a bound for the rate of convergence.

It is very important to note that we did not impose any condition on the relation between the random walks $W'$ and $W''$. They can be independent, but they can also have any kind of dependence structure, we only required them to have the same step size distributions. Usually the goal is to define for each $j=1,2,\ldots$ the joint distribution of the step sizes $X_j'$ and $X_j''$ in a way that ensures the finiteness and possibly the minimality of $T$. In other words, one would like to construct couplings of the pairs $(X_j',X_j'')$, $j=1,2,\ldots$, which guarantee that the random walks $W'$ and $W''$ should meet soon, and therefore that $\p(T>n)$ is small. One of the ways to do this is given by the so-called Mineka coupling (\cite{L} p.~44), which we now define.

Fix an arbitrary $j\in\{1,2,\ldots\}$. Set $p_{j,i}=\p(X_j=i)$, $i\in\Z$. We define the distribution of the steps $(X_j',X_j'')$ in $\Z^2$ by
\begin{align*}
&\p\!\left((X_j',X_j'')=(i-1,i)\right)=\frac{1}{2}\min\{p_{j,i-1},p_{j,i}\},\\
&\p\!\left((X_j',X_j'')=(i,i-1)\right)=\frac{1}{2}\min\{p_{j,i-1},p_{j,i}\},\\
&\p\!\left((X_j',X_j'')=(i,i)\right)=p_{j,i}-\frac{1}{2}\min\{p_{j,i-1},p_{j,i}\}-\frac{1}{2}\min\{p_{j,i},p_{j,i+1}\}.
\end{align*}
Thus the couplings force the two random walks to run at most distance 1 apart, in fact, $\{S_k:=W_k'-W_k''\}_{k\in\N}$ defines a symmetric random walk, that starts from $-1$ at time 0, and at each step either stays in place or increases or decreases by 1.

We only calculate the bound for $\p(T>n)$ resulting form the Mineka coupling in the case when the $X_j$, $j=1,2,\ldots$, are iid random variables with discrete uniform distribution on $\{1,2,\ldots,L\}$ for some integer $L\geq2$, that is $p_{j,i}=\frac{1}{L}$, $i=1,2,\ldots,L$, $j=1,2,\ldots$. In this case $$\p(X_j'-X_j''=1)=\p(X_j'-X_j''=-1)=\frac{1}{2}\sum_{i\in\Z}\min\{p_{j,i-1},p_{j,i}\}=\frac{L-1}{2L}$$
and
$$\p(X_j'-X_j''=0)=1-\sum_{i\in\Z}\min\{p_{j,i-1},p_{j,i}\}=\frac{1}{L}.$$
Using the properties of $\{S_k\}_{k\in\N}$, namely that it evolves by unit steps, that the reflection principle can be applied to it, and that $S_n$ has symmetric distribution around $-1$, we obtain
\begin{align*}
\p(T\leq n)&=\p(\max_{0\leq k\leq n}S_k\geq 0)\\
&=\p(\max_{0\leq k\leq n}S_k\geq 0,S_n=0)+\p(\max_{0\leq k\leq n}S_k\geq 0,S_n<0)+\p(\max_{0\leq k\leq n}S_k\geq 0,S_n>0)\\
&=\p(\max_{0\leq k\leq n}S_k\geq 0,S_n=0)+2\p(\max_{0\leq k\leq n}S_k\geq 0,S_n>0)\\
&=\p(S_n=0)+2\p(S_n>0)\\
&=\p(S_n=0)+\p(S_n>0)+\p(S_n<-2)\\
&=1-\p(S_n=-1).
\end{align*}
Thus $\p(T>n)\leq\max_{i\in\Z}\p(S_n=i)$, and by Lemma 4.7 of Barbour and Xia \cite{BX}, we have
$$\p(T>n)\leq\frac{1}{2}\left(n\min\left\{\frac{L-1}{L},\frac{1}{2}\right\}\right)^{-\frac{1}{2}}=\frac{1}{\sqrt{2n}}.$$
It follows by (\ref{d_TV_Tn}) that if $W_n$ is a sum of $n$ independent uniformly distributed random variables on $\{1,2,\ldots,L\}$, then
\begin{equation}\label{Mineka_uniform}
d_{\mathrm{TV}}(W_n,W_n+1)\leq\frac{1}{\sqrt{2n}}.
\end{equation}
We shall see in Section 6.1 that this inequality can be improved.


\chapter{Gumbel-like approximation}

\section{Preliminaries and results}

In this chapter we are interested in the case of the coupon collector's problem when $n\ge m+1$ is large compared to $m$, so we fix a non-negative integer $m$, and we shall look at the asymptotic behavior of the distribution function
\begin{equation}\label{g_Fnm}
F_{n,m}(x) := \p\Bigg( {{1}\over{n}}\,W_{n,m}-\sum_{k=m+1}^n
\!{1\over k}\le x \Bigg), \quad x\in\R,
\end{equation}
as $n\to\infty$.

As mentioned in the introduction, in 1961 Erd\H{o}s and R\'{e}nyi \cite{ER} proved for the case
$m=0$, a full collection, that the limiting distribution is the
Gumbel extreme value distribution, shifted by Euler's constant
$\gamma=\lim_{n\to\infty}\left(\sum_{k=1}^{n}{{1}\over{k}}-\log
n\right) = 0,577215\ldots$, so that
$$
\lim_{n\to\infty}F_{n,0}(x) = F_0(x) := \e^{-\e^{-(x+\gamma)}},
\quad x\in\R.
$$
For an arbitrary non-negative integer $m$, this beautiful result was
extended by Baum and Billingsley \cite{BB} shortly thereafter, who proved
that
\begin{align*}
\lim_{n\to\infty}F_{n,m}(x) = F_m(x) &:= {1\over
m!}\int_{-\infty}^{x} \e^{-(m+1)(y + C_m)}\,\e^{-\e^{-(y+ C_m)}}\d
y\cr &= {1\over m!}\int_{-\infty}^{x + C_m}
\e^{-(m+1)y}\,\e^{-\e^{-y}}\,\d y, \quad x\in\R,
\end{align*}
where $C_{m} := \gamma-\sum_{k=1}^{m}{{1}\over{k}}$. Much later
Cs\"{o}rg\H{o} \cite{CS} refined this general result,
proving\goodbreak\noindent that the rate of convergence in it is
surprisingly fast, namely
\begin{equation}\label{g_1.2}
\sup_{x\in\R}\big|F_{n,m}(x) - F_m(x)\big|\le D_m\,{{\log
n}\over{n}}
\end{equation}
for some constant $D_m>0$ depending only on $m$.

In this thesis, for every $m$ we give a one-term asymptotic expansion
$F_m(\cdot) + G_{n,m}(\cdot)$ that approximates $F_{n,m}(\cdot)$
with the uniform order of $1/n$ such that the explicit sequence of
functions $G_{n,m}(\cdot)$ has the uniform order of $(\log n)/n$. In
particular, it follows that the rate of convergence in (\ref{g_1.2}) can not
be improved.

To introduce $G_{n,m}(\cdot)$, consider the density function of the
limiting distribution:
$$
f_{m}(x) := F_{m}^{\prime}(x) =
{{1}\over{m!}}\,\e^{-\e^{-(x+C_{m})}}\e^{-(m+1)(x+C_{m})} =
{{1}\over{m!}}\,\e^{-e_{m}(x)}e_{m}^{m+1}(x),\quad x\in\R,
$$
where $e_m(x) := \e^{-(x+C_{m})}$, whose second derivative by simple
calculation is
\begin{align}\label{g_1.3}
f_{m}^{\prime\prime}(x)
&= f_m(x)\big[e_m^2(x) -(2m+3)e_m(x) + (m+1)^2\big]\cr
&= {\e^{-e_{m}(x)}\over m!}\,\big[e_m^{m+3}(x) - (2m+3)e_m^{m+2}(x) + (m+1)^2e_{m}^{m+1}(x)\big]
\end{align}
for all $x\in\R$. For every $k\in\N$, consider also the density
function
\begin{equation}\label{g_1.4}
h_k(x):= \left\{
           \begin{array}{ll}
             k\,\e^{-kx}, & \hbox{if $x\geq0$;} \\
             0, & \hbox{if $x<0$}
           \end{array}
         \right.
\end{equation}
of the exponential distribution with mean $1/k$, and the convolution
\begin{equation}\label{g_1.5}
[f_{m}^{\prime\prime}\!\star h_k](x) =
\int_{0}^{\infty}\!f_{m}^{\prime\prime}(x-y)h_k(y)\,\d y =
\int_{-\infty}^x\! h_k(x-y)f_{m}^{\prime\prime}(y)\,\d y, \quad
x\in\R.
\end{equation}
Then, for $n\ge m+2$, our basic sequence of functions will be
\begin{equation}\label{g_1.6}
G_{n,m}(x) = - {{1}\over{2n}}\sum_{k=m+1}^{n-1}{{1}\over{k}}
\int_{-\infty}^{x}\![f_{m}^{\prime\prime}\!\star h_{k}](u)\,\d u,
\quad x\in\R.
\end{equation}

It is natural to consider the following version of the
Baum--Billingsley theorem:
$$
F_m^{*}(x):=\lim_{n\to\infty}F_{n,m}^{*}(x)\quad\hbox{where}\quad
F_{n,m}^{*}(x) := \p\bigg( {{1}\over{n}}\,W_{n,m}-\log n \leq x
\bigg),
$$
so that, clearly,
$$
F_{n,m}^{*}(x) =
F_{n,m}\!\Bigg(x-\Bigg[\sum_{k=1}^{n}{{1}\over{k}}-\log n\Bigg] +
\sum_{k=1}^{m}{{1}\over{k}}\Bigg),\quad x\in\R,
$$
for all $n\geq m+2$, and hence
$$
F_m^{*}(x) = F_m\!\Bigg(x - \gamma + \sum_{k=1}^{m}{1\over k}\Bigg)
= F_m(x - C_m),\quad x\in\R.
$$
For every $n\ge m+2$, the corresponding version of the function in
(\ref{g_1.6}) is
$$
G_{n,m}^{*}(x) = G_{n,m}\!\Bigg(x-\Bigg[\sum_{k=1}^{n}{{1}\over{k}}
-\log n\Bigg]+\sum_{k=1}^{m}{{1}\over{k}}\Bigg),\quad x\in\R.
$$
With all asymptotic relations meant throughout as $n\to\infty$
unless otherwise specified, our main result is the following
\bigskip

\begin{theorem}\label{tetel_gumbel}
For every fixed $m\in \{0,1,2,\ldots\}$,
\begin{equation}\label{g_1.7}
\sup_{x\in\R}\big|F_{n,m}(x) - [F_{m}(x)+G_{n,m}(x)]\big| =
O\!\bigg({{1}\over{n}}\bigg),
\end{equation}
for the functions $G_{n,m}(\cdot)$ given in (\ref{g_1.6}), for which
there exist a constant $K_m>0$, a point $x_m\in\R$, a positive
function $c_{m}(\cdot)$ and a threshold function $n_m(\cdot)\in\N$,
all depending only on $m$, such that
\begin{equation}\label{g_1.8}
\sup_{x\in\R}\,\big|G_{n,m}(x)\big| \le K_m\,{\log n\over n},\quad
n\ge m+2,
\end{equation}
but
\begin{equation}\label{g_1.9}
\big|G_{n,m}(x)\big|\ge c_{m}(x)\,{\log n\over n} \quad\hbox{for
all}\quad x\in(-\infty, x_m),
\end{equation}
whenever $n\ge n_m(x)$. Furthermore,
\begin{equation}\label{g_1.10}
\sup_{x\in\R}\big|F_{n,m}^{*}(x) - [F_{m}^{*}(x)+G_{n,m}^{*}(x)]\big|
= O\!\bigg({{1}\over{n}}\bigg),
\end{equation}
where the sequence $\{G_{n,m}^{*}(\cdot)\}_{n=m+1}^{\infty}$ of
functions has the same properties as the sequence
$\{G_{n,m}(\cdot)\}_{n=m+1}^{\infty}$ in the first statement.
\end{theorem}

\bigskip

We finish this section by examining the optimality of the results of Theorem \ref{tetel_gumbel}.

It is easy to give an argument showing that any sequence of discrete probability laws corresponding to some random variables $X_n$ with finite second moments, can not be approximated with an absolutely continuous distribution in Kolmogorov distance with an error order that is smaller than
$1/d_n$, where $d_n=\#\{x\in[-2\sqrt{\var X_n},2\sqrt{\var X_n}]: \p(X_n=x)>0\}$. In order to prove this, we may assume without loss of generality that $\E(X_n)=0$ for each $n\in\N$. For a fixed $n\in\N$, by Chebisev's inequality, $\p(|X_n|\geq2\sqrt{\var X_n})\leq1/4$, which implies that $X_n$ maps into the interval $[-2\sqrt{\var X_n},2\sqrt{\var X_n}]$ with probability at least $3/4$, and thus there exists a point $x_n\in[-2\sqrt{\var X_n},2\sqrt{\var X_n}]$ for which $\p(X_n=x_n)\geq3/(4d_n)$. This means that the distribution function of $X_n$ has a jump at least $3/(4d_n)$ big at $x_n$, and hence can not be approximated at that point with a continuous function any better than $3/(8d_n)$.

Now, the distribution function $F_{n,m}$ defined in (\ref{g_Fnm}) corresponds to the discrete random variable $(W_{n,m}-\mu_n)/n$ for which $d_n$ is of order $n$, because one can calculate that $\sigma_n\sim c_mn$ with some $c_m$ constant depending only on our fixed $m$. It follows that the supremum distance of $F_{n,m}$ to any continuous function, in particular $F_{m}+G_{n,m}$, can not decrease in a faster order than $1/n$, as $n\to\infty$. This not only proves that the error order in (\ref{g_1.7}) is sharp, but also that no longer asymptotic expansion of $F_{n,m}$ than the one given by (\ref{g_1.7}) can improve the current error order $1/n$.

\section{Proofs}

We shall now verify the theorem above. Before embarking on the proof of (\ref{g_1.7}), we analyze the function
$G_{n,m}(\cdot)$ defined in (\ref{g_1.6}), to show in particular that this
formula makes sense, and prove its properties stated in (\ref{g_1.8}) and
(\ref{g_1.9}). We begin with claiming that for every $l\in\N$,
\begin{equation}\label{g_2.1}
\int_{-\infty}^{x}\!\!\e^{-e_{m}(u)}e_{m}^{\,l}(u)\,\d u =
\!\int_{e_{m}(x)}^{\infty}\!\!\e^{-v} v^{l-1}\,\d v =
\e^{-e_{m}(x)}(l-1)!\sum_{j=0}^{l-1}{{e_{m}^{j}(x)}\over{j!}},\;\,
x\in\R.
\end{equation}
for the functions $e_m(x) = \e^{-(x+C_m)}$ in (\ref{g_1.3}). Indeed, this is
true for $l=1$, and since
$$
\int_{e_{m}(x)}^{\infty}\!\e^{-v} v^{k}\,\d v =
\e^{-e_{m}(x)}e_{m}^{\,k}(x) + k\int_{e_{m}(x)}^{\infty}\!\e^{-v}
v^{k-1}\,\d v,
$$
it follows for $l=k+1$ if it holds for $k\in\N$. So, (\ref{g_2.1}) follows
by induction. Also,
\begin{equation}\label{g_2.2}
\int_{-\infty}^{\infty}\!\e^{-e_{m}(x)}e_{m}^{\,l}(x)\,\d x =
\Gamma(l) = (l-1)!\quad\hbox{for all}\quad l\in\N,
\end{equation}
and we see from (\ref{g_1.3}) that $f_{m}^{\prime\prime}(\cdot)$ is
integrable on $\R$; in fact,
$\int_{-\infty}^{\infty}\!f_{m}^{\prime\prime}(x)\,\d x = 0$.

As is well known, the convolution of two integrable functions is
integrable. Since for our convolution, from (\ref{g_1.5}),
\begin{equation}\label{g_2.3}
[f_{m}^{\prime\prime}\!\star h_{k}](x) =
k\,\e^{-kx}\int_{-\infty}^{x}\!\e^{ky}f_{m}^{\prime\prime}(y)\,\d y,
\quad x\in\R,
\end{equation}
we have $\big|[f_{m}^{\prime\prime}\!\star h_{k}](x)\big| \le
k\int_{-\infty}^{x}\!\big|f_{m}^{\prime\prime}(y)\big|\,\d y$, its
integrability follows directly by (\ref{g_1.3}) and (\ref{g_2.1}). The last
inequality also implies that $\lim_{x\to-\infty}[f_{m}^{\prime\prime}\!\star h_{k}](x) = 0$; in fact, since
$\lim_{|x|\to \infty}f_{m}^{\prime\prime}(x) = 0$ directly from
(\ref{g_1.3}), using the dominated convergence theorem in the first formula
in (\ref{g_1.5}) we also see that $\lim_{x\to
\infty}[f_{m}^{\prime\prime}\!\star h_{k}](x) = 0$. The first two of
the last three properties already make (\ref{g_1.6}) meaningful, so that,
substituting (\ref{g_2.3}) into that formula, for the derivative at each
$x\in\R$ we get
\begin{equation}\label{g_2.4}
G_{n,m}^{\,\prime}(x) = -{{1}\over{2n}}\sum_{k=m+1}^{n-1}\!{{1}\over{k}}\,[f_{m}^{\prime\prime}\!\star h_{k}](x)
=-{{1}\over{2n}}\sum_{k=m+1}^{n-1}\!\bigg[\e^{-kx}\int_{-\infty}^{x}\!\e^{kv}f_{m}^{\prime\prime}(v)\,\d v\bigg].
\end{equation}

Next we note that the derivative of the function in (\ref{g_2.1}),
\begin{equation}\label{g_2.5}
\left(\e^{-e_{m}(x)}e_{m}^{\,l}(x)\right)^{\!\prime}
=\e^{-e_{m}(x)}\big[e_{m}^{l+1}(x)-le_{m}^{l}(x)\big],
\end{equation}
is zero at $x_{0} := -\log l - C_{m}$, is positive for $x<x_0$ and
negative for $x>x_0$. Thus
\begin{equation}\label{g_2.6}
\max_{x\in\R}\left\{e^{-e_{m}(x)}e_{m}^{l}(x)\right\} =
e^{-e_{m}(x_{0})}e_{m}^{l}(x_{0}) =
\left({{l}\over{\e}}\right)^{\!l},\quad l\in\N.
\end{equation}

We also see from (\ref{g_2.5}) and (\ref{g_2.2}) that in fact the $j$-th derivative
$f_m^{(j)}(\cdot)$ of $f_m(\cdot)$ is integrable on $\R$ and
$\lim_{|x|\to \infty}f_m^{(j)}(x) = 0$ for every $j\in\{0\}\cup\N$,
not just for $j=0,1,2$, and hence, as an extension of (\ref{g_2.3}), the
convolutions
$$
\big[f_{m}^{(j)}\!\star h_{k}\big](x) =
k\,\e^{-kx}\int_{-\infty}^{x}\!\e^{ky}f_{m}^{(j)}(y)\,\d y, \quad
x\in\R,
$$
make sense as integrable functions for all $j\in\{0\}\cup\N$ and
$k\in\N$.

\bigskip

\noindent{\bf Proof of (\ref{g_1.8}).} To this end, with (\ref{g_1.6}) in mind, by (\ref{g_2.3})
we have
\begin{align*}
\bigg|\int_{-\infty}^{x}\![f_{m}^{\prime\prime}\!\star
h_{k}](u)\,\d u\bigg| &\le \int_{-\infty}^{x}\!k\,\e^{-ku}
\left|\int_{-\infty}^{u}\!\e^{kv}f_{m}^{\prime\prime}(v)\,\d
v\right|\d u\cr &= \int_{-\infty}^{x}\!k\,\e^{-ku}
\left|{\e^{ku}\over k}\,f_{m}^{\prime\prime}(u) -
\int_{-\infty}^{u}\!{\e^{kv}\over{k}}\,
f_{m}^{\prime\prime\prime}(v)\,\d v\right|\d u\cr
&\le\int_{-\infty}^{x}\!\e^{-ku}\left\{\e^{ku}\big|f_{m}^{\prime\prime}(u)\big|
+ \int_{-\infty}^{u}\!\e^{kv}
\big|f_{m}^{\prime\prime\prime}(v)\big|\,\d v\right\}\d u\cr &\le
\int_{-\infty}^{x}\!\big|f_{m}^{\prime\prime}(u)\big|\d u +
\int_{-\infty}^{x}\!
\left[\int_{-\infty}^{u}\!\big|f_{m}^{\prime\prime\prime}(v)\big|\,\d
v\right]\d u
\end{align*}
for all $x\in\R$, regardless of what $k\in\N$ is. Starting from
(\ref{g_1.3}) and using (\ref{g_2.5}), it is clear that
$|f_{m}^{\prime\prime\prime}(v)|$ in the inner integral of the
second term is bounded by a linear combination of functions of the
form $\e^{-e_{m}(v)}e_{m}^{\,l}(v)$, in which all the exponents $l$
and all the coefficients depend only on $m$. Hence, after an
application of (\ref{g_2.1}), that inner integral itself is bounded  by a
similar linear combination, in the variable $u$, that has the same
property. Thus, by another application of (\ref{g_2.1}), the sum of the two
terms is bounded by a linear combination of functions of the form
$\e^{-e_{m}(x)}e_{m}^{\,j}(x)$, in which both the coefficients and
all the exponents $j$ depend only on $m$. But all these functions
are bounded by (\ref{g_2.6}), and hence $\big|\int_{-\infty}^{x}
[f_{m}^{\prime\prime}\!\star h_{k}](u)\,\d u\big|\le 2K_m/3$ for all
$x\in\R$ and $k\in\N$, for some constant $K_m>0$ depending only on
$m$. Therefore, substituting this into the obvious term-wise bound
for (\ref{g_1.6}),
\begin{equation}\label{g_2.7}
\sup_{x\in\R}\big|G_{n,m}(x)\big|\le {K_m\over
3n}\sum_{k=m+1}^{n-1}{1\over k} \le {K_m\over 3n}\bigg[1 + \int_1^n
\!{1\over x}\,\d x\bigg] \le K_m {\log n\over n}
\end{equation}
for all $n\ge m+2\ge 2$, which is (\ref{g_1.8}).

\bigskip

Before attending to the proof of (\ref{g_1.9}), we note that replacing
$f_m^{\prime\prime}$ and $f_m^{\prime\prime\prime}$ by $f_m^{(j)}$
and $f_m^{(j+1)}$, respectively, the argument in the proof of (\ref{g_1.8})
above gives
\begin{equation}\label{g_2.8}
\int_{-\infty}^{x}\!\big|\big[f_{m}^{(j)}\!\star
h_{k}\big](u)\big|\,\d u \le K_m^{(j)}\quad\hbox{for all}\quad
x\in\R \quad\hbox{and}\quad k\in\N,
\end{equation}
for every $j\in\{0\}\cup\N$, where the constant $K_m^{(j)}>0$
depends only on $m$ and $j$.

\bigskip

\noindent{\bf Proof of (\ref{g_1.9}).} Examining the behavior of
$f_{m}^{\prime\prime}(\cdot)$ given in (\ref{g_1.3}), we see that it first
increases from $0 = f_{m}^{\prime\prime}(-\infty)$ on a half-line
and eventually reaches $0 = f_{m}^{\prime\prime}(\infty)$. Thus the
smallest value $x_m\in\R$ where $f_{m}^{\prime\prime}(\cdot)$ has a
local maximum is well defined. Consider any $x$ in the half-line
$(-\infty, x_m]$. Then the convolution $[f_{m}^{\prime\prime}\!\star
h_{k}](x)$, given in (\ref{g_2.3}) is positive since the integrand is
positive on $(-\infty, x)$. Thus the fist two displayed lines in the
proof of (\ref{g_1.8}) above become
\begin{align*}
\int_{-\infty}^{x}\![f_{m}^{\prime\prime}\!\star
h_{k}](u)\,\d u &= \int_{-\infty}^{x}\!f_{m}^{\prime\prime}(u)\,\d u
- \int_{-\infty}^{x}\!\e^{-ku}
\bigg[\int_{-\infty}^{u}\!\e^{kv}f_{m}^{\prime\prime\prime}(v)\,\d
v\bigg]\d u\cr &= \int_{-\infty}^{x}\!f_{m}^{\prime\prime}(u)\,\d u
- {1\over k}\int_{-\infty}^{x}\![f_{m}^{\prime\prime\prime}\!\star
h_{k}](u)\,\d u\cr &\ge
\int_{-\infty}^{x}\!f_{m}^{\prime\prime}(u)\,\d u - {K_m^{(3)}\over
k},
\end{align*}
where the inequality is by (\ref{g_2.8}) and the first term of the lower
bound is positive. Hence, still for the same $x\in(-\infty, x_m]$,
\begin{align*}
\big|G_{n,m}(x)\big| &=
{{1}\over{2n}}\sum_{k=m+1}^{n-1}{{1}\over{k}}
\int_{-\infty}^{x}\![f_{m}^{\prime\prime}\!\star h_{k}](u)\,\d u\cr
&\ge {1\over
2n}\Bigg[\int_{-\infty}^{x}\!f_{m}^{\prime\prime}(u)\,\d u
\sum_{k=m+1}^{n-1}{1\over{k}} - K_m^{(3)}\sum_{k=m+1}^{n-1}{1\over
k^2}\Bigg]\cr &\ge {1\over 2n}
\Bigg[\bigg\{\int_{-\infty}^{x}\!f_{m}^{\prime\prime}(u)\,\d
u\bigg\} \bigg\{\int_{m+1}^{n}{1\over x}\,\d x\bigg\} -
K_m^{(3)}\sum_{k=1}^{\infty}{1\over k^2}\Bigg]\cr &= {\log n\over n}
\Bigg[\bigg\{{1\over
2}\int_{-\infty}^{x}\!f_{m}^{\prime\prime}(u)\,\d u\bigg\} \bigg\{1
- {\log(m+1)\over \log n}\bigg\} - {\pi^2 K_m^{(3)}\over 12\log
n}\Bigg]\cr &=: {\log n\over n}\,c_{n,m}(x).
\end{align*}
Since $\lim_{n\to\infty}c_{n,m}(x) = 2c_m(x)$, where $c_m(x) =
{1\over 4}\int_{-\infty}^{x}\!f_{m}^{\prime\prime}(u)\,\d u > 0$,
there exists a threshold $n_m(x)\in\N$ such that $|G_{n,m}(x)|\ge
c_m(x)(\log n)/n$ whenever $n\ge n_m(x)$, which is the statement in
(\ref{g_1.9}).

\bigskip

We need one more preliminary remark for later use. Noticing that
$$2L_m/3 := \max_{v\in\R}|f_{m}^{\prime\prime}(v)|<\infty$$ by (\ref{g_2.6}),
it is that
\begin{equation}\label{g_2.9}
\sup_{x\in\R}\big|G_{n,m}^{\,\prime}(x)\big| \le {L_m\over
3n}\sum_{k=m+1}^{n-1}\!
\left[\e^{-kx}\int_{-\infty}^{x}\!\e^{kv}\,\d v\right] = {L_m\over
3n}\sum_{k=m+1}^{n-1}{1\over k} \le L_m{\log n\over n}
\end{equation}
for all $n\ge m+2\ge 2$, which is obtained from (\ref{g_2.4}) as in
(\ref{g_2.7}).

\bigskip

\noindent{\bf Proof of (\ref{g_1.7}).} The proof being Fourier-analytic, the
functions $G_{n,m}(\cdot)$ are first identified in terms of Fourier
transforms. We approximate the characteristic function
$\varphi_{nm}(\cdot) = \varphi_{n,m}(\cdot)$ of $F_{n,m}(\cdot)$ by
the sum of the characteristic function $\varphi_{m}(\cdot)$ of the
limiting distribution $F_{m}(\cdot)$ and a ``correcting" function.
We achieve this in four steps. Starting from $\varphi_{nm}(\cdot)$,
we define the functions $\psi_{nm,1}(\cdot)$, $\psi_{nm,2}(\cdot)$,
$\psi_{nm,3}(\cdot)$ and $\psi_{nm,4}(\cdot)$, where each of these
functions is an estimate of the preceding one --- each time obtained
by keeping only some leading terms from the series expansion of an
ingredient ---, and $\varphi_{nm}(\cdot)\approx\psi_{nm,4}(\cdot) =
\varphi_{m}(\cdot)+\psi_{nm}(\cdot)$ holds for some function
$\psi_{nm}(\cdot)$. For $t\in\R$, the approximations yield the error
functions
\begin{align}\label{g_2.10}
\rho_{nm,1}(t) &:= \varphi_{nm}(t)-\psi_{nm,1}(t),\cr
\rho_{nm,2}(t) &:= \psi_{nm,1}(t)-\psi_{nm,2}(t),\cr \rho_{nm,3}(t)
&:= \psi_{nm,2}(t)-\psi_{nm,3}(t),\cr \rho_{nm,4}(t) &:=
\psi_{nm,3}(t)-\psi_{nm,4}(t)
=\psi_{nm,3}(t)-\left[\varphi_{m}(t)+\psi_{nm}(t)\right],
\end{align}
and $G_{n,m}(\cdot)$ will be the function whose Fourier--Stieltjes
transform is $\psi_{nm}(\cdot)$.

First, since the characteristic function of the geometric
distribution with success probability $p\in(0,1)$ is
$\big(1+{{1}\over{p}}\left\{\e^{-{\rm i}t}-1\right\}\big)^{-1}$,
$t\in\R$, where $\rm i$ is the imaginary unit,
\begin{align*}
\varphi_{nm}(t) = \varphi_{n,m}(t) &:=
\int_{-\infty}^{\infty}\!\e^{{\rm i}tx}\,\d F_{n,m}(x) =
\E\Bigg(\exp\!\Bigg\{{\rm i}t\Bigg[{{1}\over{n}}W_{n}(m)
-\sum_{k=m+1}^{n}{{1}\over{k}}\Bigg]\Bigg\}\Bigg)\cr
&=\prod_{k=m+1}^{n-1}{{\e^{-{\rm i}t/k}}\over{1+{{n}\over{k}}
\left(\e^{-{\rm i}t/n}-1\right)}},\quad t\in\R,
\end{align*}
by (\ref{Wgeo}). Also, for all $t\in\R$ the limiting characteristic
function is
\begin{align}\label{g_2.11}
\varphi_m(t) &:= \int_{-\infty}^{\infty}\e^{{\rm i}tx}\,\d
F_m(x) = \E\Bigg(\exp\!\Bigg\{{\rm i}t
\Bigg[\sum_{k=m+1}^{\infty}\!\bigg(Y_{k}-{{1}\over{k}}\bigg)\Bigg]\Bigg\}\Bigg)\cr
&\;= \prod_{k=m+1}^{\infty}{{\e^{-{\rm i}t/k}}\over{1-{{{\rm
i}t}\over{k}}}} = \exp\!\left\{\sum_{k=m+1}^{\infty}\!\left[-{{{\rm
i}t}\over{k}} - \log\!\left(1-{{{\rm
i}t}\over{k}}\right)\right]\right\},
\end{align}
which follows from the observation of Baum and Billingsley \cite{BB} that
$$
F_{m}(x) =
\p\Bigg(\sum_{k=m+1}^{\infty}\!\bigg(Y_{k}-{{1}\over{k}}\bigg) \leq
x\Bigg),\quad x\in\R,
$$
itself suggested by (\ref{Wgeo}), where the $Y_{m+1}, Y_{m+2}, \ldots$ are
independent random variables such that $Y_k$ has the exponential
distribution with mean $1/k$, and hence the characteristic function
$\E\big(\e^{{\rm i}tY_k}\big) = 1/\big(1-{{{\rm i}t}\over{k}}\big)$,
$t\in\R$. Setting
\begin{equation}\label{g_2.12}
A_{nm}(t) = \exp\!\left\{\sum_{k=m+1}^{n-1}\!\left[-{{{\rm
i}t}\over{k}} - \log\!\left(1-{{{\rm
i}t}\over{k}}\right)\right]\right\},\quad t\in\R,
\end{equation}
and noticing that for every $t\in\R$,
\begin{align*}
\varphi_{nm}(t) &=
\exp\!\left\{\sum_{k=m+1}^{n-1}\!\left[-{{{\rm i}t}\over{k}} -
\log\!\bigg(1+{{n}\over{k}}\big\{\e^{-{\rm
i}t/n}-1\big\}\bigg)\right]\right\}\cr &=
A_{nm}(t)\exp\!\left\{\!-\!\sum_{k=m+1}^{n-1}\!
\log\!{{1+{{n}\over{k}}(\e^{-{\rm i}t/n}-1)}\over{1-{{{\rm
i}t}\over{k}}}}\right\}\cr &=
A_{nm}(t)\exp\!\left\{\!-\!\sum_{k=m+1}^{n-1}\!
\log\!{{1+{{n}\over{k}}\!\left[{{-{\rm i}t}\over{n}} +
{{1}\over{2}}\left({{{\rm i}t}\over{n}}\right)^2 +
\sum_{l=3}^\infty\left({{-{\rm
i}t}\over{n}}\right)^l{{1}\over{l!}}\right]}\over {1-{{{\rm
i}t}\over{k}}}}\right\}\!,
\end{align*}
for $|t|<n$ we introduce the first sequence of intermediate
approximative functions
\begin{align*}
\psi_{nm,1}(t) &:=
A_{nm}(t)\exp\!\left\{\!-\!\sum_{k=m+1}^{n-1}\!
\log{{1+{{n}\over{k}}\!\left[{{-{\rm i}t}\over{n}} +
{{1}\over{2}}\left({{{\rm i}t}\over{n}}\right)^2\right]} \over
{1-{{{\rm i}t}\over{k}}}}\right\}\cr &\;=
A_{nm}(t)\exp\!\left\{\!-\!\sum_{k=m+1}^{n-1}\!
\log\!\left(1+{{({\rm i}t)^2}\over{2n}}\,{{1}\over{k-{\rm
i}t}}\right)\right\}\cr &\;=
A_{nm}(t)\exp\!\left\{\!-\!\sum_{k=m+1}^{n-1}\! \left[{{({\rm
i}t)^2}\over{2n}}{{1}\over{k-{\rm i}t}} +
\sum_{l=2}^\infty{{(-1)^{l-1}}\over{l}} \left({{({\rm
i}t)^2}\over{2n}}\,{{1}\over{k-{\rm
i}t}}\right)^{\!l}\,\right]\right\}\!,
\end{align*}
where the expansion of the logarithm is justified because the
inequality
\begin{equation}\label{g_2.13}
\left|{{({\rm i}t)^2}\over{2n}}\,{{1}\over{k-{\rm i}t}}\right| =
{{t^2}\over{2n}}\,{1\over{\sqrt{k^2+t^2}}} < 1
\end{equation}
holds whenever $|t|<n$. For all such $t$, the second intermediate
sequence is
\begin{align*}
\psi_{nm,2}(t) &:= A_{nm}(t)
\exp\!\left\{{{t^2}\over{2n}}\sum_{k=m+1}^{n-1}{{1}\over{k-{\rm
i}t}}\right\}\cr &\;=
A_{nm}(t)\!\left\{1+{{t^2}\over{2n}}\sum_{k=m+1}^{n-1}{{1}\over{k-{\rm
i}t}} + \sum_{l=2}^\infty{{1}\over{l!}}\left({{t^2}\over{2n}}
\sum_{k=m+1}^{n-1}{{1}\over{k-{\rm i}t}}\right)^{\!l}\right\}\!,
\end{align*}
while the third and the fourth are
$$
\psi_{nm,3}(t) := A_{nm}(t)\!
\left\{1+{{t^2}\over{2n}}\sum_{k=m+1}^{n-1}{{1}\over{k-{\rm
i}t}}\right\}
$$
and
\begin{align*}
\psi_{nm,4}(t) &:=
\exp\!\left\{\sum_{k=m+1}^{\infty}\!\left[-{{{\rm i}t}\over{k}} -
\log\!\left(1-{{{\rm i}t}\over{k}}\right)\right]\right\}\!
\left\{1+{{t^2}\over{2n}}\sum_{k=m+1}^{n-1}\!{{1}\over{k-{\rm
i}t}}\right\}\cr &\;=
\exp\!\left\{\sum_{k=m+1}^{\infty}\!\left[-{{{\rm i}t}\over{k}} -
\log\!\left(1-{{{\rm
i}t}\over{k}}\right)\right]\right\}\!\left\{1-{{1}\over{2n}}
\sum_{k=m+1}^{n-1}\!{{({\rm i}t)^2}\over{k}}{{k}\over{k-{\rm
i}t}}\right\}\!,
\end{align*}
and we notice from (\ref{g_2.11}) that $\psi_{nm,4}(t) =
\varphi_{m}(t)+\psi_{nm}(t)$ for all $t\in(-n,n)$, where
\begin{equation}\label{g_2.14}
\psi_{nm}(t) := -{{1}\over{2n}}\sum_{k=m+1}^{n-1}
{{1}\over{k}}\,{{k}\over{k-{\rm i}t}}\,({\rm
i}t)^2\varphi_{m}(t),\quad t\in\R.
\end{equation}

Here $k/(k-{\rm i}t)$ is the characteristic function of the
exponential distribution with mean $1/k$, so that
$$
{{k}\over{k-{\rm i}t}} = \int_{-\infty}^{\infty}\!\e^{{\rm
i}tx}h_{k}(x)\,\d x, \quad t\in\R,
$$
where $h_k(\cdot)$ is the density function in (\ref{g_1.4}). Also, since by
(\ref{g_2.40}) below the function $t\mapsto t^j\varphi_m(t)$ is integrable
on $\R$ for every $j\in\N$, we can differentiate the density
inversion formula twice to obtain $f_{m}^{\prime\prime}(\cdot)$ of
(\ref{g_1.3}) as the inverse Fourier transform
$$
f_{m}^{\prime\prime}(x) = {1\over 2\pi}
\int_{-\infty}^{\infty}\!\e^{-{\rm i}tx} ({\rm i}t)^2
\varphi_{m}(t)\,\d t, \quad x\in\R.
$$
Since $f_{m}^{\prime\prime}(\cdot)$ is also integrable, as
established at (2.2), this can be inverted to get
$$
({\rm i}t)^2\varphi_{m}(t) = \int_{-\infty}^{\infty}\!\e^{{\rm
i}tx}f_{m}^{\prime\prime}(x)\,\d x, \quad t\in\R.
$$
The two Fourier transforms then combine to give
$$
{{k}\over{k-{\rm i}t}}\,({\rm i}t)^2\varphi_{m}(t) =
\int_{-\infty}^{\infty}\! \e^{{\rm
i}tx}\,[f_{m}^{\prime\prime}\!\star h_{k}](x)\,\d x, \quad t\in\R,
$$
for the integrable convolution in (\ref{g_2.3}). Therefore, by (\ref{g_2.4}) we
recognize (\ref{g_2.14}) as
\begin{align*}
\psi_{nm}(t) &= \int_{-\infty}^{\infty}\!\e^{{\rm i}tx}
\Bigg(\!- {{1}\over{2n}}\sum_{k=m+1}^{n-1}{{1}\over{k}}\,
[f_{m}^{\prime\prime}\!\star h_{k}](x)\Bigg)\d x\cr &=
\int_{-\infty}^{\infty}\!\e^{{\rm i}tx}\,G_{n,m}^{\,\prime}(x)\,\d x
= \int_{-\infty}^{\infty}\!\e^{{\rm i}tx}\,\d G_{n,m}(x), \quad
t\in\R,
\end{align*}
for the integrable function $G_{n,m}^{\,\prime}(\cdot)$ for which
$\lim_{|x|\to\infty}G_{n,m}^{\,\prime}(x) = 0$, so that the function
$G_{n,m}(\cdot)$ is of bounded variation on the whole line $\R$.

Then the deviation $\Delta_n :=
\sup_{x\in\R}\big|F_{nm}(x)-[F_{m}(x)+G_{nm}(x)]\big|$ in (\ref{g_1.7}) may
be estimated through Esseen's inequality (see Section 2.2.), which we use
in the form
$$
\Delta_n \le {{b}\over{2\pi}}\int_{-cn}^{cn}\!
\left|{{\varphi_{nm}(t) - [\varphi_{m}(t) +
\psi_{nm}(t)]}\over{t}}\right|\d t +
c_{b}\,{{\sup_{x\in\R}\big|f_{m}(x) +
G_{n,m}^{\,\prime}(x)\big|}\over{cn}},
$$
where $b>1$ is arbitrary and $c_b>0$ is a constant depending only on
$b$ and, due to the restriction of the arguments $t$ of the
intermediate functions $\psi_{nm,j}(t)$, $j=1,2,3$, to $(-n,n)$, the
constant $c$ is taken from the interval $(0,1)$. Since
$\max_{x\in\R}f_{m}(x) = ((m+1)/\e)^{m+1}/m!$ by (\ref{g_2.6}), we see by
(\ref{g_2.9}) that the second term here is $O(1/n)$. Thus the proof of (\ref{g_1.7})
reduces to demonstrating that the same holds for the first term as
well. This will be split in four parts according to (\ref{g_2.10}): we have
$$
\int_{-cn}^{cn}
\left|{{\varphi_{nm}(t)-[\varphi_{m}(t)+\psi_{nm}(t)]}\over{t}}\right|\d
t \le \sum_{j=1}^{4}R_{nm,j}
$$
and, introducing the set $H_{nc} = [-cn,-1)\cup(1,cn]$, it suffices
to show that
\begin{align}\label{g_2.15}
R_{nm,j} &= \int_{-cn}^{cn}
\left|{{\rho_{nm,j}(t)}\over{t}}\right|\d t =
\int_{-1}^{1}\left|{{\rho_{nm,j}(t)}\over{t}}\right|\d t +
\int_{H_{nc}}\!\left|{{\rho_{nm,j}(t)}\over{t}}\right|\d t\cr &=:
I_{nm}^{j,1} + I_{nm}^{j,2} = O\!\bigg({1\over n}\bigg)
\quad\hbox{for each}\quad j=1,2,3,4.
\end{align}
We fix $c\in(0,1)$, let $n > 1/c$ and, unless otherwise stated,
assume in all formulae containing the variable $t$ throughout that
$t\in [-cn,cn]$.

\bigskip

{\cal The case of $R_{nm,1}$.} By (\ref{g_2.10}) and the definitions
between (\ref{g_2.12}) and (\ref{g_2.14}),
\begin{align*}
|\rho_{nm,1}(t)| &=
\left|\varphi_{nm}(t)-\psi_{nm,1}(t)\right|\cr &=
\Bigg|\exp\!\left\{\sum_{k=m+1}^{n-1}\!\left[-{{{\rm i}t}\over{k}} -
\log\!\left(1-{{{\rm i}t}\over{k}}\right)\right] -
\sum_{k=m+1}^{n-1}\! \log{{1+{{n}\over{k}}(\e^{-{\rm i}t/n}-1)}
\over{1-{{{\rm i}t}\over{k}}}}\right\}\cr
&\;\;\;\;\;-\exp\!\left\{\sum_{k=m+1}^{n-1}\!\left[-{{{\rm
i}t}\over{k}} - \log\!\left(1-{{{\rm i}t}\over{k}}\right)\right] -
\sum_{k=m+1}^{n-1}\log{{1-{{{\rm i}t}\over{k}}+{{({\rm
i}t)^2}\over{2nk}}}\over {1-{{{\rm i}t}\over{k}}}}\right\}\Bigg|.
\end{align*}
The inequality
\begin{equation}\label{g_2.16}
\left|\e^z-\e^w\right|\leq{{1}\over{2}}\,\big\{|\e^z|+|\e^w|\big\}|z-w|,\quad z,w\in\C,
\end{equation}
where $\C$ denotes the complex plane, yields
\begin{equation}\label{g_2.17}
|\rho_{nm,1}(t)|\leq
{{1}\over{2}}\,\big\{|\varphi_{nm}(t)|+|\psi_{nm,1}(t)|\big\}\,\delta_{nm}^{[1]}(t),
\end{equation}
where
$$
\delta_{nm}^{[1]}(t) = \left|\sum_{k=m+1}^{n-1}
\log{{1+{{n}\over{k}}\left(\e^{-{\rm i}t/n}-1\right)} \over{1-{{{\rm
i}t}\over{k}}+{{({\rm i}t)^2}\over{2nk}}}}\right| =:
\left|\sum_{k=m+1}^{n-1}\log z_{nk}(t)\right|.
$$
We give upper bounds for each of the functions on the right-hand
side of (\ref{g_2.17}).

As usual, let $\re z$ denote the real part of $z\in\C$. Clearly,
\begin{equation}\label{g_2.18}
\delta_{nm}^{[1]}(t) \leq \sum_{k=m+1}^{n-1}\big|\log z_{nk}(t)\big|
= \sum_{k=m+1}^{n-1}\bigg|\log {{1}\over{z_{nk}(t)}}\bigg|.
\end{equation}
First we show that $\re z_{nk}(t) \ge 1/2$, so that $1/z_{nk}(t)$ is
an inner point of the circle of center $(1,0)$ and radius 1 in $\C$,
and hence its logarithm can be expanded about the point $(1,0)$. We
have
\begin{align*}
&\re z_{nk}(t)-{{1}\over{2}} = \re
{{\left(1+{{n}\over{k}}\cos{{t}\over{n}}-{{n}\over{k}}\right) - {\rm
i}{{n}\over{k}}\sin{{t}\over{n}}} \over
{\left(1-{{t^2}\over{2nk}}\right) - {\rm
i}{{t}\over{k}}}}-{{1}\over{2}}\cr \noalign{\vskip3pt}
&={{\left(1+{{n}\over{k}}\cos{{t}\over{n}}-{{n}\over{k}}\right)
\big(1-{{t^2}\over{2nk}}\big)+{{nt}\over{k^2}}\sin{{t}\over{n}}}\over
{\left(1-{{t^2}\over{2nk}}\right)^{\!2} +
{{t^2}\over{k^2}}}}-{{1}\over{2}}\cr \noalign{\vskip3pt}
&={{\big(2{{n}\over{k}}-{{t^2}\over{k^2}}\big)\cos{{t}\over{n}}+2{{tn}\over
{k^2}}\sin{{t}\over{n}}-2{{n}\over{k}}-{{t^4}\over{4n^2k^2}}+1}\over
{2\left[\left(1-{{t^2}\over{2nk}}\right)^2+{{t^2}\over{k^2}}\right]}}\cr
\noalign{\vskip3pt}
&={{\big({{tn}\over{k^2}}\sin{{t}\over{n}}-{{t^2}\over{k^2}}\cos{{t}\over{n}}\big)
+ \big({{tn}\over{2k^2}}\sin{{t}\over{n}}-{{t^4}\over{4n^2k^2}}\big)
+
\left({{tn}\over{2k^2}}\sin{{t}\over{n}}+2{{n}\over{k}}\cos{{t}\over{n}}
-2{{n}\over{k}}+1\right)}\over
{2\left[\left(1-{{t^2}\over{2nk}}\right)^2+{{t^2}\over{k^2}}\right]}}
\end{align*}
This is an even function of $t$, so we can assume that $t\geq 0$.
The denominator is obviously positive, and we are going to show that
each of the three terms in the numerator is non-negative. Beginning
with the first term, we see that
\begin{align*}
{{tn}\over{k^2}}&\sin{{t}\over{n}}-{{t^2}\over{k^2}}\cos{{t}\over{n}}
={{t}\over{k^2}}\left\{n\sum_{l=0}^{\infty}{{(-1)^l}\over{(2l+1)!}}
\left({{t}\over{n}}\right)^{\!2l+1}\!
-t\sum_{l=0}^{\infty}{{(-1)^l}\over{(2l)!}}
\left({{t}\over{n}}\right)^{\!2l}\right\}\cr
&={{t}\over{k^2}}\sum_{l=1}^{\infty}(-1)^l
\left({{1}\over{(2l+1)!}}-{{1}\over{(2l)!}}\right){{t^{2l+1}}\over{n^{2l}}}\cr
&={{t}\over{k^2}}\sum_{l\,{\rm is\, odd},\,l=1}^{\infty}
\left[\left({{1}\over{(2l)!}}-{{1}\over{(2l+1)!}}\right)
{{t^{2l+1}}\over{n^{2l}}}
-\left({{1}\over{(2l+2)!}}-{{1}\over{(2l+3)!}}\right)
{{t^{2l+3}}\over{n^{2l+2}}}\right]\cr
&\geq{{t}\over{k^2}}\sum_{l\,{\rm is\, odd},\,l=1}^{\infty}
\left({{1}\over{(2l)!}}-{{1}\over{(2l+1)!}} -
{{1}\over{(2l+2)!}}+{{1}\over{(2l+3)!}}\right){{t^{2l+1}}\over{n^{2l}}}\cr
&={{t}\over{k^2}}\sum_{l\,{\rm is\, odd},\,l=1}^{\infty}
{{(2l+3)(4l^2+4l-1)+1}\over{(2l+3)!}}\,{{t^{2l+1}}\over{n^{2l}}}
\geq 0,
\end{align*}
where the inequality is by ${{t}\over{n}}< c < 1$. Concerning the
second term, we note that $\sin x \geq {{x^2}\over{2}}$ if $0 \leq x
< 1$. Therefore, since ${{t}\over{n}}< 1$, we have
$$
{{tn}\over{2k^2}}\sin{{t}\over{n}}-{{t^4}\over{4n^2k^2}} \geq
{{tn}\over{2k^2}}{{t^2}\over{2n^2}}-{{t^4}\over{4n^2k^2}} =
{{t^2}\over{4k^2}}\left({{t}\over{n}}-\left({{t}\over{n}}\right)^{\!2}\right)
\geq 0.
$$
Finally, the third term can be settled using $t < n$ and the
inequalities $\sin x \geq x-{{x^3}\over{6}}$ and $\cos x \geq
1-{{x^2}\over{2}}$, both valid if $0 \leq x \leq 1$. Indeed,
\begin{align*}
{{tn}\over{2k^2}}\sin{{t}\over{n}}+2{{n}\over{k}}\cos{{t}\over{n}}-2{{n}\over{k}}+1
&\geq {{tn}\over{2k^2}}\left({{t}\over{n}} -
{{t^3}\over{6n^3}}\right) +
2{{n}\over{k}}\left(1-{{t^2}\over{2n^2}}\right)-2{{n}\over{k}}+1\cr
& ={{t^2}\over{2k^2}}-{{t^4}\over{12k^2n^2}}-{{t^2}\over{kn}}+1\cr
&\geq{{t^2}\over{2k^2}}-{{t^2}\over{12k^2}}-{{t}\over{k}}+1
={{5t^2-12kt+12k^2}\over{12k^2}}> 0.
\end{align*}

Returning now to (\ref{g_2.18}), we can expand the logarithm:
\begin{equation}\label{g_2.19}
\delta_{nm}^{[1]}(t) \le \!\sum_{k=m+1}^{n-1}\!\Bigg|
\sum_{j=1}^{\infty}
{{(-1)^{j-1}}\over{j}}\left[{{1}\over{z_{nk}(t)}}-1\right]^j\Bigg|
\le \!\sum_{k=m+1}^{n-1}\sum_{j=1}^{\infty}
{{1}\over{j}}\left|{{z_{nk}(t)-1}\over{z_{nk}(t)}}\right|^j.
\end{equation}
The inequalities
\begin{equation}\label{g_2.20}
\Bigg|\e^{{\rm i}u}-\sum_{j=0}^{2}{{({\rm i}u)^j}\over{j!}}\Bigg|
\leq {{|u|^{3}}\over{3!}},\;\;\,u\in\R,\,\quad\hbox{and}\quad 1-\cos
x\geq{{4}\over{\pi^2}}\,x^2,\;\;\,0\leq x\leq{{\pi}\over{2}},
\end{equation}
give the bound
\begin{align*}
\left|{{z_{nk}(t)-1}\over{z_{nk}(t)}}\right| &=
{{n}\over{k}}\left|{{\e^{-{\rm i}t/n}-1+{{{\rm i}t}\over{n}}-{{({\rm
i}t)^2} \over{2n^2}}}\over{1+{{n}\over{k}}\left(\e^{-{\rm
i}t/n}-1\right)}}\right| \le
{{n}\over{k}}{{{{|t|^3}\over{6n^3}}}\over
{\left|\left(1-{{n}\over{k}}+{{n}\over{k}}\cos{{{t}\over{n}}}\right)
- {\rm i}{{n}\over{k}}\sin{{{t}\over{n}}}\right|}}\cr &=
{{|t|^3}\over{6kn^2}}{{1}\over{\sqrt{1+2{{n}\over{k}}
\left({{n}\over{k}}-1\right)\big(1-\cos{{|t|}\over{n}}\big)}}} \le
{{|t|^3}\over{6kn^2}}{{1}\over{\sqrt{1+{{8}\over{\pi^2}}{{n}\over{k}}
\left({{n}\over{k}}-1\right){{t^2}\over{n^2}}}}}\cr &=
{{|t|^3}\over{6n^2}}{{1}\over{\sqrt{k^2-{{8t^2}\over{n\pi^2}}k +
{{8t^2}\over{\pi^2}}}}} = {{{|t|^3}\over{6n^2}}}{1\over \sqrt{k^2 +
{8t^2\over \pi^2}\left(1-{k\over n}\right)}},
\end{align*}
which, for any $k\in\{m+1, m+2,\ldots,n-1\}$ and $|t|<n$, is not
greater than
$$
{{{|t|^3}\over{6n^2}}}{1\over {\sqrt{k^2 + {8t^2\over
\pi^2}\left(1-{n-1\over n}\right)}}} \leq
{{{|t|^3}\over{6n^2}}}{1\over
{\sqrt{{{8t^2}\over{\pi^2}}{{1}\over{2}}}}} =
{{t^2}\over{n^2}}{{\pi}\over{12}} < {1\over 3}.
$$
Substituting these bounds into (\ref{g_2.19}), we easily obtain
\begin{align*}
\delta_{nm}^{[1]}(t) &\leq {|t|^3\over 6n^2}
\sum_{k=m+1}^{n-1}{1\over \sqrt{k^2 + {8t^2\over
\pi^2}\left(1-{k\over n}\right)}} \left[1 +
\sum_{j=2}^{\infty}{1\over j}\left({{{{|t|^3}\over{6n^2}}}\over
\sqrt{k^2 + {8t^2\over \pi^2}\left(1-{k\over
n}\right)}}\right)^{j-1}\right]\cr &\leq
{{|t|^3}\over{6n^2}}\sum_{k=m+1}^{n-1}{{1}\over{\sqrt{k^2}}}
\left[1+{{1}\over{2}}\sum_{j=2}^{\infty}\left({{1}\over{3}}\right)^{\!j-1}\right]
= {5\over 24}\,{{|t|^3}\over{n^2}}\sum_{k=m+1}^{n-1}{{1}\over{k}},
\end{align*}
so that by (\ref{g_2.7}),
\begin{equation}\label{g_2.21}
\delta_{nm}^{[1]}(t)\leq {{|t|^3 \log n}\over{n^2}}.
\end{equation}

Next we consider $|\varphi_{nm}(t)|$ in (\ref{g_2.17}). Since $\re\log z =
\log |z|$, $z\in\C\setminus \{0\}$,
\begin{align*}
\big|\varphi_{nm}(t)\big| &=
\left|\exp\!\left\{\sum_{k=m+1}^{n-1}\!\left[-{{{\rm i}t}\over{k}} -
\log\!\bigg(1+{{n}\over{k}}\Big(\e^{-{\rm i}t/n}-1\Big)\bigg)
\right]\right\}\right|\cr &= \exp\!\left\{-\sum_{k=m+1}^{n}\!
\log\bigg|1+{{n}\over{k}}\Big(\e^{-{\rm
i}t/n}-1\Big)\bigg|\right\}\cr &=
\exp\!\left\{-{{1}\over{2}}\sum_{k=m+1}^{n}\!
\log\!\left[1+2n{{n-k}\over{k^2}}\left(1-\cos{{|t|}\over{n}}\right)\right]\right\}.
\end{align*}
As $0\leq|t|/n\leq c<1<\pi/2$, it follows from the second inequality
in (\ref{g_2.20}) that
\begin{equation}\label{g_2.22}
\big|\varphi_{nm}(t)\big| \leq
\exp\!\left\{-{{1}\over{2}}\sum_{k=m+1}^{n}\!
\log\!\left(1+{{8t^2}\over{\pi^2}}{{n-k}\over{nk^2}}\right)\right\}.
\end{equation}
The terms in the sum are positive, so for $n > 2(m+1)$ the exponent
is bounded by
$$
-{{1}\over{2}}\sum_{k=m+1}^{\lfloor n/2 \rfloor}\!
\log\!\left[1+{{8t^2}\over{\pi^2}}{{n-{{n}\over{2}}}\over{nk^2}}\right]
\leq -{{1}\over{2}}\int_{m+1}^{n/2}\!
\log\!\left[1+{{4t^2}\over{\pi^2y^2}}\right]\!\d y =: I_{nm}^{1,0},
$$
where $\lfloor\cdot\rfloor$ denotes integer part, because the terms
in the last sum decrease as $k$ increases. Substituting in the
integral $x = 4t^2/(\pi^2y^2)$, so that $y = 2t/(\pi\sqrt{x})$,
\begin{align*}
I_{nm}^{1,0} &= -{{|t|}\over{2\pi}}
\int_{\left({{4t}\over{n\pi}}\right)^2}^{\left({{2t}\over{(m+1)\pi}}\right)^2}
{{\log(1+x)}\over{x^{3/2}}}\,\d x\cr &= -{{|t|}\over{2\pi}}
\left[{{-2\log(1+x)}\over{\sqrt{x}}}
\right]_{x=\left({{4t}\over{n\pi}}\right)^2}^{x=\big({{2t}\over{(m+1)\pi}}\big)^2}
-{{|t|}\over{\pi}}
\int_{\left({{4t}\over{n\pi}}\right)^2}^{\left({{2t}\over{(m+1)\pi}}\right)^2}\!
{{1}\over{\sqrt{x}\,(1+x)}}\,\d x\cr &=
{{m+1}\over{2}}\log\!\left(1+{{4t^2}\over{(m+1)^2\pi^2}}\right)
-{{n}\over{4}}\log\!\left(1+{{16t^2}\over{n^2\pi^2}}\right)\cr
&\;\;\;\;\,-{{2|t|}\over{\pi}}
\left[\arctan{{2|t|}\over{(m+1)\pi}}-\arctan{{4|t|}\over{n\pi}}\right]
\end{align*}
Summarizing, at this stage we have
$$
\big|\varphi_{nm}(t)\big| \leq |t|^{m+1}
{{\big({{1}\over{t^2}}+{{4}\over{(m+1)^2\pi^2}}\big)^{\!{{m+1}\over{2}}}}
\over {{\left(1+{{16t^2}\over{n^2\pi^2}}\right)}^{{{n}\over{4}}}}}\,
\exp\!\left\{\!-{{2|t|}\over{\pi}}\,
\arctan{{{{2|t|}\over{\pi}}\left({{1}\over{m+1}}-{{2}\over{n}}\right)}\over
{1+{{8t^2}\over{(m+1)n\pi^2}}}}\right\}.
$$
Using $|t|<n$ and assuming $|t|\ge 1$, we can further simplify this
to obtain
$$
|\varphi_{nm}(t)|\leq
|t|^{m+1}{{\big(1+{{4}\over{(m+1)^2\pi^2}}\big)^{\!{{m+1}\over{2}}}}\over
{{\left(1+{{16}\over{n^2\pi^2}}\right)}^{{{n}\over{4}}}}}\,
\exp\!\left\{\!-{{2|t|}\over{\pi}}\,
\arctan{{{{2|t|}\over{\pi}}\big({{1}\over{m+1}}-{{2}\over{n}}\big)}\over
{1+{{8|t|}\over{(m+1)\pi^2}}}}\right\}.
$$
If $n>2(m+1)$, then the $\arctan$ expression in the exponent is a
monotone increasing function of $|t|$ because the numerator of the
derivative $\d \{as/(1+bs)\}/\d s$ is equal to $a$ for any real
constants $a$ and $b$. Hence for $n>2(m+1)$, which was already
assumed above anyway to get to the previous bound,
\begin{equation}\label{g_2.23}
\big|\varphi_{nm}(t)\big|\leq
d_{nm}\,|t|^{m+1}\,\e^{-r_{nm}|t|},\quad 1\le |t|<n,
\end{equation}
where
$$
d_{nm} = {{\big(1+{{4}\over{(m+1)^2\pi^2}}\big)^{\!(m+1)/2}}\over
{{\left(1+{{16}\over{n^2\pi^2}}\right)}^{n/4}}} \to
\left(1+{{4}\over{(m+1)^2\pi^2}}\right)^{\!(m+1)/2} =: d_m
$$
and
$$
r_{nm} = {2\over \pi}\,
\arctan{{{{2}\over{\pi}}\big({{1}\over{m+1}}-{{2}\over{n}}\big)}\over
{1+{{8}\over{(m+1)\pi^2}}}} \to {2\over \pi}\, \arctan{2\pi\over
(m+1)\pi^2 + 8} =: r_m.
$$

Finally, for (\ref{g_2.17}), simplifying the second line of the definition
of $\psi_{nm,1}(t)$,
\begin{align*}
\big|\psi_{nm,1}(t)\big|&=
\left|\exp\!\left\{\sum_{k=m+1}^{n-1}\! \left[-{{{\rm i}t}\over{k}}
- \log\!\left(1-{{{\rm i}t}\over{k}}+{{({\rm i}t)^2}\over{2nk}}
\right)\right]\right\}\right|\cr &=
\exp\!\left\{-\sum_{k=m+1}^{n-1}\!\log\! \left|1-{{{\rm
i}t}\over{k}}+{{({\rm i}t)^2}\over{2nk}}\right|\right\}\cr &=
\exp\!\left\{-\sum_{k=m+1}^{n-1}\!
\log\sqrt{1-{{t^2}\over{nk}}+{{t^4}\over{4n^2k^2}}+{{t^2}\over{k^2}}}\,\right\}\cr
&\leq
\exp\!\left\{-\sum_{k=m+1}^{n}\!\log\sqrt{1+t^2{{n-k}\over{nk^2}}}\,\right\}\cr
&\leq \exp\!\left\{- {1\over 2}\sum_{k=m+1}^{n}\!
\log\!\left(1+{{8t^2}\over{\pi^2}}{{n-k}\over{nk^2}}\right)\right\}.
\end{align*}
The artificial factor $8/\pi^2 < 1$ was sneaked in just to get the
exact same upper bound as in (\ref{g_2.22}) for $|\psi_{nm,1}(t)|$, and
hence to conclude without any extra work that
\begin{equation}\label{g_2.24}
\big|\psi_{nm,1}(t)\big|\leq
d_{nm}\,|t|^{m+1}\,\e^{-r_{nm}|t|},\quad 1\le |t|<n,
\end{equation}
for $n>2(m+1)$, as in (\ref{g_2.23}), with the same $d_{nm}\to d_m$ and
$r_{nm}\to r_m > 0$.

Now, recalling the definition of $R_{nm,1} = I_{nm}^{1,1} +
I_{nm}^{1,2}$ in (\ref{g_2.15}), suppose that $n>\max\{2(m+1),1/c\}$. Since
$|\varphi_{nm}(t)|\leq 1$ and $|\psi_{nm,1}(t)|\leq1$, the
inequalities (\ref{g_2.17}) and (\ref{g_2.21}) yield $I_{nm}^{1,1} \le 2(\log
n)/(3n^2)$. For the other term (\ref{g_2.17}), (\ref{g_2.21}), (\ref{g_2.23}), (\ref{g_2.24}), and
the fact that the functions involved are even, imply that
$$
I_{nm}^{1,2} \leq {{2\log n}\over{n^2}}\,d_{nm}
\int_{1}^{\infty}\!t^{m+3}\,\e^{-r_{nm}t}\,\d t,
$$
so that the case $j=1$ in (\ref{g_2.15}) holds true; in fact, $R_{nm,1} =
O((\log n)/n^2)$.

\bigskip

{\cal The case of $R_{nm,2}$.} By (\ref{g_2.10}) and the formulae between
(\ref{g_2.12}) and (\ref{g_2.14}),
\begin{align*}
|\rho_{nm,2}(t)| &= |\psi_{nm,1}(t)-\psi_{nm,2}(t)|\cr
&=\Bigg|\exp\!\left\{\sum_{k=m+1}^{n-1}\!\left[-{{{\rm i}t}\over{k}}
- \log\!\left(1-{{{\rm i}t}\over{k}}\right)\right] -
\sum_{k=m+1}^{n-1}\!\log\!\left(1+{{({\rm i}t)^2}\over{2n}}\,
{{1}\over{k-{\rm i}t}}\right)\right\}\cr &\quad\;\;-
\exp\!\left\{\sum_{k=m+1}^{n-1}\!\left[-{{{\rm i}t}\over{k}} -
\log\!\left(1-{{{\rm i}t}\over{k}}\right)\right] +
{{t^2}\over{2n}}\sum_{k=m+1}^{n-1}{{1}\over{k-{\rm
i}t}}\right\}\Bigg|.
\end{align*}
The inequality (\ref{g_2.16}) now gives
\begin{equation}\label{g_2.25}
|\rho_{nm,2}(t)|\leq
{{1}\over{2}}\,\big\{|\psi_{nm,1}(t)|+|\psi_{nm,2}(t)|\big\}\,
\delta_{nm}^{[2]}(t),
\end{equation}
as an analogue of (\ref{g_2.17}), where
\begin{align*}
\delta_{nm}^{[2]}(t) &= \left|\sum_{k=m+1}^{n-1}\!
\log\!\left(1+{{({\rm i}t)^2}\over{2n}}\,{{1}\over{k-{\rm
i}t}}\right) -{{({\rm
i}t)^2}\over{2n}}\sum_{k=m+1}^{n-1}{{1}\over{k-{\rm i}t}}\right|\cr
&= \left|\sum_{k=m+1}^{n-1}\left[\sum_{j=1}^{\infty}\!
\left\{{{(-1)^{j-1}}\over{j}} \left({{({\rm i}t)^2}\over{2n(k-{\rm
i}t)}}\right)^j\right\} -{{({\rm i}t)^2}\over{2n}}{{1}\over{k-{\rm
i}t}}\right]\right|\cr &=\Bigg|\sum_{k=m+1}^{n-1}
\sum_{j=2}^{\infty}{{(-1)^{j-1}}\over{j}} \left({{({\rm
i}t)^2}\over{2n(k-{\rm i}t)}}\right)^j\Bigg|\cr &\leq
\sum_{k=m+1}^{n-1}\sum_{j=2}^{\infty}
{{1}\over{j}}\left|{{t^2}\over{2n(k-{\rm i}t)}}\right|^j
=\sum_{k=m+1}^{n-1}\sum_{j=2}^{\infty}
{{1}\over{j}}\left({{t^2}\over{2n\sqrt{k^2+t^2}}}\right)^j\cr
&\leq\sum_{k=m+1}^{n-1}{{1}\over{2}}
\left({{t^2}\over{2n\sqrt{k^2+t^2}}}\right)^{\!2}\,
\sum_{j=2}^{\infty}\left({{t^2}\over{2n\sqrt{k^2+t^2}}}\right)^{\!j-2}.
\end{align*}
Since $|t|<n$, we find that $t^2/\big(2n\sqrt{k^2+t^2}\,\big) <
n^2/\big(2n\sqrt{k^2+n^2}\,\big) < 1/2$, and so
$$
\delta_{nm}^{[2]}(t)\leq
{{1}\over{2}}\sum_{k=m+1}^{n-1}\!\!\bigg({{t^2}\over{2n\sqrt{k^2+t^2}}}\bigg)^{\!2}
\sum_{j=2}^{\infty}\!\left({{1}\over{2}}\right)^{\!j-2}\! =
{{t^4}\over{4n^2}}\sum_{k=m+1}^{n-1}{{1}\over{k^2+t^2}}
\leq{{t^4}\over{4n^2}}\sum_{k=1}^{\infty}{{1}\over{k^2}}.
$$
Therefore,
\begin{equation}\label{g_2.26}
\delta_{nm}^{[2]}(t)\leq {{\pi^2}\over{24}}\,{{t^4}\over{n^2}}.
\end{equation}

Next, recalling the notation for $A_{nm}(t)$ from (\ref{g_2.12}), we write
\begin{equation}\label{g_2.27}
\big|\psi_{nm,2}(t)\big| = \big|A_{nm}(t)\big|\,\big|B_{nm}(t)\big|
\;\;\,\hbox{with}\;\;\, B_{nm}(t) = \exp\!\left\{{{t^2}\over{2n}}
\sum_{k=m+1}^{n-1} {{1}\over{k-{\rm i}t}}\right\}\!,
\end{equation}
and estimate the two factors separately. First,
\begin{align*}
\big|A_{nm}(t)\big| &= \exp\!\left\{-\sum_{k=m+1}^{n-1}\!
\log\left|1-{{{\rm i}t}\over{k}}\right|\right\} =
\exp\!\left\{-{{1}\over{2}}\sum_{k=m+1}^{n-1}\!
\log\!\left(1+{{t^2}\over{k^2}}\right)\right\}\cr
&\leq\exp\!\left\{-{{1}\over{2}}\int_{m+1}^{n}\!
\log\!\left(1+{{t^2}\over{y^2}}\right)\d y\right\},
\end{align*}
where the inequality holds again because the terms of the last sum
decrease as $k$ increases. Substituting $x = t^2/y^2$, so that $y =
|t|/\sqrt{x}$, we obtain
\begin{align*}
\big|A_{nm}(t)\big| &\leq \exp\!\left\{\!-{{|t|}\over{4}}
\int_{\left({{t}\over{n}}\right)^2}^{\left({{t}\over{m+1}}\right)^2}
{{\log(1+x)}\over{x^{3/2}}}\,\d x\right\}\cr &=
\exp\!\left\{\!-{{|t|}\over{4}}\left[{{-2\log{(1+x)}}
\over{\sqrt{x}}}\right]_{x=\left({{t}\over{n}}\right)^2}
^{x=\left({{t}\over{m+1}}\right)^2}\! -
{{|t|}\over{2}}\int_{\left({{t}\over{n}}\right)^2}^{\left({{t}\over{m+1}}\right)^2}
{{1}\over{\sqrt{x}\,(1+x)}}\,\d x\right\}\cr
&={{\left(1+{{t^2}\over{(m+1)^2}}\right)^{\!{{m+1}\over{2}}}}\over
{\left(1+{{t^2}\over{n^2}}\right)^{{{n}\over{2}}}}}
\exp\!\left\{\!-|t|\arctan{{{|t|\left({{1}\over{m+1}}-{{1}\over{n}}\right)}\over
{1+{{t^2}\over{(m+1)n}}}}}\right\},
\end{align*}
whence
$$\big|A_{nm}(t)\big|\leq
|t|^{m+1}{{\big({{1}\over{t^2}}+{{1}\over{(m+1)^2}}\big)^{\!(m+1)/2}}
\over \left(1+{{t^2}\over{n^2}}\right)^{n/2}}\,
\exp\!\left\{\!-|t|\arctan{{{|t|(n-m-1)}\over{t^2+(m+1)n}}}\right\}.
$$
If $|t|\ge 1$, then
\begin{equation}\label{g_2.28}
\big|A_{nm}(t)\big|\leq a_{nm}\,|t|^{m+1}\,
\exp\!\left\{\!-|t|\arctan{{{|t|(n-m-1)}\over{t^2+(m+1)n}}}\right\},
\end{equation}
where
$$
a_{nm} = {{\big(1 + {{1}\over{(m+1)^2}}\big)^{\!(m+1)/2}} \over
\left(1+{1\over{n^2}}\right)^{n/2}} \to
{\left(1+{{1}\over{(m+1)^2}}\right)^{\!{{m+1}\over{2}}}} =: a_m,
$$
and the $\arctan$ expression in the exponent, as a function of
$|t|$, is monotone increasing on $(0,\sqrt{(m+1)n})$ and monotone
decreasing on $(\sqrt{(m+1)n}\,,n)$ because
$$
{\d \over \d s}\!\left({{s(n-m-1)}\over{s^2+(m+1)n}}\right) =
{{(n-m-1)\left(n(m+1)-s^2\right)}\over{\left(s^2+(m+1)n\right)^2}}.
$$
Then it follows that
\begin{equation}\label{g_2.29}
{\big|A_{nm}(t)\big|\over a_{nm}} \leq
\left\{
  \begin{array}{ll}
    |t|^{m+1}\exp\!\left\{-|t|\arctan{{{n-m-1}\over{1+(m+1)n}}}\right\}, & \hbox{if $t\in[1,\sqrt{(m+1)n}\,)$;} \\
    |t|^{m+1}\exp\!\left\{-|t|\arctan{{{n-m-1}\over{n+m+1}}}\right\}, & \hbox{if $t\in[\sqrt{(m+1)n}\,,n)$.}
  \end{array}
\right.
\end{equation}

Next we deal with the other factor in (\ref{g_2.27}), for which
\begin{align}\label{g_2.30}
\big|B_{nm}(t)\big| &=
\exp\!\left\{\re\!\!\left({{t^2}\over{2n}}
\sum_{k=m+1}^{n-1}{{1}\over{k-{\rm i}t}}\right)\right\} =
\exp\!\left\{{{t^2}\over{2n}}\sum_{k=m+1}^{n-1}{{k}\over{k^2+t^2}}\right\}\cr
&\le \exp\!\left\{{3\over 2}\,t^2\,{{\log n}\over{n}}\right\},
\quad\hbox{whenever}\quad |t|<n,
\end{align}
by (\ref{g_2.7}). This will be good enough for small $|t|$, and for moderate
$|t|$ it may be written as
\begin{equation}\label{g_2.31}
\big|B_{nm}(t)\big|\leq \exp\!\left\{{3{\sqrt{m+1}}\over{2}}\,|t|\,
{{\log n}\over{\sqrt{n}}}\right\},\quad\hbox{if}\quad |t|\le
\sqrt{(m+1)n}\,.
\end{equation}
To obtain yet a third bound, useful for large $|t|$, note that $\d
\{y/(y^2+t^2)\}/\d y = (t^2-y^2/(y^2+t^2)^2$. The behavior of this
derivative shows that for $|t|\leq m+1$ the terms of the sum in the
formula of $|B_{nm}(t)|$ decrease as $k$ increases, while for $m+1 <
|t|\leq n$ the terms increase until $k$ reaches $|t|$ and decrease
afterward. Hence
\begin{align*}
\sum_{k=m+1}^{n-1}{{k}\over{k^2+t^2}}&\leq
\int_{m}^{n}{{y}\over{y^2+t^2}}\,\d y +
{{\lfloor|t|\rfloor}\over{\lfloor|t|\rfloor^2+t^2}}
=\log\sqrt{{{n^2+t^2}\over{m^2+t^2}}} +
{{\lfloor|t|\rfloor}\over{\lfloor|t|\rfloor^2+t^2}}\cr & <
\log\sqrt{{{2n^2}\over{t^2}}}+{{|t|}\over{{{t^2}\over{2}}+t^2}} <
{{\sqrt{2}\,n}\over{|t|}}+{{2}\over{3|t|}} =
{{3\sqrt{2}\,n+2}\over{3|t|}},
\end{align*}
where, by elementary considerations concerning integer parts, the
second inequality holds provided $\sqrt{2}/(\sqrt{2}-1)\le |t|< n$.
Thus, substituting this bound into (\ref{g_2.30}),
\begin{equation}\label{g_2.32}
\big|B_{nm}(t)\big| \leq
\exp\!\left\{|t|\,{{3\sqrt{2}n+2}\over{6n}}\right\},\quad\hbox{if}\quad
{\sqrt{2}\over \sqrt{2}-1}\le |t|< n.
\end{equation}

Introduce now the sets $S_{nm} =
\big[-\sqrt{(m+1)n}\,,-1\big)\cup\big(1,\sqrt{(m+1)n}\,\big]$ and
$T_{nm} =
\big[-cn,-\sqrt{(m+1)n}\,\big)\cup\big(\sqrt{(m+1)n}\,,cn\big]$, so
that $H_{nc} = [-cn,-1)\cup(1,cn] = S_{nm}\cup T_{nm}$ in the second
term of $R_{nm,2} = I_{nm}^{2,1} + I_{nm}^{2,2}$ in (\ref{g_2.15}). Since
$|\psi_{nm,1}(t)|\leq 1$ and $|A_{nm}(t)|\leq 1$, (\ref{g_2.25}), (\ref{g_2.26}) and
(\ref{g_2.27}), combined with (\ref{g_2.30}), give
$$
I_{nm}^{2,1} \leq {\pi^2\over 48}\,{1\over n^2} \int_{-1}^{1}\!|t|^3
\bigg[1 + \exp\!\bigg\{{3\over 2}\,t^2\,{{\log
n}\over{n}}\bigg\}\bigg]\d t \leq {\pi^2\over 96}\, {1 +
\exp\!\left\{{3\over 2}\,{{\log n}\over{n}}\right\}\over n^2}
=O\!\bigg({{1}\over{n^2}}\bigg),
$$
while by (\ref{g_2.25}) and (\ref{g_2.26}) only, dropping the factor $\pi^2/48 \le
1$ for simplicity,
\begin{align*}
I_{nm}^{2,2} &\le {1\over
n^2}\int_{H_{nc}}\!\!|t|^3|\psi_{nm,1}(t)|\,\d t + {1\over
n^2}\int_{S_{nm}}\!\!|t|^3|\psi_{nm,2}(t)|\,\d t + {1\over
n^2}\int_{T_{nm}}\!\!|t|^3|\psi_{nm,2}(t)|\,\d t\cr &=:
I_{nm}^{2,2,1} + I_{nm}^{2,2,2} + I_{nm}^{2,2,3}.
\end{align*}
Here $I_{nm}^{2,2,1} = O(1/n^2)$ by (\ref{g_2.24}). Also, by (\ref{g_2.27}), (\ref{g_2.29})
and (\ref{g_2.31}),
$$
I_{nm}^{2,2,2}\le
{{2}\over{n^2}}\int_{1}^{\sqrt{(m+1)n}}t^3|\psi_{nm,2}(t)|\,\d t
\leq {{2a_{nm}}\over{n^2}}
\int_{1}^{\sqrt{(m+1)n}}t^{m+4}\,\e^{-b_{nm}t}\,\d t,
$$
where $a_{nm}\to a_m$ and
$$
b_{nm} := \arctan{{n-m-1}\over{1+(m+1)n}}-{{\sqrt{m+1}}\over{2}}\,
{{\log n}\over{\sqrt{n}}} \to \arctan{{1}\over{m+1}}>0,
$$
so that $I_{nm}^{2,2,2} = O(1/n^2)$ as well. Finally, by (\ref{g_2.27}),
(\ref{g_2.29}) and (\ref{g_2.32}),
$$
I_{nm}^{2,2,3} \le {{2}\over{n^2}} \int_{\sqrt{(m+1)n}}^{cn}
t^3|\psi_{nm,2}(t)|\,\d t \leq {{2a_{nm}}\over{n^2}}
\int_{\sqrt{(m+1)n}}^{n} t^{m+4}\,\e^{-c_{nm}t}\,\d t
$$
with the same $a_{nm}$, where
$$
c_{nm} := \arctan{{n-m-1}\over{n+m+1}}-{{3\sqrt{2}n+2}\over{6n}} \to
{{\pi}\over{4}}-{{\sqrt{2}}\over{2}}>0.078\,.
$$
Thus $I_{nm}^{2,2,3} = o(1/n^2)$. So, the case $j=2$ of (\ref{g_2.15})
holds; in fact $R_{nm,2} = O(1/n^2)$.

\bigskip
{\cal The case of $R_{nm,3}$.} Recalling (\ref{g_2.10}) and the formulae
after (\ref{g_2.12}), we have $|\rho_{nm,3}(t)| =
\left|\psi_{nm,2}(t)-\psi_{nm,3}(t)\right| =
|A_{nm}(t)|\,|C_{nm}(t)|$, where $A_{nm}(t)$ is given in (\ref{g_2.12}), and
already occurs also in (\ref{g_2.27}), and
$$
C_{nm}(t) =
\sum_{l=2}^{\infty}{{1}\over{l!}}\!\left[{{t^2}\over{2n}}
\sum_{k=m+1}^{n-1}{{1}\over{k-{\rm i}t}}\right]^{l}.
$$
Clearly, $C_{nm}(0) =0$, and for $t\not= 0$,
\begin{align*}
\big|C_{nm}(t)\big| &\leq
\sum_{l=2}^{\infty}{{1}\over{l!}}\left[{{t^2}\over{2n}}
\sum_{k=m+1}^{n-1}\left|{{1}\over{k-{\rm i}t}}\right|\right]^{l} =
\sum_{l=2}^{\infty}{{1}\over{l!}}\left[{{t^2}\over{2n}}
\sum_{k=m+1}^{n-1}{{1}\over{\sqrt{k^2+t^2}}}\right]^{l}\cr &\leq
\sum_{l=2}^{\infty}{{1}\over{l!}}\left[{{t^2}\over{2n}}
\int_{0}^{n}\!{{1}\over{\sqrt{y^2+t^2}}}\,\d y\right]^{l} =
\sum_{l=2}^{\infty}{{1}\over{l!}}
\left[{{t^2}\over{2n}}\int_{0}^{n/|t|}\! {1\over{\sqrt{x^2+1}}}\,\d
x\right]^{l}\cr &=
\sum_{l=2}^{\infty}{{1}\over{l!}}\left[{{t^2}\over{2n}}\, {\rm
rsh}\!{\left({{n}\over{|t|}}\right)}\right]^l =
\left[{{t^2}\over{2n}}\, {\rm
rsh}\!{\left({{n}\over{|t|}}\right)}\right]^2
\sum_{j=0}^{\infty}{{1}\over{j!}}\left[{{t^2}\over{2n}}\, {\rm
rsh}\!{\left({{n}\over{|t|}}\right)}\right]^j,
\end{align*}
where, for momentary gain of space, ${\rm rsh}(x) := {\rm
arcsinh}(x) = \log\big(x+\sqrt{x^2+1}\,\big)$ is the reverse
(inverse) to the hyperbolic function $\sinh(x) = (\e^{x} -
\e^{-x})/2$, $x\in\R$. Thus,
\begin{equation}\label{g_2.33}
\big|\rho_{nm,3}(t)\big|\leq
\big|A_{nm}(t)\big|\,{{t^4}\over{4n^2}}\, {\rm
rsh}^{2}\!{\left({{n}\over{|t|}}\right)}
\exp\!\left\{{{t^2}\over{2n}}\,{\rm
rsh}\!\left({{n}\over{|t|}}\right)\right\}.
\end{equation}
Now we can estimate $R_{nm,3} = I_{nm}^{3,1} + I_{nm}^{3,2}$ in
(\ref{g_2.15}).

Concerning the first term, notice first that for the derivatives of
the functions $v_n(s) = s^3\,{\rm rsh}^2(n/s)$ and $w_n(s) =
s^2\,{\rm rsh}(n/s)/(2n)$, for $s\in(0,1]$, we have
$$
v_n^{\prime}(s) = s^2\,{\rm rsh}\!\left({{n}\over{s}}\right)\!
\bigg[3\,{\rm rsh}\!\left({{n}\over{s}}\right) -{2n\over
\sqrt{s^2+n^2}}\bigg] \ge s^2\,{\rm
rsh}\!\left({{n}\over{s}}\right)\!\big[3\,{\rm rsh}(1) - 2\big]
$$
which is positive since $3\,{\rm rsh}(1) - 2 > 0.64$, and
$$
w_n^{\prime}(s) = {s\over 2n} \bigg[2\,{\rm
rsh}\!\left({{n}\over{s}}\right) -{n\over \sqrt{s^2+n^2}}\bigg] \ge
{s\over 2n}\big[2\,{\rm rsh}(1) - 1\big] > {s\over 2n}\,{76\over
100} > 0,
$$
so both $v_n(\cdot)$ and $\exp\{w_n(\cdot)\}$ are monotone
increasing on the interval $(0,1)$. Hence, using (\ref{g_2.8}), the fact
that $|A_{nm}(t)|\leq 1$ and the evenness of the functions involved,
we can bound the integrand by its value at $1$ to obtain
$$
I_{nm}^{3,1}\le {1\over 2n^2}\int_{0}^{1}v_n(s)\,\e^{w_n(s)}\,\d s
\le {\e^{\,{\rm arcsinh}(n)/(2n)}\over 2}\,{{\rm arcsinh}^2(n)\over
n}.
$$
The asymptotic equality ${\rm arcsinh}(n)\sim \log n$ then shows
that $I_{nm}^{3,1} = O\big((\log n)^2/n^2\big)$.

For $I_{nm}^{3,2}$, using (\ref{g_2.33}), (\ref{g_2.28}) with the $a_{nm}\to a_m$
given there, the evenness of the integrand and the fact that the
function $t\mapsto {\rm arcsinh}(n/t)$ is decreasing, we get
\begin{equation}\label{g_2.34}
I_{nm}^{3,2}\leq {{a_{nm}}\over{2}}\,{{\rm arcsinh}^2(n)\over n^2}
\int_{1}^{n}t^{m+4}\,\e^{-tg_{nm}(t)}\,\d t,
\end{equation}
where, with $t\in[1, n)$ everywhere in what remains of the present
case $j=3$ of (\ref{g_2.15}),
$$
g_{nm}(t) = \arctan{{t(n-m-1)}\over{t^2+n(m+1)}} -
{{t}\over{2n}}\;{\rm arcsinh}\!\left({{n}\over{t}}\right).
$$
Here
$$
g_{nm}^{\,\prime}(t) =
{{(n-m-1)\left[n(m+1)-t^2\right]}\over{\left[n(m+1)+t^2\right]^2 +
t^2(n-m-1)^2}} - {{1}\over{2}}\!\left({{{\rm arcsinh}(n/t)}\over{n}}
-{{1}\over{\sqrt{t^2+n^2}}}\right).
$$
The first term of this expression is positive if $t < \sqrt{n(m+1)}$
and negative if $t > \sqrt{n(m+1)}$, while the second term is
negative for all $t\in (0,n)$, since this term takes on a negative
value at $t=n$ and is an increasing function of $t$ on this
interval:
$$
-{{1}\over{2}}\! \left({{{\rm
arcsinh}(n/t)}\over{n}}-{{1}\over{\sqrt{t^2+n^2}}}\right)^{\!\prime}
= {{n^2}\over{2t(t^2+n^2)^{3/2}}}>0,\quad\hbox{if}\quad 0<t<n.
$$
In particular, $g_{nm}(\cdot)$ is a decreasing function on the
interval $(\sqrt{n(m+1)},n)$, and hence $g_{nm}(t)\ge g_{nm}(n)$ for
all $t\in(\sqrt{n(m+1)},n)$, where this value is
$$g_{nm}(n) =
\arctan{{n-m-1}\over{n+m+1}}-{{{\rm arcsinh}(1)}\over{2}} \to
\arctan(1)-{{{\rm arcsinh}(1)}\over{2}} > 0.34.
$$
Therefore, for all $n$ sufficiently large,
\begin{equation}\label{g_2.35}
\int_{\sqrt{n(m+1)}}^{n}t^{m+4}\,\e^{-tg_{nm}(t)}\,\d t
\leq\int_{\sqrt{n(m+1)}}^{n}t^{m+4}\,\e^{-(0.3)t}\,\d t = o(1).
\end{equation}
For $t\in[1,\sqrt{n(m+1)}\,)$ we replace the leading $t$ of the
second term in $g_{nm}(t)$ by $\sqrt{n(m+1)}$, and then notice that
the resulting lower bound is an increasing function of $t$ on the
interval $[1,\sqrt{n(m+1)}\,)$. This way we obtain
\begin{align*}
g_{nm}(t)&\geq \arctan{{t(n-m-1)}\over{t^2+n(m+1)}} -
{\sqrt{m+1}\over{2}}\, {{\rm
arcsinh}\!\left({{n}\over{t}}\right)\over \sqrt{n}}\cr &\geq
\arctan{{n-m-1}\over{1+n(m+1)}} - {\sqrt{m+1}\over{2}}\, {{\rm
arcsinh}(n)\over \sqrt{n}},
\end{align*}
and this converges to $\arctan(1/(m+1))>0$. Hence for all $n$ large
enough,
\begin{equation}\label{g_2.36}
\int_{1}^{\sqrt{n(m+1)}}t^{m+4}\,\e^{-tg_{nm}(t)}\,\d t \leq
\int_{1}^{\sqrt{n(m+1)}}t^{m+4}\,
\e^{-\left({{1}\over{2}}\arctan{{1}\over{m+1}}\right)\,t}\,\d t.
\end{equation}
Now (\ref{g_2.36}), (\ref{g_2.35}) and (\ref{g_2.34}) together give $I_{nm}^{3,2} =
O\big((\log n)^2/n^2\big)$ again. Thus $R_{nm,3} = O\big((\log
n)^2/n^2\big)$ as well, and so the case $j=3$ in (\ref{g_2.15}) is amply
satisfied.

\bigskip

{\cal The case of $R_{nm,4}$.} Using (\ref{g_2.10})--(\ref{g_2.12}) and the
formulae above (\ref{g_2.14}), we see that $|\rho_{nm,4}(t)| =
|\psi_{nm,3}(t)-\psi_{nm,4}(t)| =
|D_{nm}(t)|\,|\varphi_{m}(t)-A_{nm}(t)|$, where
\begin{align}\label{g_2.37}
\big|D_{nm}(t)\big| &=
\left|1+{{t^2}\over{2n}}\sum_{k=m+1}^{n-1}{{1}\over{k-{\rm
i}t}}\right| \le
1+{{t^2}\over{2n}}\sum_{k=m+1}^{n-1}{{1}\over{\sqrt{k^2+t^2}}}\cr
&\le 1+{{t^2}\over{2n}}\sum_{k=1}^{2n}{{1}\over{|t|}} = 1+|t| \le
\left\{
  \begin{array}{ll}
    2, & \hbox{if $|t|\leq 1$;} \\
    2|t|, & \hbox{if $|t|>1$.}
  \end{array}
\right.
\end{align}
With this factor done, inequality (\ref{g_2.16}) gives another analogue of
(\ref{g_2.17}), namely
\begin{equation}\label{g_2.38}
|\rho_{nm,4}(t)|\leq {{|D_{nm}(t)|}\over{2}}\big\{|\varphi_{m}(t)| +
|A_{nm}(t)|\big\} \,\delta_{mn}^{[4]}(t),
\end{equation}
where, since the assumption $|t|\le cn$ for some $c\in(0,1)$ ensures
that the logarithms can be expanded for all $k=n, n+1,\ldots$, as
noted at (\ref{g_2.13}),
\begin{align*}
\delta_{mn}^{[4]}(t) &=
\left|\sum_{k=n}^{\infty}\left[-{{{\rm i}t}\over{k}}
-\log\!\left(1-{{{\rm i}t}\over{k}}\right)\right]\right| \leq
\sum_{k=n}^{\infty}\sum_{j=2}^{\infty}{{1}\over{j}}{{|t|^j}\over{k^j}}
\leq
\sum_{j=2}^{\infty}{{|t|^j}\over{j}}\int_{n-1}^{\infty}{{1}\over{y^j}}\,\d
y\cr & = \sum_{j=2}^{\infty}{{|t|^j}\over{j(j-1)}}\,
{{1}\over{n^{j-1}}}\,{{n^{j-1}}\over{(n-1)^{j-1}}} \leq
{{t^2}\over{n}}\,{{n}\over{n-1}}\!
\left[{{1}\over{2}}+{{1}\over{6}}\sum_{j=3}^{\infty}
\left({{|t|}\over{n}}\,{{n}\over{n-1}}\right)^{\!j-2}\right].
\end{align*}
Let $n$ be large enough to make $n/(n-1)< (1+c)/(2c)$, where $c$ is
as above. Then
$$
{{|t|}\over{n}}{{n}\over{n-1}}<
c\,{{1+c}\over{2c}}={{1+c}\over{2}}<1,
$$
and so
$$
\sum_{j=3}^{\infty}
\left({{|t|}\over{n}}{{n}\over{n-1}}\right)^{\!j-2} \le
{{1}\over{1-{{1+c}\over{2}}}}-1 = {{1+c}\over{1-c}}.
$$
Hence, for all $n$ large enough,
\begin{equation}\label{g_2.39}
\delta_{mn}^{[4]}(t) \leq {{t^2}\over{n}}\,{{1+c}\over{2c}}
\left[{{1}\over{2}}+{{1}\over{6}}\,{{1+c}\over{1-c}}\right]
\leq{{t^2}\over{n}}\,{{(1+c)^2}\over{c(1-c)}}.
\end{equation}

Next, by classical results on the $\Gamma$ function (\cite{N}, \S 8, for
example) we see that
\begin{align*}
\left|\varphi_{m}(t)\right|&=
\left|\prod_{k=1}^{m}\e^{{\rm i}t/k}\!\left(1-{{{\rm
i}t}\over{k}}\right)\right| \left|\prod_{k=1}^{\infty}{{\e^{-{\rm
i}t/k}}\over{1-{{{\rm i}t}\over{k}}}}\right| = \left|\Gamma(1-{\rm
i}t)\e^{-{\rm i}t\gamma}\right| \prod_{k=1}^{m}\left|\e^{{\rm
i}t/k}\!\left(1-{{{\rm i}t}\over{k}}\right)\right|\cr
&=\left|\Gamma(1-{\rm i}t)\right|
\sqrt{\prod_{k=1}^{m}\left(1+{{t^2}\over{k^2}}\right)} =
{{\sqrt{2\pi|t|}}\over{\sqrt{1-\e^{-2\pi|t|}}}}\,\e^{-{{\pi}\over{2}}|t|}
\sqrt{\prod_{k=1}^{m}\left(1+{{t^2}\over{k^2}}\right)},
\end{align*}
starting from (\ref{g_2.11}), which implies
\begin{equation}\label{g_2.40}
\left|\varphi_{m}(t)\right| \leq
{{\sqrt{2\pi}}\over{\sqrt{1-\e^{-2\pi}}}}\,\sqrt{|t|(1+t^2)^m}\,\e^{-
\pi |t|/2}, \quad\hbox{whenever}\quad |t|\in[1,\infty).
\end{equation}
Also, since $(m+1)n + 1\ge n + m + 1$, combining the two cases in
(\ref{g_2.29}) we have
\begin{equation}\label{g_2.41}
|A_{nm}(t)|\leq a_{nm}|t|^{m+1}
\exp\!\left\{-|t|\arctan{{n-m-1}\over{1+(m+1)n}}\right\},
\quad\hbox{if}\quad |t|\in[1,cn],
\end{equation}
and we are ready to deal with $R_{nm,4} = I_{mn}^{4,1} +
I_{mn}^{4,2}$ in (\ref{g_2.15}).

The inequalities $|\varphi_{m}(t)|\leq 1$ and $|A_{nm}(t)|\leq 1$
and (\ref{g_2.37})--(\ref{g_2.39}) imply
$$I_{mn}^{4,1}
\leq {{1}\over{n}}\,{{2(1+c)^2}\over{c(1-c)}}\int_{-1}^{1}|t|\d t =
{{2(1+c)^2}\over{c(1-c)}}\,{{1}\over{n}}
$$
for all $n$ large enough. Also, collecting the five bounds from
(\ref{g_2.37})--(\ref{g_2.41}), for all sufficiently large $n$ we obtain
\begin{align*}
I_{mn}^{4,2} &\leq{{a_{nm}}\over{n}}\,{{2(1+c)^2}\over{c(1-c)}}
\int_{1}^{\infty}t^{m+3}
\exp\!\left\{-|t|\arctan{{n-m-1}\over{1+(m+1)n}}\right\}\d t\cr
&\;\;\;\,+ {{1}\over{n}}\,
{{2(1+c)^2}\sqrt{2\pi}\over{c(1-c)}{\sqrt{1-\e^{-2\pi}}}}
\int_{1}^{\infty}t^{5/2}\,\big(1+t^2\big)^{m/2}\,\e^{-\pi t/2}\,\d
t.
\end{align*}
Since $a_{nm}\to a_m$ and $\arctan([n-m-1]/[1+(m+1)n])\to
\arctan(1/(m+1))>0$, we conclude that $I_{mn}^{4,2} = O(1/n)$.
Therefore, $R_{nm,4} = O(1/n)$, establishing the case $j=4$ of
(\ref{g_2.15}) and thus completing the proof of (\ref{g_1.7}) in the
theorem.

\bigskip

\noindent{\bf Proof of (\ref{g_1.10}).} Let
$E_{n}=\sum_{k=1}^{n}{{1}\over{k}}-\log n-\gamma$, for which, by a
classical asymptotic expansion due to Euler,
$$
E_{n} = {{1}\over{2n}}-{{1}\over{12n^2}}+{{\theta_{n}}\over{120n^4}}
\quad\hbox{for some}\quad \theta_{n}\in(0,1).
$$
Using the forms of the three ingredients given before the statement
of the theorem, for the deviation $\Delta_n^* := \sup_{x\in\R}
\big|F_{nm}^{*}(x) - [F_{m}^{*}(x) + G_{n,m}^{*}(x)]\big|$ in
question we obtain
\begin{align*}
\Delta_n^* &= \sup_{x\in\R}\big|F_{n,m}(x-C_m-E_n) -
F_m(x-C_m) - G_{n,m}(x-C_m-E_n)\big|\cr
&=\sup_{y\in\R}\big|F_{n,m}(y)- F_m(y+E_n) - G_{n,m}(y)\big|\cr &\le
\sup_{y\in\R}\big|F_{n,m}(y) - [F_m(y) + G_{n,m}(y)]\big| +
\sup_{y\in\R}\big|F_{m}(y)-F_{m}(y+E_{n})\big|.
\end{align*}
Hence (\ref{g_1.7}) and the inequality, obtained by the mean value theorem,
$$\sup_{y\in\R}\big|F_{m}(y)-F_{m}(y+E_{n})\big|
\le |E_{n}| \, \max_{y\in\R}f_{m}(y) = O\!\bigg({{1}\over{n}}\bigg)
$$
yield the desired statement in (\ref{g_1.10}). Also, the functions
$G_{n,m}^{\diamond}(x) := G_{n,m}(x - C_m)$ clearly inherit the
stated order properties of $G_{n,m}(x)$, $x\in\R$, since the shift
$C_m$ is constant. Then so do the functions $G_{n,m}^{*}(x) =
G_{n,m}^{\diamond}(x - E_n)$, $x\in\R$, because
$$
\sup_{x\in\R}\big|G_{n,m}^{\diamond}(x) - G_{n,m}^{*}(x)\big| \le
|E_{n}| \, \sup_{x\in\R}\big|G_{n,m}^{\,\prime}(x)\big| =
O\!\bigg({\log n\over n^2}\bigg),
$$
where the last bound is due to the inequality in
(\ref{g_2.9}). $\blacksquare$


\chapter{Normal approximation}

In this chapter we prove an error bound for normal approximation to the coupon collector's standardized waiting time. We introduce the distribution
functions
$$
F_{n,m}(x):=\p\bigg({{W_{n,m}-\mu_n}\over{\sigma_n}}\leq
x\bigg),\quad x\in\R.
$$
As mentioned in the Introduction, Baum and Billingsley \cite{BB} showed that if the $m$ goes to infinity along with $n$,
but slowly enough to let the sequence $(n-m)/\sqrt{n}$ tend to
infinity as-well, then the standardized
$W_{n,m}$ is asymptotically normal:
$$
\lim_{n\to\infty}F_{n,m}(x)=\Phi(x),\quad\hbox{where}\quad
\Phi(x):={{1}\over{\sqrt{2\pi}}}\int_{-\infty}^{x}\e^{-s^2/2}\d s,\;\;\;
x\in\R.
$$
The following theorem gives a bound for the rate of convergence in this central limit
theorem.

\bigskip

\begin{theorem}\label{t_normalis} For all $n\geq3$ and $1\leq m\leq n-2$, we have
\begin{equation}\label{normal_approx1}
\sup_{x\in \R}\big|F_{n,m}(x)-\Phi(x)\big|\leq C\frac{n}{m}\frac{1}{\sigma_n},
\end{equation}
where $C=9.257$.
\end{theorem}

\bigskip

One can check that the bound given by Theorem \ref{t_normalis} goes to 0 iff $m$ goes to infinity along with $n$, but slowly enough to let the sequence $(n-m)/\sqrt{n}$ tend to infinity as-well, which is in accord with the central limit theorem stated above.
Indeed, this follows easily from the
asymptotic formulae for $\sigma_n^2$, given by Baum and Billingsley \cite{BB}:

If ${{m}\over{n}}\to d$ for some $d\in(0,1)$, so that
${{n-m}\over{n}}\to 1-d$, then $\sigma_n^2\sim n{{1-d+d\log
d}\over{d}}$.

If ${{m}\over{n}}\to1$, so that ${{n-m}\over{n}}\to 0$, and
${{(n-m)^2}\over{n}}\to\infty$, then
$\sigma_n^2\sim{1\over2}{{(n-m)^2}\over{n}}$.

If ${{m}\over{n}}\to 0$, so that ${{n-m}\over{n}}\to 1$, and
$m\to\infty$, then $\sigma_n^2\sim {{n^2}\over{m}}$.

\medskip

\noindent These asymptotic relations then give the following typical
examples:

\medskip

If $m\sim dn$ for some $0<d<1$, then $\frac{m}{n}\sigma_n\sim {\rm constant}\cdot\sqrt{n}$.

If $m\sim n-n^\alpha$ for some ${1\over2}<\alpha<1$, then $\frac{m}{n}\sigma_n\sim{\rm constant}\cdot n^{\alpha-{1\over2}}$.

If $m\sim n^\beta$ for some $0<\beta<1$, then $\frac{m}{n}\sigma_n\sim{\rm constant}\cdot n^{\beta/2}$.

If $m\sim \log n$, then $\frac{m}{n}\sigma_n\sim{\rm constant}\cdot\sqrt{\log n}$.


\bigskip

\noindent{\bf Proof.} We estimate the supremum distance between the
distribution function $F_{n,m}$ and the limiting
distribution function $\Phi$ in terms of their characteristic
functions, using Esseen's smoothing inequality. Since the
characteristic function of the geometric distribution with success
probability $p\in(0,1)$ is $p\e^{\i t}/(1-q\e^{\i t})$, where $\i$ is the
imaginary unit and $q=1-p$, by (\ref{Wgeo}) we have
\begin{align*}
\varphi_{n,m}(t)&:=\int_{-\infty}^{\infty}\e^{\i xt}\d F_{n,m}(x)=
\E\!\left(\exp\!\left\{{\i t\over\sigma_n}\left[W_{n,m}-n\sum_{k=m+1}^{n}{1\over
k}\right]\right\}\right)\cr &=\prod_{k=m+1}^{n-1}{{{k\over
n}\e^{\i t/\sigma_n}\e^{-\i tn/k\sigma_n}}\over{1-{n-k\over
n}\e^{\i t/\sigma_n}}},
\end{align*}
while the limiting characteristic function is
$\int_{-\infty}^{\infty}\e^{\i xt}\d\Phi(x)=\e^{-t^2/2}$, $t\in\R$.
Choosing the main parameter in Esseen's inequality (see Section 2.2.) to
be $c_n\sigma_n$,
where, with any fixed $c\in(0,1)$, the sequence $c_n(m) = c_n$ is
given by
\begin{equation}\label{n_3}
c_n(m) :=
\min\!\left\{1,{{c(m+1)}\over{\sqrt{n(n-m-1)}}}\right\},
\end{equation}
the inequality in the case of our distribution
functions takes on the following form:
\begin{equation}\label{n_4}
\sup_{x\in\R}\big|F_{n,m}(x)-\Phi(x)\big|\leq{{b}\over{2\pi}}
\int_{-c_n\sigma_n}^{c_n\sigma_n}\left|{{\varphi_{n,m}(t)-\e^{-t^2/2}}\over{t}}\right|\d t
+{{c_b}\over{\sqrt{2\pi}}}{{1}\over{c_n\sigma_n}},
\end{equation}
where $b>1$ is arbitrary and $c_b$ is a positive constant depending
only on $b$.

Since we restricted the domain of the characteristic functions to
$(-c_n\sigma_n,c_n\sigma_n)$, from now on we assume
that this interval is the domain of all formulae containing the
variable $t$. We emphasize that by the definition of $c_n$ in (\ref{n_3})
this means that, on the one hand, $|t|<\sigma_n$, and, on the other
hand, for any $c\in(0,1)$ chosen in $c_n$ and
$k\in\{m+1,\ldots,n-1\}$,
\begin{equation}\label{n_5}
|t|<{{c(m+1)}\over{\sqrt{n(n-m-1)}}}\sigma_n\leq{ck\over{\sqrt{n(n-k)}}}\sigma_n
<{k\over{\sqrt{n(n-k)}}}\sigma_n.
\end{equation}

We estimate the deviation $|\varphi_{n,m}(t)-\e^{-t^2/2}|$ in the
integrand on the right-hand side of the inequality in (\ref{n_4}) through
the following heuristic steps:
\begin{align*}
&\varphi_{n,m}(t)
=\exp\!\left\{\sum_{k=m+1}^{n-1}\!\left[{{n-k}\over
k}{{-\i t}\over{\sigma_n}}+\log {k\over n}-\log\!\left(1-{{n-k}\over
n}\e^{\i t/\sigma_n}\right)\right]\right\}\cr
&\approx\exp\!\left\{\sum_{k=m+1}^{n-1}\!\left[{{n-k}\over
k}{{-\i t}\over{\sigma_n}}+\log {k\over n}-\log\!\left(1-{{n-k}\over
n}\left(1+{{\i t}\over{\sigma_n}}-{{t^2}\over{2\sigma_n^2}}\right)\right)\right]\right\}\cr
&=\exp\!\left\{\sum_{k=m+1}^{n-1}\!\left[{{n-k}\over
k}{{-\i t}\over{\sigma_n}} -\log\!\left(1-{{n-k}\over
k}\left({{\i t}\over{\sigma_n}}-{{t^2}\over{2\sigma_n^2}}\right)\right)\right]\right\}
:=\varphi_{n,m}^{[1]}(t)\cr
&\approx\exp\!\left\{\sum_{k=m+1}^{n-1}\!\left[{{n-k}\over
k}{{-\i t}\over{\sigma_n}}+{{n-k}\over
k}\Bigg({{\i t}\over{\sigma_n}}-{{t^2}\over{2\sigma_n^2}}\Bigg)\!+{1\over
2}{{(n-k)^2}\over
k^2}\left({{\i t}\over{\sigma_n}}-{{t^2}\over{2\sigma_n^2}}\right)^2\right]\right\}\cr
&=\exp\!\left\{\sum_{k=m+1}^{n-1}\!\left[-{{n(n-k)}\over{2k^2}}{{t^2}
\over{\sigma_n^2}}-{{(n-k)^2}\over{2k^2}}{{\i t^3}\over{\sigma_n^3}}
+{{(n-k)^2}\over{8k^2}}{{t^4}\over{\sigma_n^4}}\right]\right\}\cr
&=\exp\!\left\{-{{t^2}\over{2}}\right\}\exp\!\left\{\sum_{k=m+1}^{n-1}\!{{(n-k)^2}\over{2k^2}}\left({{-\i t^3}\over{\sigma_n^3}}
+{{t^4}\over{4\sigma_n^4}}\right)\right\}
=:\varphi_{n,m}^{[2]}(t)\cr
&\approx\exp\!\left\{-{{t^2}\over{2}}\right\}.
\end{align*}
At each approximation a certain function was replaced with the first
few terms of its series expansion. At the second one this was done
with a logarithmic expression, whose expansion about 1 exists,
because for an arbitrary term of the sum in
$\varphi_{n,m}^{[1]}(\cdot)$, that is, for an arbitrary
$k\in\{m+1,\ldots,n-1\}$,
$$
{{n-k}\over{k}}\left|{{\i t}\over{\sigma_n}}-{{t^2}\over{2\sigma_n^2}}\right|
={n-k\over k}{|t|\over\sigma_n}\sqrt{1+{t^2\over{4\sigma_n^2}}}
<{n-k\over k}{{ck}\over{\sqrt{n(n-k)}}}\sqrt{1+{k^2\over{4n(n-k)}}}
$$
by (\ref{n_5}), which gives
\begin{equation}\label{n_6}
{{n-k}\over{k}}\left|{{\i t}\over{\sigma_n}}-{{t^2}\over{2\sigma_n^2}}\right|
<c{{2n-k}\over{2n}}\leq c<1.
\end{equation}

Now, the errors resulting from the first two of our three
approximations can be estimated applying the following inequality
$$
\left|\e^{z_1}-\e^{z_2}\right|\leq{{1}\over{2}}\big\{|\e^{z_1}|+|\e^{z_2}|\big\}|z_1-z_2|
$$
for arbitrary complex numbers $z_1$ and $z_2$. This yields
\begin{equation}\label{n_7}
\big|\varphi_{n,m}(t)-\varphi_{n,m}^{[1]}(t)\big|
\leq{{1}\over{2}}\Big\{|\varphi_{n,m}(t)|+\big|\varphi_{n,m}^{[1]}(t)\big|\Big\}
\delta_{n,m}^{[1]}(t),
\end{equation}
where
$$
\delta_{n,m}^{[1]}(t)=\left|\sum_{k=m+1}^{n-1}\!\left[\log\!\left(1-{{n-k}\over
n}\e^{\i t/\sigma_n}\right)\!-\log\!\left(1-{{n-k}\over
n}\left(1+{{\i t}\over{\sigma_n}}-{{t^2}\over{2\sigma_n^2}}\right)\right)\right]\right|,
$$
and
\begin{equation}\label{n_8}
\big|\varphi_{n,m}^{[1]}(t)-\varphi_{n,m}^{[2]}(t)\big|
\leq{{1}\over{2}}\Big\{|\varphi_{n,m}^{[1]}(t)|+\big|\varphi_{n,m}^{[2]}(t)\big|\Big\}
\delta_{n,m}^{[2]}(t),
\end{equation}
where
\begin{align*}
\quad \delta_{n,m}^{[2]}(t)
&=\Bigg|\sum_{k=m+1}^{n-1}\!\Bigg[\log\!\left(1-{{n-k}\over
k}\left({{\i t}\over{\sigma_n}}-{{t^2}\over{2\sigma_n^2}}\right)\right)\cr
&\quad +{{n-k}\over k}
\Bigg({{\i t}\over{\sigma_n}}-{{t^2}\over{2\sigma_n^2}}\Bigg)\!+{1\over
2}{{(n-k)^2}\over
k^2}\left({{\i t}\over{\sigma_n}}-{{t^2}\over{2\sigma_n^2}}\right)^2\Bigg]\Bigg|\cr
&=\left|\sum_{k=m+1}^{n-1}\sum_{j=3}^{\infty}\!\left[{{-1}\over{j}}
\left({{n-k}\over{k}}\right)^j\left({{\i t}\over{\sigma_n}}-{{t^2}
\over{2\sigma_n^2}}\right)^j\right]\right|.
\end{align*}

Summarizing, for the estimation of the integral
$$I_{n,m}:=\int_{-c_n\sigma_n}^{c_n\sigma_n}\left|{{\varphi_{n,m}(t)-\e^{-t^2/2}}\over{t}}\right|\d t$$
in (\ref{n_4}), we use the intermediate approximative functions
$\varphi_{n,m}^{[1]}(\cdot)$ and $\varphi_{n,m}^{[2]}(\cdot)$,
and the inequalities above concerning their differences, obtaining
have
\begin{align}\label{n_9}
I_{n,m}&\leq
{{1}\over{2}}\int_{-c_n\sigma_n}^{c_n\sigma_n}
\left|{{\varphi_{n,m}(t)}\over{t}}\right|\delta_{n,m}^{[1]}(t)\d t
+ {{1}\over{2}}\int_{-c_n\sigma_n}^{c_n\sigma_n}
\left|{{\varphi_{n,m}^{[1]}(t)}\over{t}}\right|\delta_{n,m}^{[1]}(t)\d t\cr
&\,\,\,\,\,\,+ {{1}\over{2}}\int_{-c_n\sigma_n}^{c_n\sigma_n}
\left|{{\varphi_{n,m}^{[1]}(t)}\over{t}}\right|\delta_{n,m}^{[2]}(t)\d t
+ {{1}\over{2}}\int_{-c_n\sigma_n}^{c_n\sigma_n}
\left|{{\varphi_{n,m}^{[2]}(t)}\over{t}}\right|\delta_{n,m}^{[2]}(t)\d t\cr
&\,\,\,\,\,\,+ \int_{-c_n\sigma_n}^{c_n\sigma_n}
\left|{{\varphi_{n,m}^{[2]}(t)-\e^{-t^2/2}}\over{t}}\right|\d t.
\end{align}
Now we give upper bounds for each of the functions occurring in the
integrals above.

First we consider $|\varphi_{n,m}(t)|$. By simple computation
\begin{align*}
\big|\varphi_{n,m}(t)\big|&=\left|\exp\!\left\{\sum_{k=m+1}^{n-1}\left[{{n-k}\over{k}}{{-it}\over{\sigma_n}}+\log
{k\over n}-\log\!\left(1-{{n-k}\over
n}\e^{\i t/\sigma_n}\right)\right]\right\}\right|\cr
&=\exp\!\left\{\sum_{k=m+1}^{n-1}\left[\log {k\over
n}-\log\!\left|1-{{n-k}\over
n}\e^{\i t/\sigma_n}\right|\right]\right\}\cr
&=\exp\!\left\{\sum_{k=m+1}^{n-1}\left[\log {k\over
n}-\log\sqrt{1+\left({{n-k}\over n}\right)^2-2{{n-k}\over
n}\cos{{t}\over{\sigma_n}}}\right]\right\}\cr
&=\exp\!\left\{\sum_{k=m+1}^{n-1}\left[\log {k\over
n}-\log\sqrt{\left(k\over n\right)^2+2{{n-k}\over
n}\left(1-\cos{{t}\over{\sigma_n}}\right)}\right]\right\}\cr
&=\exp\!\left\{-{{1}\over{2}}\sum_{k=m+1}^{n-1}\log\!\left(1+2{{n(n-k)}\over
k^2}\left(1-\cos{{t}\over{\sigma_n}}\right)\right)\right\}.
\end{align*}
Since $t/\sigma_n<1<\pi/2$, we can continue with applying the
inequality $1-\cos x \geq {{4}\over{\pi^2}}x^2$, true for
$x\in(0,\pi/2)$, and obtain
$$
|\varphi_{n,m}(t)|\leq\exp\!\left\{-{{1}\over{2}}\sum_{k=m+1}^{n-1}\log\!\left(1+{8\over\pi^2}{{n(n-k)}\over{k^2}}{{t^2}\over{\sigma_n^2}}
\right)\right\}.
$$
Now we see from (\ref{n_5}) that we can use the inequality
\begin{equation}\label{n_10}
\log(1+x)\geq x-{{x^2}\over{2}},\;\; \hbox{if}\;\; x\in(0,1),
\end{equation}
which yields
$$
|\varphi_{n,m}(t)|\leq\exp\!\left\{-{{1}\over{2}}\sum_{k=m+1}^{n-1}\left[{{n(n-k)}\over{k^2}}{{t^2}\over{\sigma_n^2}}\left(
{8\over\pi^2}-{32\over\pi^4}{{n(n-k)}\over{k^2}}{{t^2}\over{\sigma_n^2}}\right)\right]\right\}.
$$
Using once again the bound in (\ref{n_5}), on the second $t^2$ in the
expression above, and recalling (\ref{0var}), we easily get
\begin{equation}\label{n_11}
|\varphi_{n,m}(t)|\leq
\exp\!\left\{-{{20}\over{\pi^4}}t^2\right\}.
\end{equation}

Next, in a completely analogous way,
\begin{align*}
\big|\varphi_{n,m}^{[1]}(t)\big|&=
\left|\exp\!\left\{\sum_{k=m+1}^{n-1}\!\left[{{n-k}\over
k}{{-\i t}\over{\sigma_n}} -\log\!\left(1-{{n-k}\over
k}\left({{\i t}\over{\sigma_n}}-{{t^2}\over{2\sigma_n^2}}\right)\right)\right]\right\}\right|\cr
&=\exp\!\left\{-\sum_{k=m+1}^{n-1}\log\!\left|1-{{n-k}\over
k}\left({{\i t}\over{\sigma_n}}-{{t^2}\over{2\sigma_n^2}}\right)\right|\right\}\cr
&=\exp\!\left\{-{{1}\over{2}}\sum_{k=m+1}^{n-1}\log\!\left(1+{{n(n-k)}\over{k^2}}{{t^2}\over{\sigma_n^2}}
+{{(n-k)^2}\over{k^2}}{{t^4}\over{4\sigma_n^4}}\right)\right\}\cr
&\leq\exp\!\left\{-{{1}\over{2}}\sum_{k=m+1}^{n-1}\log\!\left(1+{{n(n-k)}\over{k^2}}{{t^2}\over{\sigma_n^2}}\right)\right\}.
\end{align*}
Again, (\ref{n_5}) allows us to use the inequality in (\ref{n_10}) to obtain
$$
\big|\varphi_{n,m}^{[1]}(t)\big|\leq\exp\!\left\{-{{1}\over{2}}\sum_{k=m+1}^{n-1}\!
\left[{{n(n-k)}\over{k^2}}{{t^2}\over{\sigma_n^2}}\left(1-{{n(n-k)}\over{2k^2}}{{t^2}\over{\sigma_n^2}}\right)\right]\right\},
$$
and now, another application of (\ref{n_5}), and recognizing (\ref{0var}), yields
\begin{equation}\label{n_12}
\big|\varphi_{n,m}^{[1]}(t)\big|\leq\exp\!\left\{-{{t^2}\over{4}}\right\}.
\end{equation}

Also, by elementary considerations and (\ref{0var}) again,
\begin{align*}
\big|\varphi_{n,m}^{[2]}(t)\big|
&=\left|\exp\!\left\{-{{t^2}\over{2}}\right\}\exp\!\left\{\sum_{k=m+1}^{n-1}
{{(n-k)^2}\over{2k^2}}\left({{-\i t^3}\over{\sigma_n^3}}+{{t^4}\over{4\sigma_n^4}}\right)\right\}\right|\cr
&=\exp\!\left\{-{{t^2}\over{2}}\right\}\exp\!\left\{\sum_{k=m+1}^{n-1}{{(n-k)^2}\over{k^2}}{{t^4}\over{8\sigma_n^4}}\right\}\cr
&\leq\exp\!\left\{-{{t^2}\over{2}}\right\}\exp\!\left\{\sum_{k=m+1}^{n-1}{{n(n-k)}\over{k^2}}{{t^4}\over{8\sigma_n^4}}\right\}
=\exp\!\left\{-{{t^2}\over{2}}\right\}\exp\!\left\{{{t^4}\over{8\sigma_n^2}}\right\}\cr
&\leq\exp\!\left\{-{{t^2}\over{2}}\right\}\exp\!\left\{{{t^2}\over{8}}\right\},
\end{align*}
where the last inequality follows from the assumption that
$|t|<\sigma_n$. Therefore
\begin{equation}\label{n_13}
\big|\varphi_{n,m}^{[2]}(t)\big|\leq\exp\!\left\{-{{3t^2}\over{8}}\right\}.
\end{equation}

Next, we see for $\delta_{n,m}^{[1]}(t)$ in (\ref{n_7}) that
\begin{align*}
\delta_{n,m}^{[1]}(t)
&=\left|\sum_{k=m+1}^{n-1}\!\log{{1-{{n-k}\over
n}\e^{\i t/\sigma_n}}\over{1-{{n-k}\over
n}\left(1+{{\i t}\over{\sigma_n}}-{{t^2}\over{2\sigma_n^2}}\right)}}\right|\cr
&=\left|\sum_{k=m+1}^{n-1}\!\log\!\left(1+{{{{n-k}\over
n}\left(1+{{\i t}\over{\sigma_n}}-{{t^2}\over{2\sigma_n^2}}\right)-{{n-k}\over
n}\e^{\i t/\sigma_n}}\over {1-{{n-k}\over
n}\left(1+{{\i t}\over{\sigma_n}}-{{t^2}\over{2\sigma_n^2}}\right)}}\right)\right|\cr
&=\left|\sum_{k=m+1}^{n-1}\!\log\!\left(1+{{1+{{\i t}\over{\sigma_n}}-{{t^2}\over{2\sigma_n^2}}-\e^{\i t/\sigma_n}}\over
{{k\over{n-k}}-{{\i t}\over{\sigma_n}}+{{t^2}\over{2\sigma_n^2}}}}\right)\right|
=:\!\left|\sum_{k=m+1}^{n-1}\!\log\!\left(1+z_{n,k}(t)\right)\right|.
\end{align*}
Since $|\e^{\i u}-(1+\i u-{{u^2}\over{2}})|\leq{{|u|^3}\over{6}}$ for all
$u\in\R$ and $|t|<\sigma_n$,
$$
|z_{n,k}(t)|\leq{{{{|t|^3}\over{6\sigma_n^3}}}\over
{\sqrt{{{k^2}\over{(n-k)^2}}+{{t^4}\over{4\sigma_n^4}}+{{n}\over{n-k}}{{t^2}\over{\sigma_n^2}}}}}
\leq{{{{|t|^3}\over{6\sigma_n^3}}}\over
{\sqrt{{{t^4}\over{4\sigma_n^4}}}}}
\leq{{{1}\over{3}}{{t^2}\over{2\sigma_n^2}}\over
{\sqrt{{{t^4}\over{4\sigma_n^4}}}}}={1\over3}<1,
$$
the logarithmic expression can be expanded, so that
\begin{align*}
\delta_{n,m}^{[1]}(t)
&=\left|\sum_{k=m+1}^{n-1}\sum_{j=1}^{\infty}
{{(-1)^{j+1}z_{n,k}^j(t)}\over{j}}\right|
\leq\sum_{k=m+1}^{n-1}\sum_{j=1}^{\infty}{{\left|z_{n,k}(t)\right|^j}\over{j}}\cr
&\leq\sum_{k=m+1}^{n-1}\!\left[z_{n,k}(t)\left\{1+{1\over2}\sum_{j=1}^{\infty}\left|z_{n,k}(t)\right|^j\right\}\right].
\end{align*}
Using the upper bounds just given for $|z_{n,k}(t)|$, we obtain
\begin{align*}
\delta_{n,m}^{[1]}(t)
&\leq\sum_{k=m+1}^{n-1}\!\left[{{{{|t|^3}\over{6\sigma_n^3}}}\over
{\sqrt{{{k^2}\over{(n-k)^2}}+{{t^4}\over{4\sigma_n^4}}+{{n}\over{n-k}}{{t^2}\over{\sigma_n^2}}}}}
\left\{1+{{1}\over{2}}\sum_{j=1}^{\infty}\!\left({{1}\over
{3}}\right)^j\right\}\right]\cr
&={{5}\over{24}}{{|t|^3}\over{\sigma_n^3}}\sum_{k=m+1}^{n-1}{{1}\over{\sqrt{{{k^2}\over{(n-k)^2}}+{{t^4}\over{4\sigma_n^4}}
+{{n}\over{n-k}}{{t^2}\over{\sigma_n^2}}}}}
\leq{{5}\over{24}}{{|t|^3}\over{\sigma_n^3}}\sum_{k=m+1}^{n-1}{{n-k}\over{k}}\cr
&\leq{{5}\over{24}}{{|t|^3}\over{\sigma_n^3}}\sum_{k=m+1}^{n-1}{{n(n-k)}\over{k^2}},
\end{align*}
so that recalling (\ref{0var}) again,
\begin{equation}\label{n_14}
\delta_{n,m}^{[1]}(t)\leq{{5}\over{24}}{{|t|^3}\over{\sigma_n}}.
\end{equation}

We now turn to $\delta_{n,m}^{[2]}(t)$ in (\ref{n_8}). Clearly,
\begin{align*}
\delta_{n,m}^{[2]}(t)&\leq\sum_{k=m+1}^{n-1}\sum_{j=3}^{\infty}\left[{{1}\over{j}}\left({{n-k}\over{k}}\right)^j
\left|{{-\i t}\over{\sigma_n}}+{{t^2}\over{2\sigma_n^2}}\right|^j\right]\cr
&\leq{1\over3}\sum_{k=m+1}^{n-1}\left[{{(n-k)^3}\over{k^3}}
\left|{{-\i t}\over{\sigma_n}}+{{t^2}\over{2\sigma_n^2}}\right|^3\left\{\sum_{j=0}^{\infty}\left({{n-k}\over{k}}\right)^j
\left|{{-\i t}\over{\sigma_n}}+{{t^2}\over{2\sigma_n^2}}\right|^j\right\}\right].
\end{align*}
By (\ref{n_6}), the infinite sum here is not greater than
$\sum_{j=0}^{\infty}c^j= 1/(1-c)$, for any fixed $c\in(0,1)$ chosen
in the definition of $c_n$, and this, together with the inequality
$$
\left|{{-\i t}\over{\sigma_n}}+{{t^2}\over{2\sigma_n^2}}\right|^3
={|t|^3\over\sigma_n^3}\left(1+{t^2\over{4\sigma_n^2}}\right)^{3/2}
\leq{|t|^3\over\sigma_n^3}\left({5\over4}\right)^{3/2},
$$
which is true because $|t|<\sigma_n$ was assumed, gives
$$
\delta_{n,m}^{[2]}(t)
\leq{{5\sqrt{5}}\over{24(1-c)}}{|t|^3\over\sigma_n^3}
\sum_{k=m+1}^{n-1}{{(n-k)^3}\over{k^3}} =:
{{5\sqrt{5}}\over{24(1-c)}}{|t|^3\over\sigma_n^3}\,r_{n,m}.
$$
For the remaining sum $r_{n,m}$ here, by (\ref{0var}) we obtain
\begin{align*}
r_{n,m}
&=\sum_{k=m+1}^{n-1}{{(n-k)^2}\over{nk}}{{n(n-k)}\over{k^2}}
\leq{{(n-m-1)^2}\over{n(m+1)}}\sum_{k=m+1}^{n-1}{{n(n-k)}\over{k^2}}\cr
&={{(n-m-1)^2}\over{n(m+1)}}\sigma_n^2
\leq{{n-m-1}\over{m+1}}
\sigma_n^2\leq{{\sqrt{n(n-m-1)}}\over{m+1}}\sigma_n^2
\leq{{\sigma_n^2}\over{c_n}},
\end{align*}
and we conclude
\begin{equation}\label{n_15}
\delta_{n,m}^{[2]}(t)\leq{{5\sqrt{5}}\over{24(1-c)}}{|t|^3\over{c_n\sigma_n}}.
\end{equation}

It remains to deal with the deviation
$|\varphi_{n,m}^{[2]}(t)-\e^{-t^2/2}|$. We have
\begin{align*}
\big|\varphi_{n,m}^{[2]}(t)-\e^{-t^2/2}\big|
&=\exp\!\left\{-{{t^2}\over{2}}\right\}
\left|\exp\!\left\{\sum_{k=m+1}^{n-1}{{(n-k)^2}\over{k^2}}\left({{-\i t^3}\over{2\sigma_n^3}}+{{t^4}\over{8\sigma_n^4}}\right)\right\}-1\right|\cr
&=\exp\!\left\{-{{t^2}\over{2}}\right\}
\left|\sum_{j=1}^{\infty}{{1}\over{j!}}
\left[\sum_{k=m+1}^{n-1}{{(n-k)^2}\over{k^2}}
\left({{-it^3}\over{2\sigma_n^3}}+{{t^4}\over{8\sigma_n^4}}\right)\right]^j\right|\cr
&\leq\exp\!\left\{-{{t^2}\over{2}}\right\}
\sum_{j=1}^{\infty}{{1}\over{j!}}\left[\sum_{k=m+1}^{n-1}{{(n-k)^2}\over{k^2}}
\left|{{-it^3}\over{2\sigma_n^3}}+{{t^4}\over{8\sigma_n^4}}\right|\right]^j\cr
&\leq\exp\!\left\{-{{t^2}\over{2}}\right\}s_{n,m}\!(t)
\sum_{j=0}^{\infty}{{s_{n,m}^j\!(t)}\over{j!}},
\end{align*}
where
$$
s_{n,m}\!(t)=\sum_{k=m+1}^{n-1}\!{{(n-k)^2}\over{k^2}}
\left|{{-\i t^3}\over{2\sigma_n^3}}+{{t^4}\over{8\sigma_n^4}}\right|
=\sum_{k=m+1}^{n-1}\!{{(n-k)^2}\over{k^2}}
\sqrt{{{t^6}\over{4\sigma_n^6}}+{{t^8}\over{64\sigma_n^8}}}\,.
$$
We give two different bounds for $s_{n,m}\!(t)$. First, since
$|t|<\sigma_n$, by (\ref{0var}) we can write
$$
s_{n,m}\!(t)={|t|^3\over{2\sigma_n^3}}\sqrt{1+{{t^2}\over{16\sigma_n^2}}}\!\!\sum_{k=m+1}^{n-1}\!\!{{(n-k)^2}\over{k^2}}\leq
{{\sqrt{17}|t|^3}\over{8\sigma_n^3}}\!\!\sum_{k=m+1}^{n-1}\!\!{{n(n-k)}\over{k^2}}
={\sqrt{17}\over8}{|t|^3\over{\sigma_n}}.
$$
For the second bound, we first estimate $(n-k)|t|/(n\sigma_n)$.
Applying (\ref{n_5}), and the fact that the maximum of the function
$\sqrt{x(1-x)}$ is $1/2$ on the interval [0,1], we get
$$
{{n-k}\over{n}}{|t|\over\sigma_n}
\leq{{n-k}\over{n}}{k\over{\sqrt{n(n-k)}}}
=\sqrt{{{n-k}\over{n}}}{k\over{n}}
\leq\sqrt{{{n-k}\over{n}}{k\over{n}}} \leq{1\over2}.
$$
This inequality, along with $|t|<\sigma_n$ and (\ref{0var}), gives
$$
s_{n,m}(t)=
{t^2\over{2\sigma_n^2}}\sqrt{1+{{t^2}\over{16\sigma_n^2}}}\sum_{k=m+1}^{n-1}\!\left[{{n(n-k)}\over{k^2}}
{{n-k}\over{n}}{|t|\over\sigma_n}\right]\leq{\sqrt{17}\over16}t^2.
$$
Continuing, we apply the first bound of $s_{n,m}(t)$ to the
function before the sum and the second one to each term of the sum,
to get
\begin{align}\label{n_16}
\big|\varphi_{n,m}^{[2]}(t)-\e^{-t^2/2}\big|
&\leq{\sqrt{17}\over8}{|t|^3\over{\sigma_n}}
\exp\!\left\{-\left({1\over2}-{\sqrt{17}\over16}\right)t^2\right\}\cr
&\leq{\sqrt{17}\over8}{|t|^3\over{\sigma_n}}\exp\!\left\{-{3\over16}t^2\right\}.
\end{align}

Now, collecting the bounds, we can return to the estimation of the
integrals in (\ref{n_9}). We obtain
$$
{1\over2}\int_{-c_n\sigma_n}^{c_n\sigma_n}\left|{{\varphi_{n,m}(t)}\over{t}}\right|\delta_{n,m}^{[1]}(t)\,\d
t
\leq{{1}\over{\sigma_n}}{{5}\over{48}}\int_{-\infty}^{\infty}\exp\!\left\{-{{20t^2}\over{\pi^4}}\right\}t^2\,\d
t ={{1}\over{\sigma_n}}{{\sqrt{5}\pi^{13/2}}\over{3840}},
$$
by (\ref{n_11}) and (\ref{n_14});
$$
{1\over2}\int_{-c_n\sigma_n}^{c_n\sigma_n}\left|{{\varphi_{n,m}^{[1]}(t)}\over{t}}\right|\delta_{n,m}^{[1]}(t)\,\d
t
\leq{{1}\over{\sigma_n}}{{5}\over{48}}\int_{-\infty}^{\infty}\exp\!\left\{-{{t^2}\over4}\right\}t^2\,\d
t ={{1}\over{\sigma_n}}{{5\sqrt{\pi}}\over{12}},
$$
by (\ref{n_12}) and (\ref{n_14});
\begin{align*}
{1\over2}\int_{-c_n\sigma_n}^{c_n\sigma_n}\left|{{\varphi_{n,m}^{[1]}(t)}\over{t}}\right|\delta_{n,m}^{[2]}(t)\,\d
t
&\leq{{1}\over{c_n\sigma_n}}{{5\sqrt{5}}\over{48(1-c)}}\int_{-\infty}^{\infty}\exp\!\left\{-{{t^2}\over{4}}\right\}t^2\,\d
t\cr &={{1}\over{c_n\sigma_n}}{{5\sqrt{5\pi}}\over{12(1-c)}},
\end{align*}
by (\ref{n_12}) and (\ref{n_15});
\begin{align*}
{1\over2}\!\int_{-c_n\sigma_n}^{c_n\sigma_n}\!\left|{{\varphi_{n,m}^{[2]}(t)}\over{t}}\right|\delta_{n,m}^{[2]}(t)\d
t
&\leq{{1}\over{c_n\sigma_n}}{{5\sqrt{5}}\over{48(1-c)}}\int_{-\infty}^{\infty}\!\exp\!\left\{-{{3t^2}\over8}\right\}t^2\,\d
t\cr &={1\over{c_n\sigma_n}}{{5\sqrt{30\pi}}\over{54(1-c)}},
\end{align*}
by (\ref{n_13}) and (\ref{n_15}); and finally
$$
\int_{-c_n\sigma_n}^{c_n\sigma_n}\left|{{\varphi_{n,m}^{[2]}(t)-\e^{-t^2/2}}\over{t}}\right|\,\d
t\leq
{{1}\over{\sigma_n}}{{\sqrt{17}}\over{8}}\int_{-\infty}^{\infty}\exp\!\left\{-{{3t^2}\over{16}}\right\}t^2\,\d
t ={{1}\over{\sigma_n}}{{4\sqrt{51\pi}}\over9}
$$
by (\ref{n_16}). Since $1/\sigma_n\leq 1/(c_n\sigma_n)$  by the definition
of $c_n$ in (\ref{n_3}), substitution of all these bounds into (\ref{n_9}) yields
$I_{n,m}\leq \overline{C}/(c_n\sigma_n)$, where
$$
\overline{C}={{\sqrt{5}\,\pi^{13/2}}\over{3840}}+
{{5\sqrt{\pi}}\over{12}}+ {{5\sqrt{5\pi}}\over{12(1-c)}}+
{{5\sqrt{30\pi}}\over{54(1-c)}}+ {{4\sqrt{51\pi}}\over9}<
{2.5503\over{1-c}}+7.3566
$$
Writing this back in (\ref{n_4}), we obtain the inequality
$$
\sup_{x\in \R}\big|F_{n,m}(x)-\Phi(x)\big|\leq{{\widetilde{C}}\over{c_n\,\sigma_n}}
$$
with the constant
$$
\widetilde{C}=\left({{2.5503}\over{1-c}}+7.3566\right){b\over{2\pi}}+{c_b\over\sqrt{2\pi}}.
$$ According to Cs\"{o}rg\H{o} \cite{CS}, the minimum of $c_b$ is 4.439 occurring at
$b=1.868\,$, and with these values $\widetilde{C}>0.7583/(1-c)+3.9581$. Also, it is easy to see from the definition of $c_n$ in (\ref{n_3}) that $c_n\geq c\frac{m}{n}$, thus
$$
\sup_{x\in \R}\big|F_{n,m}(x)-\Phi(x)\big|\leq\frac{0.7583/(1-c)+3.95811}{c}\frac{n}{m}\frac{1}{\sigma_n}.
$$
Now minimizing $\frac{0.7583/(1-c)+3.95811}{c}$ over $c\in(0,1)$, we finally obtain the inequality of the theorem with $C=9.257$. $\blacksquare$


\chapter{Poisson approximation}

In the first section of this chapter, we leave the coupon collector's problem and concern Poisson approximation to the distribution
of sums of independent nonnegative integer valued random variables in general.  We
complement the classical Poisson convergence theorem
of Gnedenko and Kolmogorov \cite{GK}, in the setting of triangular arrays, with error
bounds that are expressed in terms of total variation distance, which was defined in Section 2.1. as
\begin{equation}\label{tav}
d_{\mathrm{TV}}(X,Y)=\sup_{A\subset \Z_+}|\p(X\in A)-\p(Y\in A)|,
\end{equation}
for any two random variables $X$ and $Y$ that map into $\Z_+:=\{0,1,\ldots\}$.

For each~$n$, we approximate the distribution of the $n$-th row sum with a Poisson
distribution whose mean~$\lambda_n$ is defined only in terms of the distributions
of the random variables in the~$n$-th row, namely $\lambda_n$ equals the sum of the probabilities $\p(X\neq0)$, where $X$ runs over the random variables of the~$n$-th row. We do not assume the existence of
moments, as is the case in analogous results proved by Barbour and Hall \cite{BH},
and our lower bounds are much simpler in form to theirs, being of precisely
the same form, up to a constant, as our upper bound, provided that the
means~$\lambda_n$ are bounded away from infinity.

We then continue in the second section of the chapter with an application of these results to the coupon collector's problem.
We recall from Section 1.2. that Baum and Billingsley proved in \cite{BB} (using the method of characteristic functions) that if
\begin{equation}\label{m}
m\to\infty \quad\textrm{ and }\quad\frac{n-m}{\sqrt{n}}\to\sqrt{2\lambda}\quad\textrm{ for some }\lambda>0\textrm{ constant, as } n\to\infty,
\end{equation}
then $W_{n,m}-(n-m)$ converges in distribution to the Poisson law with mean $\lambda$.
We express this problem as a special case of the Poisson limit theorem above, and immediately obtain the corresponding Poisson approximation results. An even stronger result can be proved in this special case: due to the combinatorial structure of the problem, one can determine explicitly the first order term in the error of the approximation, and this is what we shall do in Section 5.3.

Finally, in Section 5.4. we give another Poisson approximation result to the coupon collector's waiting time. This time the mean $\lambda_n'$ of the approximating Poisson law is chosen to match the mean of the waiting time. The result is proven with the help of Stein's method.

\section{Poisson approximation in a general Poisson limit theorem}

In \cite{GK} (p.~132) Gnedenko and Kolmogorov give necessary and sufficient conditions for sums of independent infinitesimal random variables to converge to the Poisson law. In case of nonnegative integer valued random variables their limit theorem can be stated as follows.

\bigskip

\begin{theorem}\label{t_GK}{\em{\bf (Gnedenko, Kolmogorov)}}
Let $\{Y_{n1}, Y_{n2}, \ldots, Y_{nr_n}\}_{n\in\N}$ be a triangular array of row-wise independent nonnegative integer valued random variables such that
\begin{align}
&\min_{1\leq k\leq r_n}\p(Y_{nk}=0)\rightarrow1,\quad n\to\infty\label{Tfeltetel1},\\
&\sum_{k=1}^{r_n}\p(Y_{nk}\geq 1)\rightarrow\lambda,\quad \lambda>0 \textrm{ constant},\quad n\to\infty\label{Tfeltetel2}\\
&\sum_{k=1}^{r_n}\p(Y_{nk}\geq2)\to 0,\quad n\to\infty. \label{Tfeltetel3}
\end{align}
Then
$$Y_n:=\sum_{k=1}^{r_n}Y_{nk}\xrightarrow{\cal D}N_{\lambda}$$
as $n\to\infty$, where $N_{\lambda}$ is a Poisson random variable with parameter $\lambda$.
\end{theorem}

\bigskip
We shall refine the obvious approximation of the $Y_n$-s that the limit theorem suggests by approximating the distribution of each of the $Y_n$ random variables not with the limiting Poisson distribution, but with a Poisson distribution that has a suitably chosen parameter that depends on $n$, namely by the distribution of $N_{\lambda_n}\sim \mathrm{Po}(\lambda_n)$, where
$$\lambda_n=\sum_{k=1}^{r_n}\p(Y_{nk}\geq 1).$$

\bigskip

\begin{theorem}\label{t_upper_bound}{\em{\bf (The upper bound.)}}
For any triangular array $\{Y_{n1}, Y_{n2}, \ldots, Y_{nr_n}\}_{n\in\N}$ of row-wise independent nonnegative integer valued random variables
$$d_{\mathrm{TV}}({\cal D}(Y_n),\mathrm{Po}(\lambda_n))
\leq\sum_{k=1}^{r_n}\left[\p(Y_{nk}\geq 2)+\p(Y_{nk}\geq 1)^2\right].$$
\end{theorem}

\bigskip

\noindent\textbf{Proof.}  The proof follows the argument in \cite{BHJ} p.~181. For each $k=1, 2,\ldots,r_n$, $n\in\N$, we define the random variable
$$I_{nk}:=
\left\{
  \begin{array}{ll}
    0, & \hbox{if $Y_{nk}=0$;} \\
    1, & \hbox{if $Y_{nk}\geq1$.}
  \end{array}
\right.
$$
Thus for each $n\in\N$, $I_n:=\sum_{k=1}^{r_n}I_{nk}$ is a sum of independent Bernoulli random variables with success probabilities $q_{nk}:=\p(Y_{nk}\geq1)$, $k=1, 2,\ldots,r_n$. By Le Cam's inequality \cite{Cam}
$$d_{\mathrm{TV}}({\cal D}(I_n), \mathrm{Po}(\lambda_n))\leq\sum_{k=1}^{r_n}q_{nk}^2=\sum_{k=1}^{r_n}\p(Y_{nk}\geq1)^2.$$
Also, for any two random variables $X$ and $Y$ defined on the same probability space the coupling inequality (\ref{coupl_ineq}) says that
$$d_{\mathrm{TV}}({\cal D}(X), {\cal D}(Y))\leq\p(X\neq Y),$$
hence we have
\begin{align*}
d_{\mathrm{TV}}\left({\cal D}(Y_n),{\cal D}(I_n)\right)
&\leq\p\left(\sum_{k=1}^{r_n}Y_{nk}\neq\sum_{k=1}^{r_n}I_{nk}\right)
=\p\left(\cup_{k=1}^{r_n}\{Y_{nk}\neq I_{nk}\}\right)\\
&\leq\sum_{k=1}^{r_n}\p(Y_{nk}\neq I_{nk})
=\sum_{k=1}^{r_n}\p(Y_{nk}\geq2).
\end{align*}
Putting these two bounds together in
$$d_{\mathrm{TV}}({\cal D}(Y_n),\mathrm{Po}(\lambda_n))\leq d_{\mathrm{TV}}\left({\cal D}(Y_n),{\cal D}(I_n)\right) +d_{\mathrm{TV}}({\cal D}(I_n), \mathrm{Po}(\lambda_n)),$$
the assertion of the theorem follows. $\blacksquare$

\bigskip

\begin{theorem}\label{t_lower bound}{\em{\bf(The lower bound.)}}
If $\{Y_{n1}, Y_{n2}, \ldots, Y_{nr_n}\}_{n\in\N}$ is a triangular array of row-wise independent nonnegative integer valued random variables such that\\ $\min_{1\leq k\leq r_n}\p(Y_{nk}=0)\geq \frac{3}{4}$ for all $n\in\N$, then
$$d_{\mathrm{TV}}({\cal D}(Y_n),\mathrm{Po}(\lambda_n))
\geq\frac{1}{10}\left(\prod_{k=1}^{r_n}\p(Y_{nk}=0)\right)
\sum_{k=1}^{r_n}\left[\p(Y_{nk}\geq 2)+\p(Y_{nk}\geq 1)^2\right].$$
\end{theorem}

\bigskip

Before turning to the proof of Theorem \ref{t_lower bound} we prove a simple result we are going to need later on.

\bigskip

\begin{proposition} If $0\leq y_i\leq x_i\leq1$ for all $i=1, 2,\ldots, n$, $n\in\N$, then
$$\left(\prod_{i=1}^{n}y_i\right)\sum_{i=1}^{n}(x_i-y_i)\leq\prod_{i=1}^{n}x_i-\prod_{i=1}^{n}y_i\leq\sum_{i=1}^{n}(x_i-y_i).$$
\end{proposition}

\noindent\textbf{Proof.} Defining $y_0:=1$, we can write the difference of the two products in the form of a telescopic sum, thus
\begin{align*}
\prod_{i=1}^{n}x_i-\prod_{i=1}^{n}y_i
&=\sum_{k=1}^{n}[y_1\cdots y_{k-1} x_{k}\cdots x_n - y_1\cdots y_{k} x_{k+1}\cdots x_n]\\
&=\sum_{k=1}^{n}(x_k-y_k)(y_1\cdots y_{k-1}x_{k+1}\cdots x_{n}).
\end{align*}
Due to our assumption on the $y_i$-s the last expression can be bounded form above and from below by
$$y_1\cdots y_n \sum_{k=1}^{n}(x_k-y_k)\leq\sum_{k=1}^{n}(x_k-y_k)(y_1\cdots y_{k-1}x_{k+1}\cdots x_{n})\leq\sum_{k=1}^{n}(x_k-y_k),$$
and the assertion follows. $\blacksquare$

\bigskip

\noindent\textbf{Proof of Theorem \ref{t_lower bound}.} We introduce the notations $\p(Y_{nk}=0)=p_{nk}$ and $\p(Y_{nk}=1)=(1-p_{nk})\tilde{p}_{nk}$, $k=1,2,\ldots,r_n$, thus $\lambda_n=\sum_{k=1}^{r_n}(1-p_{nk})$, $n\in\N$. We are going to prove the theorem by approximating the following elementary lower bound for the total variation distance of the distributions considered:
\begin{equation}\label{d_TV}
d_{\mathrm{TV}}({\cal D}(Y_n),\mathrm{Po}(\lambda_n))\geq \frac{1}{2}|\p(Y_n=0)-\p(N_{\lambda_n}=0)|+\frac{1}{2}|\p(Y_n=1)-\p(N_{\lambda_n}=1)|,
\end{equation}
which can be justified by taking $A=\{0\}$ and $A=\{1\}$ in (\ref{tav}).

We start by bounding the difference of the point probabilities at 0. Since
$$\p(N_{\lambda_n}=0)=\e^{-\lambda_n}=\prod_{k=1}^{r_n}\e^{-\left(1-p_{nk}\right)},$$
$$\p(Y_n=0)=\prod_{k=1}^{r_n}p_{nk},$$
and
$\e^{-\left(1-p_{nk}\right)}\geq p_{nk}$ for all $k=1, 2,\ldots,r_n$, $n\in\N$,
applying the Proposition above yields
\begin{equation*}
|\p(N_{\lambda_n}=0)-P(Y_n=0)|\geq\left(\prod_{k=1}^{r_n}p_{nk}\right)\sum_{k=1}^{r_n}\left[\e^{-\left(1-p_{nk}\right)}-p_{nk}\right].
\end{equation*}
Since $1-p_{nk}\leq 1$, we have
\begin{equation*}
\e^{-\left(1-p_{nk}\right)}
\geq 1-\left(1-p_{nk}\right)+\frac{1}{2}\left(1-p_{nk}\right)^2-\frac{1}{6}\left(1-p_{nk}\right)^3
\geq p_{nk}+\frac{1}{3}\left(1-p_{nk}\right)^2
\end{equation*}
for $k=1, 2,\ldots,n$, $n\in\N$, which yields
\begin{equation}\label{point_prob_0}
|\p(N_{\lambda_n}=0)-P(Y_n=0)|\geq\frac{1}{3}\left(\prod_{k=1}^{r_n}p_{nk}\right)\sum_{k=1}^{r_n}\left(1-p_{nk}\right)^2.
\end{equation}

This inequality implies the assertion of Theorem \ref{t_lower bound} in the case when $3\sum_{k=1}^{r_n}\left(1-p_{nk}\right)^2\!\\ \geq2\sum_{k=1}^{r_n}\left(1-p_{nk}\right)(1-\tilde{p}_{nk})$, because we can bound 3/5th of the sum in the display above using this assumption. In fact in this case we obtain a better bound than the one we aimed at. Otherwise, if $2\sum_{k=1}^{r_n}\left(1-p_{nk}\right)(1-\tilde{p}_{nk})\geq3\sum_{k=1}^{r_n}\left(1-p_{nk}\right)^2$, we need to examine the point probabilities at 1 too to improve our current bound.

We have
$$\p(N_{\lambda_n}=1)=\lambda_n \e^{-\lambda_n}=\sum_{k=1}^{r_n}(1-p_{nk})\exp\left\{-\sum_{k=1}^{r_n}(1-p_{nk})\right\},$$
and since for an arbitrary $n\in\N$ $Y_n=1$ iff for $k=1, 2,\ldots, r_n$ exactly one of the $Y_{nk}$-s takes on 1 and the rest take on 0,
$$\p(Y_n=1)=\left(\prod_{k=1}^{r_n}p_{nk}\right)\sum_{k=1}^{r_n}\frac{(1-p_{nk})\tilde{p}_{nk}}{p_{nk}}.$$
Some elementary algebra gives
\begin{multline*}
\p(N_{\lambda_n}=1)-\p(Y_n=1)
=\prod_{k=1}^{r_n}p_{nk}\left(\sum_{k=1}^{r_n}\left(1-p_{nk}\right)(1-\tilde{p}_{nk})-\sum_{k=1}^{r_n}\frac{(1-p_{nk})^2\tilde{p}_{nk}}{p_{nk}}\right)+\\
+\left(\exp\left\{-\sum_{k=1}^{r_n}(1-p_{nk})\right\}-\prod_{k=1}^{r_n}p_{nk}\right)\sum_{k=1}^{r_n}(1-p_{nk}),
\end{multline*}
where in the second term we recognize the point probabilities at 0. Using $\displaystyle\frac{\tilde{p}_{nk}}{p_{nk}}\leq\frac{1}{\displaystyle\min_{1\leq k\leq r_n}p_{nk}}$ and the fact that the difference of the 0 probabilities in the formula above is always positive we obtain
$$
\p(N_{\lambda_n}=1)-\p(Y_n=1)
\geq\prod_{k=1}^{r_n}p_{nk}\left(\sum_{k=1}^{r_n}\left(1-p_{nk}\right)(1-\tilde{p}_{nk})-\frac{1}{\displaystyle\min_{1\leq k\leq r_n}p_{nk}}\sum_{k=1}^{r_n}(1-p_{nk})^2\right).
$$
From this by inequality (\ref{point_prob_0}) we obtain
\begin{multline*}
\p(N_{\lambda_n}=0)-\p(Y_n=0)+\p(N_{\lambda_n}=1)-\p(Y_n=1)\geq\\
\left(\prod_{k=1}^{n}p_{nk}\right)\left(\sum_{k=1}^{r_n}\left(1-p_{nk}\right)(1-\tilde{p}_{nk})+\left[\frac{1}{3}-\frac{1}{\displaystyle\min_{1\leq k\leq r_n}p_{nk}}\right]\sum_{k=1}^{r_n}(1-p_{nk})^2\right).
\end{multline*}
Now $\frac{1}{3}-\frac{1}{\displaystyle\min_{1\leq k\leq r_n}p_{nk}}\geq -1$ in the range of $n$ for which the assumption $\min_{1\leq k\leq r_n}p_{nk}\geq \frac{3}{4}$ of the Theorem holds, thus
\begin{multline*}
\p(N_{\lambda_n}=0)-\p(Y_n=0)+\p(N_{\lambda_n}=1)-\p(Y_n=1)\geq\\ \left(\prod_{k=1}^{n}p_{nk}\right)\left(\sum_{k=1}^{r_n}\left(1-p_{nk}\right)(1-\tilde{p}_{nk})-\sum_{k=1}^{r_n}(1-p_{nk})^2\right),
\end{multline*}
and it can be seen that the latter bound is at most
$$\frac{1}{5}\left(\prod_{k=1}^{n}p_{nk}\right)\left(\sum_{k=1}^{r_n}\left(1-p_{nk}\right)(1-\tilde{p}_{nk})+\sum_{k=1}^{r_n}(1-p_{nk})^2\right)$$
for all $n$ such that $2\sum_{k=1}^{r_n}\left(1-p_{nk}\right)(1-\tilde{p}_{nk})\geq3\sum_{k=1}^{r_n}\left(1-p_{nk}\right)^2$. This together with (\ref{d_TV}) proves the Theorem. $\blacksquare$

\bigskip

Theorems \ref{t_lower bound} and \ref{t_upper_bound} together state that the order of the error of our Poisson approximation for the random variables in Theorem \ref{t_GK} is
$$\sum_{k=1}^{r_n}\left[\p(Y_{nk}\notin\{0,1\})+\p(Y_{nk}\geq 1)^2\right].$$
Barbour and Hall have proved similar results in \cite{BH} using Stein's method: they approximate a sum $\sum_{j=1}^{n}Y_j$ of independent nonnegative integer valued random variables with a Poisson variable that has mean $\sum_{j=1}^{n}\p(Y_j=1)$ or $\sum_{j=1}^{n}\E(Y_j)$. (Note that the parameter of our approximating Poisson random variable is between these two values.) Their bounds are expressed differently, and involve second moments of the random variables~$Y_j$. Moreover, their lower bounds would yield no useful information at all in the application to be considered in the next section.

We also obtain the following result.

\bigskip

\begin{corollary} For the rate of convergence in Theorem \ref{t_GK} we have the upper bound
$$d_{\mathrm{TV}}({\cal D}(Y_n),\mathrm{Po}(\lambda))
\leq\sum_{k=1}^{r_n}\left[\p(Y_{nk}\geq2)+\p(Y_{nk}\geq 1)^2\right]
+\left|\sum_{k=1}^{r_n}\p(Y_{nk}\geq 1)-\lambda\right|,\quad n\in\N.$$
\end{corollary}

\bigskip

\noindent\textbf{Proof.} Since
$$d_{\mathrm{TV}}({\cal D}(Y_n),\mathrm{Po}(\lambda))\leq
d_{\mathrm{TV}}({\cal D}(Y_n),\mathrm{Po}(\lambda_n))+d_{\mathrm{TV}}(\mathrm{Po}(\lambda_n), \mathrm{Po}(\lambda)),$$
the assertion follows from Theorem \ref{t_upper_bound} and because for any $N_{\nu_1}\sim \textrm{Poisson}(\nu_1)$ and $N_{\nu_2}\sim \textrm{Poisson}(\nu_2)$, where $0<\nu_1<\nu_2$, we have
$$d_{\mathrm{TV}}({\cal D}(N_{\nu_1}),{\cal D}(N_{\nu_2}))\leq \min\left\{1,\nu_2^{-1/2}\right\}(\nu_2-\nu_1).$$
For reference see for example Remark 1.1.4. in \cite{BHJ}. $\blacksquare$

\bigskip
\bigskip


\section{Coupon collecting with an approximately Poisson \\distributed waiting time -- application of the general results}

We begin this section by examining how the coupon collector's problem defined in the introduction fits in the framework of the previous section. The equality in distribution in (\ref{Wgeo}) can be reformulated for $\widetilde{W}_{n,m}:=W_{n,m}-(n-m)$ as
\begin{equation}\label{hullamos}
\widetilde{W}_{n,m}\;\eq{\cal D}\;\sum_{i=m+1}^{n}\widetilde{X}_{n,i},
\end{equation}
where the $\widetilde{X}_{n,i}$, $i=m+1,\ldots,n$, random variables are independent, and $\widetilde{X}_{n,i}+1$ has geometric distributions with success probability $i/n$, $i\in\{m+1,\ldots, n\}$, $n\in\N$. The triangular array $\{\widetilde{X}_{n,m+1},\ldots,\widetilde{X}_{n,n}\}_{n\in\N}$ satisfies the conditions of Theorem \ref{t_GK}: the variables of the array are infinitesimal, i.e.~they satisfy condition (\ref{Tfeltetel1}): for any $0<\varepsilon<1$
$$\max_{m+1\leq i\leq n}\p(\widetilde{X}_{n,i}>\varepsilon)
=\left[1-\min_{m+1\leq i\leq n}\p(\widetilde{X}_{n,i}=0)\right]
=\left[1-\min_{m+1\leq i\leq n}\frac{i}{n}\right]
=\frac{n-m+1}{n}\to 0,$$
by (\ref{m}); and according to (\ref{lambda_n}) and (\ref{lambda_n,j}) in the proposition below, they also satisfy conditions (\ref{Tfeltetel2}) and (\ref{Tfeltetel3}).

\bigskip

\begin{proposition}\label{p_lambdak}
If $\{m=m(n)\}_{n\in\N}$ is a sequence of integers that satisfies (\ref{m}), then
\begin{align}
&\lambda_n=\lambda_{n,1}:=\sum_{i=m+1}^{n}\left(1-\frac{i}{n}\right)\to\lambda,\quad\textrm{ and}\label{lambda_n}\\
&\lambda_{n,j}:=\sum_{i=m+1}^{n}\left(1-\frac{i}{n}\right)^j\leq\lambda_n\left(\frac{2\lambda_n}{n}\right)^{\frac{j-1}{2}},
\quad\textrm{and}\quad\lambda_{n,j}\to0,\quad j=2,3,\ldots\label{lambda_n,j}
\end{align}
\end{proposition}

\bigskip

\noindent\textbf{Proof.} (\ref{lambda_n}) is true, because
\begin{equation*}
\lambda_n=\sum_{i=m+1}^{n}\left(1-\frac{i}{n}\right)
=n-m-\frac{1}{n}\left[\frac{n(n+1)}{2}-\frac{m(m+1)}{2}\right]
=\frac{(n-m)(n-m-1)}{2n}\to\lambda
\end{equation*}
by (\ref{m}). By taking the square root of both sides of the equality above it can be deduced that
\begin{equation}\label{csillag}
\frac{n-m-1}{\sqrt{n}}\leq\sqrt{2\lambda_n}.
\end{equation}
Now we prove the first assertion of (\ref{lambda_n,j}) by induction. For an arbitrary $j=2,3,\ldots$ we bound $\lambda_{n,j}$ as follows:
$$\lambda_{n,j}=\sum_{i=m+1}^{n}\left(1-\frac{i}{n}\right)^j
\leq\frac{n-m-1}{n}\sum_{i=m+1}^{n}\left(1-\frac{i}{n}\right)^{j-1}
=\frac{n-m-1}{n}\lambda_{n,j-1}
$$
Since for $j=2$ this gives $\lambda_{n,2}\leq\lambda_n\sqrt{\frac{2\lambda_n}{n}}$ by (\ref{csillag}), we have the first part of (\ref{lambda_n,j}) in this case. If we have the same result for some $j>2$, then it holds true for $j+1$ as well by the argument above, (\ref{csillag}) and the inductional hypothesis. Since $\lambda_n\to\lambda$ by (\ref{lambda_n}), the second part of (\ref{lambda_n,j}) follows from the first. $\blacksquare$

\bigskip

Thus we see that the limit theorem proved by Baum and Billingsley \cite{BB} concerning the coupon collector's problem is a special case of the Gnedeno--Kolmogorov theorem. If we apply the results of the previous section to $\widetilde{W}_{n,m}$, we obtain the following.

\bigskip

\begin{corollary}\label{kov_p1}
If $\{m=m(n)\}_{n\in\N}$ is a sequence of integers that satisfies (\ref{m}), then the error of the approximation of the coupon collector's $\widetilde{W}_{n,m}$ waiting time with the Poisson random variable $N_{\lambda_n}$, that has mean $\lambda_n=\sum_{i=m+1}^{n}\left(1-\frac{i}{n}\right)$, is of order $\sum_{i=m+1}^{n}\left(1-\frac{i}{n}\right)^2$. In fact, for all $n$ such that $\min_{m+1\leq i\leq n}\frac{i}{n}\geq \frac{3}{4}$,
\begin{equation*}
\frac{1}{5}\left(\prod_{i=m+1}^{n}\frac{i}{n}\right)\sum_{i=m+1}^{n}\left(1-\frac{i}{n}\right)^2
\leq d_{\mathrm{TV}}({\cal D}(\widetilde{W}_{n,m}),{\cal D}(N_{\lambda_n}))
\leq 2\sum_{i=m+1}^{n}\left(1-\frac{i}{n}\right)^2.
\end{equation*}
\end{corollary}

\bigskip

\begin{corollary}\label{kov_p2}
For the rate of convergence in the Poisson limit theorem concerning the coupon collector's problem we have the upper bound
\begin{equation*}
d_{\mathrm{TV}}({\cal D}(\widetilde{W}_{n,m}),{\cal D}(N_{\lambda}))\leq2\sum_{i=m+1}^{n}\left(1-\frac{i}{n}\right)^2
+\left|\sum_{i=m+1}^{n}\left(1-\frac{i}{n}\right)-\lambda\right|.
\end{equation*}
\end{corollary}


\section{Coupon collecting with an approximately Poisson \\distributed waiting time -- combinatorial approach}

Now we undertake a combinatorial approach to the coupon collector's problem, which will yield us a stronger result than the one of Corollary \ref{kov_p1}. Namely, we shall derive the first asymptotic correction of the $\p(\widetilde{W}_{n,m}=k)$, $k=0,1,\ldots$, probabilities to the corresponding Poisson point probabilities. We state the result in the following theorem. We note that in principal the method presented in the proof can be extended to determine higher order terms in the asymptotic expansion.

\bigskip

\begin{theorem}\label{t_poisson_sorfejtes}
If $\{m=m(n)\}_{n\in\N}$ is a sequence of nonnegative integers that satisfies (\ref{m}) and $\lambda_n$ and $\lambda_{n,2}$ are defined as in (\ref{lambda_n}) and (\ref{lambda_n,j}), then
\begin{align*}
&\p(\widetilde{W}_{n,m}=0)=\e^{-\lambda_n}
-\e^{-\lambda_n}\frac{\lambda_{n,2}}{2}
+O\!\left(\frac{1}{n}\right),\\
&\p(\widetilde{W}_{n,m}=1)=\e^{-\lambda_n}\lambda_n
-\e^{-\lambda_n}\lambda_n\frac{\lambda_{n,2}}{2}
+O\!\left(\frac{1}{n}\right),\\
&\p(\widetilde{W}_{n,m}=k)
=\e^{-\lambda_n}\frac{\lambda_n^k}{k!}
+\e^{-\lambda_n}\left(\frac{\lambda_n^{k-2}}{(k-2)!}-\frac{\lambda_n^{k}}{k!}\right)\frac{\lambda_{n,2}}{2}
+O\!\left(\frac{1}{n}\right),\quad k\geq2.
\end{align*}
\end{theorem}

\bigskip

We note that $\lambda_{n,2}=\frac{(2\lambda_n)^{3/2}}{3\sqrt{n}} + O\!\left(\frac{1}{n}\right)$. Indeed,
\begin{align*}
\lambda_{n,2}
&=\sum_{i=m+1}^{n}\left(1-\frac{i}{n}\right)^2\\
&=n-m-\frac{2}{n}\left[\frac{n(n+1)}{2}-\frac{m(m+1)}{2}\right]+\frac{1}{n^2}\left[\frac{n(n+1)(2n+1)}{2}-\frac{m(m+1)(2m+1)}{2}\right]\\
&=\frac{(n-m)(n-m-1)(n-m-\frac{1}{2})}{3n^2}\\
&=\frac{(2\lambda_n)^{3/2}}{3\sqrt{n}}+\left(\frac{(n-m)(n-m-1)}{3n^2}\left[n-m-\frac{1}{2}-\sqrt{(n-m)(n-m-1)}\right]\right),
\end{align*}
where we used the fact that $\lambda_n=\frac{(n-m)(n-m-1)}{2n}$, which we calculated in the proof of Proposition \ref{p_lambdak}, and the second term in the formula above is $O\!\left(\frac{1}{n}\right)$ by (\ref{m}).

\bigskip

\noindent\textbf{Proof of Theorem \ref{t_poisson_sorfejtes}.} We are going to represent each possible outcome of the collector's sampling with a sequence of integers the following way: let us suppose that while sampling (with replacement), the collector labels the distinct coupons he draws form 1 to $n-m$ in the order he obtains them in the course of time, and after each draw he writes down the label of the coupon just drawn. So he begins the enumeration of labels with a 1 after the first draw, and each number that he writes to the end of his list after a draw is either the label already on the coupon he just got (if he had drawn the same one before), or it is the label he gives the coupon at that moment, which would be the smallest positive integer he has not yet used in the process of sampling and labeling. In the first case we call the new member of the sequence "superfluous", while in the second case we call it a "first appearance".

We fix an arbitrary $k\in\N$, and we suppose that $n$ so big that $n-m>k$ holds. Now $\widetilde{W}_{n,m}=k$ means that the collector had $k$ "superfluous" draws, thus the corresponding representing sequence contains $n-m$ "first appearances" and $k$ "superfluous" members. We categorize all such outcomes according to how the $k$ "superfluous" draws are split into blocks by the $n-m$ "first appearances" in the representing sequences: to each vector $\underline{k}=(k_{m+1},k_{m+2},\ldots,k_{n-1})$, where $k_i\in\Z_+$, $i=m+1,\ldots,n-1$, and $\sum_{i=m+1}^{n-1}k_i=k$, correspond the sequences where there are $k_{n-1}$ "superfluous" members between the 1st and 2nd "first appearances", $k_{n-2}$ "superfluous" members between the 2nd and 3rd "first appearances", and so on, $k_{m+1}$ "superfluous" members between the $(n-m-1)$th and $(n-m)$th "first appearances". (This is the same as saying that $\widetilde{X}_{ni}=k_i$, for all $i=m+1,\ldots,n$.) The probability of getting such a sequence is
\begin{multline*}
\frac{n}{n}\left(1-\frac{n-1}{n}\right)^{k_{n-1}}
\frac{n-1}{n}\left(1-\frac{n-2}{n}\right)^{k_{n-2}}\cdots
\left(1-\frac{m+1}{n}\right)^{k_{m+1}}\frac{m+1}{n}\\
=\left(\prod_{i=m+1}^{n}\frac{i}{n}\right)\prod_{i=m+1}^{n-1}\left(1-\frac{i}{n}\right)^{k_i}.
\end{multline*}
It follows that
\begin{equation}\label{p(w=k)}
\p(\widetilde{W}_{n,m}=k)=\left(\prod_{i=m+1}^{n}\frac{i}{n}\right)\sum_{\underline{k}\in I_k}\prod_{i=m+1}^{n-1}\left(1-\frac{i}{n}\right)^{k_i},
\end{equation}
where
$$I_k:=\left\{\underline{k}\in \Z_+^{n-m-1} : \sum_{i=m+1}^{n-1}k_i=k\right\}.$$

\bigskip

Now we are going to examine the sum in (\ref{p(w=k)}) above, which we denote by $S_{n,m,k}=S_k$. For $k=0$ it is an empty sum, and thus it equals 1 by definition. Now let us suppose that $k>2$, we are going to return to the cases $k=0$ and 1 later on. For an arbitrary such $k$ we see that
$$I_k=\cup_{l=1}^{k} I_{k,l},\quad\textrm{ where }\quad
I_{k,l}=\{\underline{k}\in I_k : \underline{k} \textrm{ has exactly } l \textrm{ nonzero components}\}, l=1,\ldots,k,$$
and we correspondingly define $S_{k,l}$ to be the part of $S_{k}$ that contains the summands over $\underline{k}\in I_{k,l}$, thus we have
\begin{equation}\label{s_nm}
S_{k}=\sum_{\underline{k}\in I_k}\prod_{i=m+1}^{n-1}\left(1-\frac{i}{n}\right)^{k_i}
=\sum_{l=1}^{k}\sum_{\underline{k}\in I_{k,l}}\prod_{i=m+1}^{n-1}\left(1-\frac{i}{n}\right)^{k_i}
=\sum_{l=1}^{k}S_{k,l}.
\end{equation}

To determine the limit of $S_k$ we examine the asymptotic behavior of the $S_{k,l}$ expressions separately. We fix an arbitrary $l=1,\ldots,k$, and with $|A|$ denoting the cardinality of an arbitrary set $A$, we now calculate $|I_{k,l}|$. We can think of the vectors in $I_k$ as the results of distributing $k$ 1-s in $n-m-1$ spaces in all possible ways: to each of these distributions correspond a vector in $I_k$ whose $i$th component is the number of 1-s put in the $i$th space, $i=m+1,\ldots,n$. To produce a vector in $I_{k,l}$ we first choose $l$ different spaces, and we put a 1 in each of them, then we distribute the remaining $k-l$ 1-s in these previously chosen $l$ spaces that already have a 1, but this time any such space can be chosen more than once. This gives
$$|I_{k,l}|={n-m-1\choose l}{k-1\choose k-l},\quad l=1,\ldots,k.$$
We obviously bound $S_{k,l}$ from above if we replace each of the factors in its products by the largest one of them, namely by $1-\frac{m+1}{n}$. This together with the just calculated formula gives
\begin{equation*}
S_{k,l}
\leq{n-m-1\choose l}{k-1\choose k-l}\left(1-\frac{m+1}{n}\right)^k\\
\leq(k-1)!\left(\frac{n-m-1}{\sqrt{n}}\right)^{k+l}\left(\frac{1}{\sqrt{n}}\right)^{k-l}.
\end{equation*}
Hence by (\ref{csillag}) we have
\begin{equation}\label{s_maradek}
S_{k,l}\leq\frac{(k-1)!}{l!(l-1)!}\sqrt{2\lambda_n}^{k+l}\left(\frac{1}{\sqrt{n}}\right)^{k-l}\quad\textrm{and}\quad
\sum_{l=1}^{l'}S_{k,l}\leq k!\min\left\{1,(2\lambda_n)^{k}\right\}\left(\frac{1}{\sqrt{n}}\right)^{k-l'}
\end{equation}
for any $l'\in\{1,\ldots,k\}$. We see from the first inequality that $S_{k,l}$ goes to 0 for $l=1,\ldots,k-1$, but it gives a constant upper bound for $l=k$. We are going to examine the latter case more carefully. Notice that the components of a vector in $I_{k,k}$ are all 0-s and 1-s, thus for any $\underline{k}\in I_{k,k}$ $\frac{1}{k_{m+1}!k_{m+2}!\dots k_{n-1}!}=1$. Using this and the decomposition of the index set $I_k=\cup_{l=1}^{k}I_{k,l}$ we obtain
\begin{multline*}
S_{k,k}=\frac{1}{k!}\sum_{\underline{k}\in I_k}\frac{k!}{k_{m+1}!k_{m+2}!\dots k_{n-1}!}\prod_{i=m+1}^{n-1}\left(1-\frac{i}{n}\right)^{k_i}-\\
-\sum_{l=1}^{k-1}\sum_{\underline{k}\in I_{k,l}}\frac{1}{k_{m+1}!k_{m+2}!\dots k_{n-1}!}\prod_{i=m+1}^{n-1}\left(1-\frac{i}{n}\right)^{k_i}.
\end{multline*}
The first term of $S_{k,k}$ is equal to $\frac{1}{k!}\left[\sum_{i=m+1}^{n}\left(1-\frac{i}{n}\right)\right]^k$ by the polynomial theorem, thus we have
\begin{equation}\label{s_fotag}
S_{k,k}=\frac{\lambda_n^k}{k!}
-\sum_{l=1}^{k-1}\sum_{\underline{k}\in I_{k,l}}\frac{1}{k_{m+1}!k_{m+2}!\dots k_{n-1}!}\prod_{i=m+1}^{n-1}\left(1-\frac{i}{n}\right)^{k_i}.
\end{equation}
It follows that $\lim_{n\to\infty}S_{k,k}=\frac{\lambda^k}{k!}$, because we have (\ref{lambda_n}), and the sum above can be bounded by $\sum_{l=1}^{k-1}S_{k,l}$, which goes to 0 by (\ref{s_maradek}). Thus putting together our results for the expressions $S_{k,l}$ in (\ref{s_nm}), we conclude that the part of $S_k$ that counts -- in the sense that it asymptotically contributes a positive constant to $S_k$ --, is $S_{k,k}$, which is the part of the sum in the defining formula of $S_k$ that corresponds to the 0 - 1 vectors of the $I_k$ index set.

If we write (\ref{s_fotag}) into (\ref{s_nm}), we obtain the following formula for $S_k$:
\begin{equation}\label{s_k}
S_{k}=\frac{\lambda_n^k}{k!}+\sum_{l=1}^{k-1}R_{k,l},
\end{equation}
where
$$R_{k,l}=\sum_{\underline{k}\in I_{k,l}}\left(1-\frac{1}{k_{m+1}!k_{m+2}!\dots k_{n-1}!}\right)\prod_{i=m+1}^{n-1}\left(1-\frac{i}{n}\right)^{k_i}.$$
Our aim is to determine the first order term of the error when we approximate $S_k$ by $\frac{\lambda_n^k}{k!}$. Since $R_{k,l}\leq S_{k,l}$ for each $l=1,\ldots,k-1$, and for the latter expressions we have the bounds of (\ref{s_maradek}), we see that $\sum_{l=1}^{k-1}R_{k,l}=O\left(\frac{1}{\sqrt{n}}\right)$, and the same, but more detailed argument also gives
\begin{equation}\label{maradek}
\sum_{l=1}^{k-2}R_{k,l}\leq\sum_{l=1}^{k-2}S_{k,l}\leq k!\min\left\{1,(2\lambda_n)^{k}\right\}\frac{1}{n}.
\end{equation}
Thus the leading term of the error $\left|S_k-\frac{\lambda_n^k}{k!}\right|$ is of order $\frac{1}{\sqrt{n}}$, and it comes from the term $R_{k,k-1}$.

Before examining $R_{k,k-1}$ we introduce some notations for further use. As an analogue of the set $I_{k,l}$ we define $I_{k-2,l}$ to be the set of vectors $\underline{k}\in \Z_+^{n-m-1}$ such that $\sum_{i=m+1}^{n-1}k_i=k-2$ and $\underline{k}$ has exactly $l$ nonzero components, $l=1,\ldots,k-2$. Also, as an analogue of the expressions $S_{k,l}$ and $S_k$ we define $S_{k-2,l}$ and $S_{k-2}$ by the formulas in (\ref{s_nm}) with $k$ replaced by $k-2$. Finally we introduce
$$I_{k-2,k-2}^j=\left\{\underline{k}\in I_{k-2,k-2}: k_j=0\right\},\quad j=m+1,\ldots,n.$$

We now return to $R_{k,k-1}$. The corresponding index set $I_{k,k-1}$ contains vectors that have exactly one component equal to 2, $k-2$ components equal to 1, and the rest 0. Thus we have
$$R_{k,k-1}
=\frac{1}{2}\sum_{\underline{k}\in I_{k,k-1}}\prod_{i=m+1}^{n-1}\left(1-\frac{i}{n}\right)^{k_i}$$
We can write $R_{k,k-1}$ in another form, if we first sum according to the component of the vectors in $I_{k,k-1}$ which equals 2:
\begin{align*}
R_{k,k-1}
&=\frac{1}{2}\sum_{j=m+1}^{n}\left(1-\frac{j}{n}\right)^2
\left(\sum_{\underline{k}\in I_{k-2,k-2}^j}\prod_{i=m+1}^{n-1}\left(1-\frac{i}{n}\right)^{k_i}\right)\\
=\frac{1}{2}&\sum_{j=m+1}^{n}\left(1-\frac{j}{n}\right)^2
\left[\sum_{\underline{k}\in I_{k-2,k-2}}\prod_{i=m+1}^{n-1}\left(1-\frac{i}{n}\right)^{k_i}
-\displaystyle{\sum_{\underline{k}\in I_{k-2,k-2}\backslash I_{k-2,k-2}^j}}\prod_{i=m+1}^{n-1}\left(1-\frac{i}{n}\right)^{k_i}\right]
\end{align*}
We recognize $S_{k-2,k-2}$ in the first sum in the brackets, thus we can replace it by the formula in (\ref{s_fotag}) with $k-2$ in the place of $k$. As for the second sum in the brackets, we see that $k_j=1$, so there is a $1-\frac{j}{n}$  factor in each of the products, which we can bring before the brackets. These considerations lead to
\begin{align}
R_{k,k-1}
=&\frac{1}{2}\sum_{j=m+1}^{n}\left(1-\frac{j}{n}\right)^2\frac{\lambda_n^{k-2}}{(k-2)!}\nonumber\\
&-\frac{1}{2}\sum_{j=m+1}^{n}\left(1-\frac{j}{n}\right)^2\sum_{l=1}^{k-3}\sum_{\underline{k}\in I_{k-2,l}}\frac{1}{k_{m+1}!k_{m+2}!\dots k_{n-1}!}\prod_{i=m+1}^{n-1}\left(1-\frac{i}{n}\right)^{k_i}\nonumber\\
&-\frac{1}{2}\sum_{j=m+1}^{n}\left(1-\frac{j}{n}\right)^3
\displaystyle{\sum_{\underline{k}\in I_{k-2,k-2}\backslash I_{k-2,k-2}^j}}\prod_{i=m+1,i\neq j}^{n-1}\left(1-\frac{i}{n}\right)^{k_i}\nonumber\\
=&:\frac{1}{2}\sum_{j=m+1}^{n}\left(1-\frac{j}{n}\right)^2\frac{\lambda_n^{k-2}}{(k-2)!}-R_{k,k-1}^1-R_{k,k-1}^2\label{R_k,k-1}
\end{align}
Now we bound the last two expressions. First,
\begin{equation}\label{R_k,k-1,1}
0\leq R_{k,k-1}^1
\leq\frac{1}{2}\sum_{j=m+1}^{n}\left(1-\frac{j}{n}\right)^2\sum_{l=1}^{k-3}S_{k-2,l}
\leq\frac{\lambda_n^{3/2}(k-2)!\min\left\{1,(2\lambda_n)^{k-2}\right\}}{\sqrt{2}}\frac{1}{n}
\end{equation}
by (\ref{lambda_n,j}) and the second inequality in (\ref{s_maradek}) with $k$ replaced by $k-2$. Next,
\begin{equation}\label{R_k,k-1,2}
0\leq R_{k,k-1}^2
\leq\frac{n-m-1}{2n}\sum_{j=m+1}^{n}\left(1-\frac{j}{n}\right)^2S_{k-2,k-2}
\leq\frac{2^{k-2}\lambda_n^k}{(k-2)!}\frac{1}{n}
\end{equation}
by (\ref{csillag}), (\ref{lambda_n,j}) and the first inequality in (\ref{s_maradek}) with $k$ replaced by $k-2$ and $l=k-2$.

We conclude that if we write (\ref{R_k,k-1}) into (\ref{s_k}), we obtain
\begin{equation*}
S_k=\frac{\lambda_n^k}{k!}
+\frac{1}{2}\sum_{i=m+1}^{n}\left(1-\frac{i}{n}\right)^2\frac{\lambda_n^{k-2}}{(k-2)!}
+R_{k,k-1}^1+R_{k,k-1}^2+\sum_{l=1}^{k-2}R_{k,l},
\end{equation*}
where $R_{k,k-1}^1+R_{k,k-1}^2+\sum_{l=1}^{k-2}R_{k,l}=O\left(\frac{1}{n}\right)$ by (\ref{R_k,k-1,1}), (\ref{R_k,k-1,2}), (\ref{maradek}) and the fact that $\lambda_n\to\lambda$ by (\ref{lambda_n}). Thus
\begin{equation}\label{egy}
S_k=\frac{\lambda_n^k}{k!}
+\frac{1}{2}\sum_{i=m+1}^{n}\left(1-\frac{i}{n}\right)^2\frac{\lambda_n^{k-2}}{(k-2)!}.
+O\left(\frac{1}{n}\right)
\end{equation}

\bigskip

Now we return to (\ref{p(w=k)}), and approximate the product $\prod_{i=m+1}^{n}\frac{i}{n}$ in it by $\e^{-\lambda_n}$. Using the definition of $\lambda_n$ in (\ref{lambda_n}) and the expansion formula of the logarithm function the error of the approximation can be written in the form
\begin{align*}
\e^{-\lambda_n}-\prod_{i=m+1}^{n}\frac{i}{n}
&=\exp\left\{-\sum_{i=m+1}^{n}\left(1-\frac{i}{n}\right)\right\}
-\exp\left\{\sum_{i=m+1}^{n}\log\left[1-\left(1-\frac{i}{n}\right)\right]\right\}\nonumber \\
&=\e^{-\lambda_n}\left(1-\exp\left\{-\sum_{j=2}^{\infty}\frac{1}{j}\lambda_{n,j}\right\}\right)\\
=\e^{-\lambda_n}&\left(\frac{1}{2}\sum_{i=m+1}^{n}\left(1-\frac{i}{n}\right)^2+\sum_{j=3}^{\infty}\frac{1}{j}\lambda_{n,j}
-\left[\exp\left\{-\sum_{j=2}^{\infty}\frac{1}{j}\lambda_{n,j}\right\}
-1+\sum_{j=2}^{\infty}\frac{1}{j}\lambda_{n,j}\right]\right),
\end{align*}
where the expressions $\lambda_{n,j}$ are defined as in (\ref{lambda_n,j}). Thus we have
\begin{equation}\label{hiba_e}
\e^{-\lambda_n}-\prod_{i=m+1}^{n}\frac{i}{n}
=\e^{-\lambda_n}\frac{1}{2}\sum_{i=m+1}^{n}\left(1-\frac{i}{n}\right)^2+R_n,
\end{equation}
where
$$R_n=\e^{-\lambda_n}\left(\sum_{j=3}^{\infty}\frac{1}{j}\lambda_{n,j}
-\left[\exp\left\{-\sum_{j=2}^{\infty}\frac{1}{j}\lambda_{n,j}\right\}
-1+\sum_{j=2}^{\infty}\frac{1}{j}\lambda_{n,j}\right]\right),$$
and we are going to show that $R_n=O\left(\frac{1}{n}\right)$.

We are going to bound the sum in the exponent in $R_n$. Since $\lambda_n\to\lambda$ by (\ref{lambda_n}), there exists a threshold number $n_0$ such that for all $n\geq n_0$ we have $\sqrt{\frac{2\lambda_n}{n}}<\frac{1}{2}$. This with inequality (\ref{lambda_n,j}) yields
\begin{equation}\label{lambda_n,j_becslese}
\sum_{j=j_0}^{\infty}\frac{1}{j}\lambda_{n,j}
\leq\lambda_n\left(\frac{2\lambda_n}{n}\right)^{\frac{j_0-1}{2}}\!\sum_{j=j_0}^{\infty}\!\left(\sqrt{\frac{2\lambda_n}{n}}\right)^{j-j_0}
\!\!\!\leq\lambda_n\left(\frac{2\lambda_n}{n}\right)^{\frac{j_0-1}{2}}\!\sum_{j=j_0}^{\infty}\!\left(\frac{1}{2}\right)^{j-j_0}
=2\lambda_n\left(\frac{2\lambda_n}{n}\right)^{\frac{j_0-1}{2}}
\end{equation}
for all $n\geq n_0$. Let us suppose that $n$ satisfies this condition from now on.

Now we bound $|R_n|$. First we apply the triangle inequality, then the inequality $|\e^{-x}-1+x|\leq\frac{x^2}{2}$ valid for all positive real $x$ with $x=\sum_{j=2}^{\infty}\frac{1}{j}\lambda_{n,j}$ , and finally use inequality (\ref{lambda_n,j_becslese}) with $j_0=2$ and 3. Thus we obtain
\begin{align*}
|R_n|
&\leq \e^{-\lambda_n}\left(\left|\sum_{j=3}^{\infty}\frac{1}{j}\lambda_{n,j}\right|
+\left|\exp\left\{-\sum_{j=2}^{\infty}\frac{1}{j}\lambda_{n,j}\right\}
-1+\sum_{j=2}^{\infty}\frac{1}{j}\lambda_{n,j}\right|\right)\\
&\leq \e^{-\lambda_n}\left(\sum_{j=3}^{\infty}\frac{1}{j}\lambda_{n,j}
+\frac{1}{2}\left(\sum_{j=2}^{\infty}\frac{1}{j}\lambda_{n,j}\right)^2\right)\\
&\leq \e^{-\lambda_n}\left(\frac{4\lambda_n^2}{n}
+\frac{1}{2}\left(2\lambda_n\sqrt{\frac{2\lambda_n}{n}}\right)^2\right)=\e^{-\lambda_n}4\lambda_n^2(\lambda_n+1)\frac{1}{n}.
\end{align*}
Recalling (\ref{hiba_e}) we see that we proved
\begin{equation}\label{ketto}
\e^{-\lambda_n}-\prod_{i=m+1}^{n}\frac{i}{n}
=\e^{-\lambda_n}\frac{1}{2}\sum_{i=m+1}^{n}\left(1-\frac{i}{n}\right)^2+O\!\left(\frac{1}{n}\right).
\end{equation}

\bigskip

Finally, recalling (\ref{p(w=k)}) we have
$$\p(\widetilde{W}_{n,m}=0)=\left(\prod_{i=m+1}^{n}\frac{i}{n}\right)=\e^{-\lambda_n}-\left(\e^{-\lambda_n}-\prod_{i=m+1}^{n}\frac{i}{n}\right)$$
for $k=0$,
$$\p(\widetilde{W}_{n,m}=1)=\left(\prod_{i=m+1}^{n}\frac{i}{n}\right)\lambda_n=\e^{-\lambda_n}\lambda_n
-\left(\e^{-\lambda_n}-\prod_{i=m+1}^{n}\frac{i}{n}\right)\lambda_n$$
for $k=1$, and
$$\p(\widetilde{W}_{n,m}=k)=\left(\prod_{i=m+1}^{n}\frac{i}{n}\right)S_k
=\e^{-\lambda_n}S_k-\left(\e^{-\lambda_n}-\prod_{i=m+1}^{n}\frac{i}{n}\right)S_k$$
for $k\geq2$. We obtain the first assertion of Theorem \ref{t_poisson_sorfejtes} if we write (\ref{egy}) and (\ref{ketto}) into these expressions. The second assertion follows from the first and (\ref{lambda_n}). $\blacksquare$


\bigskip

\section{Poisson approximation -- matching the means}

As mentioned at the beginning of this chapter, we shall now approximate the coupon collector's shifted waiting time $\widetilde{W}_{n,m}$ with another Poisson law, namely with the one that has the same mean as $\widetilde{W}_{n,m}$. One can easily calculate that in the range of parameters $n$ and $m$ for which the Poisson limit theorem of Section 1.2. holds true, the error order of this new approximation, given in the theorem below, is $1/n$, which is clearly better than the error order $1/\sqrt{n}$ given by Corollary \ref{kov_p1} or Theorem \ref{t_poisson_sorfejtes} for the same case. As we shall see, the proof of Theorem \ref{t_p_mean} is based on Stein's method, and heavily uses the fact that the means of the compared probability measures coincide. We note that the argument presented in the proof of Theorem \ref{t_upper_bound} would not work here.

\begin{theorem}\label{t_p_mean}
For the coupon collector's shifted waiting time $\widetilde{W}_{n,m}$ with $\lambda'_n=\E\widetilde{W}_{n,m}=\sum_{i=m+1}^{n}\left(\frac{n}{i}-1\right)$, we have
\begin{equation}
d_{\mathrm{TV}}({\cal D}(\widetilde{W}_{n,m}),\mathrm{Po}(\lambda_n'))\leq
8\left(1\wedge\sqrt{\frac{2}{\e\lambda_n'}}\right)\sum_{i=m+1}^{n}\left(\frac{n-i}{i}\right)^3.
\end{equation}
\end{theorem}

\bigskip

\noindent{\bf Proof.}
Recalling Section 2.3,  we apply the Stein-Chen method for Poisson approximation. By (\ref{stein_d_TV_formula}), we get the following formula:
\begin{equation}\label{d_TV_2}
d_{\mathrm{TV}}({\cal D}(\widetilde{W}_{n,m}),\mathrm{Po}(\lambda_n'))
=\sup_{A\subset\Z_+}|\E\{\lambda_n' f_A(\widetilde{W}_{n,m}+1)-\widetilde{W}_{n,m}f_A(\widetilde{W}_{n,m})\}|,
\end{equation}
where $f_A$ is the solution to the Stein equation (\ref{stein_egyenlet}), and by (\ref{stein_mokorlat}), we know that
\begin{equation}\label{fnorma_2}
\sup_{k\in \Z_+}|f_A(k)|\leq 1\wedge \sqrt{\frac{2}{\e\lambda_n'}}.
\end{equation}

Recalling the distributional equalities in (\ref{Wgeo}) and (\ref{hullamos}), we introduce
$$W_{n,m}^i:=W_{n,m}-X_i,\quad i\in\{m+1,\ldots, n\}$$
and
$$\widetilde{W}_{n,m}^i:=\widetilde{W}_{n,m}-\widetilde{X}_i,\quad i\in\{m+1,\ldots, n\}.$$
Taking an arbitrary $A\subset\Z_+$, the two terms on the right hand side of (\ref{d_TV_2}) can be written in the form
\begin{align*}
\E\left\{\lambda_n' f_A(\widetilde{W}_{n,m}+1)\right\}
&=\E(\widetilde{W}_{n,m})\E\left\{f_A(\widetilde{W}_{n,m}+1)\right\}\\
&=\sum_{i=m+1}^{n}\E(\widetilde{X}_i)\E\left\{f_A(\widetilde{W}_{n,m}+1)\right\}\\
&=\sum_{i=m+1}^{n}\sum_{k=1}^{\infty}k\p(\widetilde{X}_i=k)\E\left\{f_A(\widetilde{W}_{n,m}+1)\right\}
\end{align*}
by $(\ref{hullamos})$, and
\begin{align*}
\E\left\{\widetilde{W}_{n,m}f_A(\widetilde{W}_{n,m})\right\}
&=\sum_{i=m+1}^{n}\E\left\{\widetilde{X}_{i}f_A(\widetilde{W}_{n,m})\right\}\\
&=\sum_{i=m+1}^{n}\sum_{k=1}^{\infty}\E\left\{kf_A(\widetilde{W}_{n,m}^i+k) | \widetilde{X}_i=k\right\}\p( \widetilde{X}_i=k)\\
&=\sum_{i=m+1}^{n}\sum_{k=1}^{\infty}k\p(\widetilde{X}_i=k)\E\left\{f_A(\widetilde{W}_{n,m}^i+k)\right\},
\end{align*}
where we used $(\ref{hullamos})$ again and the independence of $\widetilde{X}_i$ and $\widetilde{W}_{n,m}^i$.
Putting these together, we get
\begin{align*}
&\E\{\lambda_n' f_A(\widetilde{W}_{n,m}+1)-\widetilde{W}_{n,m}f_A(\widetilde{W}_{n,m})\}=\\
&=\sum_{i=m+1}^{n}\sum_{k=1}^{\infty}k\p(\widetilde{X}_i=k)
\E\left\{f_A(\widetilde{W}_{n,m}+1)-f_A(\widetilde{W}_{n,m}^{i}+k)\right\}\\
&=\sum_{i=m+1}^{n}\sum_{k=1}^{\infty}k\p(\widetilde{X}_i=k)
\E\left\{f_A(\widetilde{W}_{n,m}^{i}+X_i)-f_A(\widetilde{W}_{n,m}^{i}+k)\right\}\\
&=\sum_{i=m+1}^{n}\sum_{k=1}^{\infty}k\p(\widetilde{X}_i=k)
\left(\sum_{l=1}^{\infty}
\E\left\{f_A(\widetilde{W}_{n,m}^{i}+l)-f_A(\widetilde{W}_{n,m}^{i}+k)|X_i=l\right\}\p(X_i=l)\right)\\
&=\sum_{i=m+1}^{n}\sum_{k=1}^{\infty}\sum_{l=1 \atop l\neq k}^{\infty}k
\p(X_i=k+1)\p(X_i=l)\E\left\{f_A(\widetilde{W}_{n,m}^{i}+l)-f_A(\widetilde{W}_{n,m}^{i}+k)\right\},
\end{align*}
where at the last step we used the independence again. Thus by $(\ref{fnorma_2})$
\begin{align*}
|\E\{\lambda_n'
&f_A(\widetilde{W}_{n,m}+1)-\widetilde{W}_{n,m}f_A(\widetilde{W}_{n,m})\}|
\leq\\
&\leq\sum_{i=m+1}^{n}\sum_{k=1}^{\infty}\sum_{l=1 \atop l\neq k}^{\infty}k
\p(X_i=k+1)\p(X_i=l)\E\left\{|f_A(\widetilde{W}_{n,m}^{i}+l)-f_A(\widetilde{W}_{n,m}^{i}+k)|\right\}\\
&\leq2\left(1\wedge\sqrt{\frac{2}{\e\lambda_n'}}\right)\sum_{i=m+1}^{n}\sum_{k=1}^{\infty}\sum_{l=1
\atop l\neq k}^{\infty}k\p(X_i=k+1)\p(X_i=l).
\end{align*}
Since this inequality holds for each $A\subset\Z_+$, it yields the same upper bound for the supremum of the expectations on left hand side taken over the set of functions $f_A$, $A\subset\Z_+$, so by $(\ref{d_TV_2})$, we have
\begin{equation}\label{d_TV2_2}
d_{\mathrm{TV}}({\cal D}(\widetilde{W}_{n,m}),\mathrm{Po}(\lambda_n'))\leq2\left(1\wedge\sqrt{\frac{2}{\e\lambda_n'}}\right)
\sum_{i=m+1}^{n}\sum_{k=1}^{\infty}\sum_{l=1\atop l\neq k}^{\infty}k\p(X_i=k+1)\p(X_i=l).
\end{equation}

Keeping in mind the infinite series' sums
$\sum_{j=1}^{\infty}x^{j-1}(1-x)=1$, $\sum_{j=1}^{\infty}jx^j=\frac{x}{(1-x)^2}$ and $\sum_{j=1}^{\infty}jx^{2j-1}=\frac{x}{(1-x)^2(1+x)^2}$ for any $0<x<1$, we start the calculation of the expression above:
\begin{align*}
\sum_{k=1}^{\infty}\sum_{l=1 \atop l\neq k}^{\infty}k&\p(X_i=k+1)\p(X_i=l)
=\sum_{k=1}^{\infty}k\left(1-\frac{i}{n}\right)^k\frac{i}{n}\left[1-\left(1-\frac{i}{n}\right)^{k-1}\frac{i}{n}\right]\\
&=\frac{i}{n}\sum_{k=1}^{\infty}k\left(1-\frac{i}{n}\right)^k-\left(\frac{i}{n}\right)^2
\sum_{k=1}^{\infty}k\left(1-\frac{i}{n}\right)^{2k-1}\\
&=\frac{i}{n}\frac{1-\frac{i}{n}}{\left[1-\left(1-\frac{i}{n}\right)\right]^2}-\left(\frac{i}{n}\right)^2
\frac{1-\frac{i}{n}}{\left[1-\left(1-\frac{i}{n}\right)\right]^2\left[1+\left(1-\frac{i}{n}\right)\right]^2}\\
&=\frac{1-\frac{i}{n}}{\frac{i}{n}}-\frac{1-\frac{i}{n}}{\left(2-\frac{i}{n}\right)^2}\\
&=\frac{n-i}{i}-\frac{n(n-i)}{(2n-i)^2}\\
&=\frac{(n-i)(4n^2-4ni+i^2-ni)}{i(2n-i)^2}\\
&=\frac{(n-i)^2(4n-i)}{i(2n-i)^2}
=\left(\frac{n-i}{i}\right)^3\frac{i^2(4n-i)}{(2n-i)^2(n-i)}.
\end{align*}
Now for an arbitrary $i\in[m+1,\ldots,n]$ the sequence $\frac{i^2(4n-i)}{(2n-i)^2(n-i)}$ can be bounded from above by
4. If we use this bound for the triple sum in (\ref{d_TV2_2}), we obtain the theorem. $\blacksquare$



\chapter{Compound Poisson approximation}

According to our goals set out in the introduction, in this chapter we approximate the distribution of the appropriately centered coupon collector's waiting time with a compound Poisson measure $\pi_{\mu,a}$ defined at the end of Section 1.3. Based on the distributional equality in (\ref{hullamos}), we shall apply general results of translated compound Poisson approximation of sums of independent integer valued random variables, which has been studied in a series of papers. Using Stein's method, \cite{BX} and \cite{BC} give bounds for the errors of such approximations in total variation distance. Their upper bounds are expressed with the help of the first three moments of the summands $X_1,X_2,\ldots,X_n$ and the critical ingredient $d_{\mathrm{TV}}\left({\cal D}\left(W_n\right),{\cal D}\left(W_n+1\right)\right)$, where $W_n=\sum_{j=1}^{n}X_j$.

The expression $d_{\mathrm{TV}}\left({\cal D}\left(W_n\right),{\cal D}\left(W_n+1\right)\right)$
is usually bounded by the Mineka coupling introduced in Section 2.4, which typically yields a
bound of order $1/\sqrt{n}$. If the $X_j$'s are roughly similar in magnitude, this is
comparable with the order $O(1/\sqrt{\var W_n})$ expected for the error in the central
limit theorem.  However, if the distributions of the~$X_j$ become progressively more
spread out as~$j$ increases, then $\var W_n$ may grow faster than~$n$, and then
$1/\sqrt n$ is bigger than the ideal order $O(1/\sqrt{\var W_n})$. In fact, this is the situation in the case when we chose $W_n$ to be the coupon collector's waiting time.

In the first part of this chapter we shall introduce a new coupling which allows us to improve the bounds obtained by the
Mineka coupling in such cases. Then, in the second part of the chapter, with the help of this new coupling, we shall prove
$d_{\mathrm{TV}}\left({\cal D}\left(W_{n,m}\right),{\cal D}\left(W_{n,m}+1\right)\right)
\\=O\left(1/\sqrt{\var W_{n,m}}\right)$, and therefore that a translated compound
Poisson approximation to the collector's waiting time $W_{n,m}$, with ideal
error rate, can be obtained in all ranges of $n$ and~$m$ in which a central or Poisson
limit theorem can be proved.


\section{An extension of Mineka's coupling inequality}

We saw in (\ref{Mineka_uniform}) that if $V_r$ is a sum of iid discrete uniform random variables on the finite interval $\{1,2,\ldots,2l-1,2l\}$ for some $l>1$ integer, then
$$
d_{\mathrm{TV}}({\cal D}(V_r),{\cal D}(V_r+1))\leq\frac{1}{\sqrt{2r}},
$$
where $1/\sqrt{r}>1/(l\sqrt{r})=1/\sqrt{\var V_r}$. However, in the following lemma we show that $d_{\mathrm{TV}}({\cal D}(V_r),{\cal D}(V_r+1))=1/\sqrt{\var V_r}$ can be established, if we use a new coupling instead of Mineka's coupling.

\bigskip

\begin{lemma}\label{l_uniform_coupling}
Let $U_1, U_2,\ldots, U_r$, $r\geq 2$, be independent identically distributed random variables with discrete uniform distribution on $\{1,2,\ldots,2l-1,2l\}$ for some integer $l\geq1$. If $V_r=\sum_{j=1}^{r}U_j$, then
\begin{equation*}
d_{\mathrm{TV}}({\cal D}(V_r),{\cal D}(V_r+1))\leq\frac{1}{l\sqrt{r}}.
\end{equation*}
\end{lemma}

\noindent {\bf Proof.} We construct a coupling of $(V_r,V_r+1)$. Let $U_1$ be an arbitrary random variable of uniform distribution on $\{1,2,\ldots,2l\}$. If $U_1\in\{1,2,\ldots,2l-1\}$, then define
$$U_1'= U_1+1\quad \textrm{and} \quad U'_j=U_j,\; 2\leq j\leq r,$$
where $U_1,\ldots,U_r$ are independent; while if $U_1=2l$, then put
$$U_1'=1\quad\textrm{ and }\quad U_j=\widetilde U_j+lI_j,\quad U_j'=\widetilde U_j+l(1-I_j),\quad 2\leq j\leq r,$$
where $\widetilde U_j$ has uniform distribution on $\{1,\ldots,l\}$, $I_j$ takes on the values 0 and 1, each with probability $1/2$, and $\widetilde U_j$, $I_j$, $2\leq j\leq r$, are independent, also of $U_1$.

Introducing $V_s:=\sum_{j=1}^{s}U_j$ and $V'_s:=\sum_{j=1}^{s}U'_j$, $s\in\{1,\ldots,r\}$, we see that
$$S_s:=(V_s+1)-V_s'=\left\{\begin{array}{ll}
                                                               0, & \hbox{if $U_1\in\{1,2,\ldots,2l-1\}$} \\
                                                               2l+\sum_{j=2}^{s}(U_j-U_j'), & \hbox{if $U_1=2l$,}
                                                             \end{array}
                                                           \right.
$$
where \vspace{-0.4 cm}
    $$U_j-U_j'=\left\{
                 \begin{array}{ll}
                   l, & \hbox{with probability $1/2$,} \\
                   -l, & \hbox{with probability $1/2$.}
                 \end{array}
               \right.
    $$
Thus if $U_1=2l$, $\left(S_s\right)_{s=1}^{r}$ can be regarded as a symmetric random
walk that starts from $2l$ in time step one, and then at each subsequent time step increases or decreases by $l$. Define $T$ to be the first time the random walk hits 0, that is
$$T:=\inf\left\{s\geq2: S_s=0\right\}=\inf\left\{s\geq2: \sum_{j=2}^{s}(U_j-U_j')=-2l\right\}.$$
By the reflection principle and symmetry,
\begin{align*}
\p(T>r|U_1=2l)
&=1-\p(T\leq r|U_1=2l)\\
&=1-\p\left(S_r=0|U_1=2l\right)-2\p\left(S_r<0|U_1=2l\right)\\
&=1-\p\left(S_r=0|U_1=2l\right)-\p\left(S_r<0|U_1=2l\right)-\p\left(S_r>4l|U_1=2l\right)\\
&=\sum_{k=1}^{4}\p\left(S_r=kl|U_1=2l\right)\leq 2\max_{k\in\Z}\p\left(S_r=kl|U_1=2l\right),
\end{align*}
and by Lemma 4.7 of Barbour and Xia \cite{BX}, we have
$$\max_{k\in\Z}\p\left(S_r=kl|U_1=2l\right)\leq\frac{1}{\sqrt{2}}\frac{1}{\sqrt{r-1}},$$
thus
\begin{equation}\label{T}
\p(T>r|U_1=2l)\leq\frac{2}{\sqrt{r}}.
\end{equation}

Now for $j,s\in\{1,\ldots,r\}$ put
$$U_j'':=\left\{
                           \begin{array}{ll}
                             U_j', & \hbox{if $1\leq j\leq T$,} \\
                             U_j, & \hbox{$j>T$,}
                           \end{array}
                         \right.\quad\textrm{and}\quad V_s'':=\sum_{j=1}^{s}U_j''.
$$
Of course $(V_s)_{s=1}^{r}$, $(V'_s)_{s=1}^{r}$ and $(V''_s)_{s=1}^{r}$ all have the same distribution, thus $(V''_r,V_r+1)$ is a coupling of $(V_r,V_r+1)$, therefore
$$d_{\mathrm{TV}}(V_r,V_r+1)\leq \p(V_r+1\neq V_r'')=\p(T>r)$$
by the coupling inequality. Since
$$\p(T>r)=\p(U_1=2l)\p(T>r|U_1=2l)\leq\frac{1}{l\sqrt{r}}$$
by (\ref{T}), the proof is complete. $\blacksquare$

\medskip

Now we show how the result of Lemma \ref{l_uniform_coupling} concerning sums of iid uniform random variables can be used to obtain similar results for sums of arbitrary integer valued random variables. The idea is to embed the uniform random variables in the ones we want to prove the result for.

\medskip

\begin{proposition}\label{p_coupling}
If $X_1,X_2,\ldots,X_n$, $n\geq2$, are independent integer valued random variables and $W=\sum_{j=1}^{n}X_n$, then
$$d_{\mathrm{TV}}({\cal D}(W),{\cal D}(W+1))\leq\frac{4}{l\sqrt{nlp}}+\frac{8d_n}{nlp},$$
where $l\in\{2,4,6,\ldots\}$ and $p\leq\min\{\p(X_j=k): k=1,\ldots,l, j=1,\ldots,n\}$ are arbitrary and $d_n=d_{\mathrm{TV}}({\cal D}(X_n),{\cal D}(X_n+1))$.
\end{proposition}

\medskip

\noindent{\bf Proof.} We write each of the variables $X_1,\ldots,X_{n}$ in the form
\begin{equation}\label{X_j_decomposition}
X_{j}=I_{j}U_{j}+(1-I_{j})R_{j},\quad j=1,\ldots,n,
\end{equation}
where $I_{j}$, $U_{j}$ and $R_{j}$, $j=1,\ldots,n$, are all independent random variables defined on a common probability space, and for each $j=1,\ldots,n$: $U_j$ has discrete uniform distribution on $\left\{1,2,\ldots,l\right\}$ for some even integer $l$; $I_j$ is a Bernoulli random variable with parameter $lp$, where $p\leq\min\{\p(X_j=k): k=1,\ldots,l, j=1,\ldots,n\}$ is fixed; and
$$
\p(R_j=k)=\left\{
             \begin{array}{ll}
               \frac{\p(X_j=k)-p}{1-lp}, & \hbox{$1\leq k\leq l$,} \\
               \frac{\p(X_j=k)}{1-lp}, & \hbox{otherwise,}
             \end{array}
           \right.\quad k\in\Z.
$$

Since ${\cal D}(X_j|I_{j}=1)={\cal D}(U_j)$ and ${\cal D}(X_j|I_{j}=0)={\cal D}(R_j)$, for any $\delta_{1},\ldots,\delta_{n-1}\in\{0,1\}$ and $\rho_{1},\ldots,\rho_{n-1}\in\Z$ we have
\begin{equation*}
{\cal D}\left(\sum_{j=1}^{n-1}X_j\big|I_j=\delta_j, R_j=\rho_j, j=1,\ldots,n-1\right)={\cal D}(V_r+\rho),
\end{equation*}
where $r=\sum_{j=1}^{n-1}\delta_j\,$, $V_r = \sum_{j=1}^r U_j'$, where the $U_j'$ are independent copies of $U_1$, and are independent of everything else, and $\rho=\sum_{j=1}^{n-1}(1-\delta_j)\rho_j$.

Now we apply the inequality
\begin{equation}\label{Z_123}
d_{\mathrm{TV}}({\cal D}(Z_1),{\cal D}(Z_2))\leq \E\{d_{\mathrm{TV}}({\cal D}(Z_1|Z_3),{\cal D}(Z_2|Z_3))\}
\end{equation}
true for any random elements $Z_1,Z_2$ and $Z_3$ defined on the same probability space. We obtain
\begin{align*}
&d_{\mathrm{TV}}({\cal D}(W),{\cal D}(W+1))\leq\\
&\leq \E\!\left\{\!d_{\mathrm{TV}}\!\!\left({\cal D}\!\left(\sum_{j=1}^{n} X_j|I_j,R_j,j=1,\ldots,n-1\right)\!,
{\cal D}\!\left(\sum_{j=1}^{n} X_j+1|I_j,R_j,j=1,\ldots,n-1\right)\right)\!\right\}\\
&= \E\{d_{\mathrm{TV}}({\cal D}(V_T + X_n + R | T,R), {\cal D}(V_T + X_n + R + 1 | T,R))\},
\end{align*}
where $T = \sum_{j=1}^{n-1} I_j$ and $R = \sum_{j=1}^{n-1} (1 - I_j)R_j$
are independent of $(U_j',\,j\ge1)$ and of $X_n$.  Hence
\begin{equation}\label{d_+bound}
   d_{\mathrm{TV}}({\cal D}(W),{\cal D}(W+1)) \ \le\
  	\E\{d_{\mathrm{TV}}({\cal D}(V_T + X_n | T), {\cal D}(V_T + X_n + 1 | T))\},
\end{equation}
since total variation distance is invariant under translation.

Now, since $T$, $X_n$ and $(U_j',\,j\ge1)$ are independent, we have
\begin{multline*}
  	d_{\mathrm{TV}}({\cal D}(V_T + X_n | T=t), {\cal D}(V_T + X_n + 1 | T=t))\\
		\ \le\ \min\{d_{\mathrm{TV}}({\cal D}(V_t), {\cal D}(V_t+1)),d_{\mathrm{TV}}({\cal D}(X_n), {\cal D}(X_n+1))\},
\end{multline*}
and Lemma \ref{l_uniform_coupling} provides the bound
\begin{equation*}
d_{\mathrm{TV}}({\cal D}(V_t),{\cal D}(V_t+1))\leq f(t):=\left\{
                                         \begin{array}{ll}
                                           \frac{2}{l\sqrt{t}}, & \hbox{if $t>0$,} \\
                                           1, & \hbox{if $t=0$.}
                                         \end{array}
                                       \right.
\end{equation*}
Writing $d_n=d_{\mathrm{TV}}({\cal D}(X_n),{\cal D}(X_n+1))$ we thus obtain from (\ref{d_+bound}) that
\begin{align*}
d_{\mathrm{TV}}({\cal D}(W),{\cal D}(W+1))
&\leq \E\{d_{\mathrm{TV}}({\cal D}(V_T + X_n | T), {\cal D}(V_T + X_n + 1 | T))\}\\
&\leq \E\left\{\min\left[f(T); d_{n}\right]\right\}\\
&\leq \E\left\{\frac{2}{l\sqrt{T}}\Big|T\geq\frac{\E T}{2}\right\}\p\left(T\geq\frac{\E T}{2}\right)+d_{n}\p\left(T<\frac{\E T}{2}\right)\\
&\leq \frac{2\sqrt{2}}{l\sqrt{\E T}}+d_n\p\left(T<\frac{\E T}{2}\right).
\end{align*}
Since $T$ has distribution Bin$(n-1,lp)$, $\E T=(n-1)lp\geq\frac{1}{2}nlp$, and by Chebishev's inequality
$$\p\left(T<\frac{\E T}{2}\right)
\leq \p\left(|T-\E T|>\frac{\E T}{2}\right)
\leq\frac{4\var T}{(\E T)^2}
\leq\frac{4}{(n-1)lp}
\leq\frac{8}{nlp},$$
thus
\begin{equation}\label{INEQ}
d_{\mathrm{TV}}({\cal D}(W),{\cal D}(W+1))\leq\frac{4}{l\sqrt{nlp}}+\frac{8d_n}{nlp}.
\end{equation}$\blacksquare$

\medskip

\begin{remark}
Since total variation distance is invariant under translation, there is no loss of generality in supposing that the $l$-intervals begin at~$1$.
\end{remark}

\medskip

\begin{remark}
The choice of $(p,l)$ depends very much on the problem.
\end{remark}

\medskip

The constants in the upper bound of Proposition \ref{p_coupling} can be improved by refining the method proposed in the proof. One could embed not one, but many uniform random variables in the $X_j$-s by splitting the whole line into the $l$-blocks $(\{(m-1)l,\ldots,ml\})_{m\in\Z}$ and defining a uniform variable corresponding to each block. Thus one could use potential overlaps from the whole distribution and not just the interval $\{1,\ldots,l\}$, when bounding $d_{\mathrm{TV}}({\cal D}(W),{\cal D}(W+1))$.

More precisely, each $X_j$, $j=1,\ldots,n$, can be given in the form
$$X_j=I_{j0}R_j+\sum_{i=1}^{\infty}I_{ji}(U_{ji}+(i-1)l),$$
where all random variables in the decompositions are defined on a common probability space, and for each $j=1,\ldots,n$ the following hold true: $U_{ji}$ has discrete uniform distribution on $\{1,\ldots,l\}$, $i=1,2,\ldots$, for some fixed even integer $l$; $I_{j0}\sim\textrm{Bernoulli}\left(1-\sum_{i=1}^{\infty}lp_i\right)$, $I_{ji}\sim\textrm{Bernoulli}(lp_i)$, where $p_i\leq\min\{\p(X_j=k):k=(i-1)l,\ldots,il, j=1,\ldots,n\}$ is fixed, $i=1,2,\ldots$,
and these Bernoulli variables depend on each other in a way that for each outcome exactly one of them is 1 and the rest are 0; all the other variables in the decompositions are independent of each other and of the $I_{ij}$-s; and $R_j$ is defined to make the distribution of the decomposition equal the distribution of $X_j$.

Then, to bound $d_{\mathrm{TV}}({\cal D}(W),{\cal D}(W+1))$ we would use (\ref{Z_123}), conditioning on all the $I_{ji}$-s and $R_j$-s, which would give us (\ref{d_+bound}) with $T=\sum_{j=1}^{n-1}\sum_{i=1}^{\infty}I_{ji}$. In this case $\E T=(n-1)l\sum_{i=1}^{\infty}p_i$ and $\var T=(n-1)l\left(\sum_{i=1}^{\infty}p_i\right)\left(1-l\sum_{i=1}^{\infty}p_i\right)$, hence we would obtain (\ref{INEQ}) with $p$ replaced by $\sum_{i=1}^{\infty}p_i$.


\section{Compound Poisson approximation in the range of the central and Poisson limit theorems}

We return to the coupon collector's problem. Taking advantage of the decomposition in (\ref{hullamos}), we apply a theorem of Barbour and Xia \cite{BX}
on translated compound Poisson approximation in total variation distance to the
distributions of sums of independent integer valued random variables. One of the
elements in their approximation error is (almost)
$d_{\mathrm{TV}}({\cal D}(W_{n,m}), {\cal D}(W_{n,m}+1))$, to bound which we invoke our
proposition of the previous section. Recalling that in Section 1.3, for $\mu, a>0$, we defined the compound Poisson distribution $\pi_{\mu,a}$
to be the distribution of $Z_1+2Z_2$, where $Z_1\sim\mathrm{Po}(\mu)$ and
$Z_2\sim\mathrm{Po}(a/2)$ are independent, we have the following result:

\medskip

\begin{theorem}\label{t_CP}
For any fixed $n\geq2$ and $2\leq m\leq n-4$, if
\begin{align}\nonumber
&\mu=\sigma_n^2-2\langle \sigma_n^2-\mu_n\rangle,\\ \label{acmu}
&a=\langle \sigma_n^2-\mu_n\rangle\quad\textrm{and}\\
&c=\lfloor \sigma_n^2-\mu_n\rfloor,\nonumber
\end{align}
where $\langle x\rangle$ and $\lfloor x\rfloor$ denote the fractional and integer part of $x$ respectively,
then there exists a positive constant $C$ such that
\begin{equation}\label{discrete_approx1}
d_{\mathrm{TV}}\Big({\cal D}\left(W_{n,m}+c\right),\pi_{\mu,a}\Big)
\leq \frac{C}{\sigma_n}\left(\frac{\left\lfloor\sigma_n^2-\mu_n-(n-m)\right\rfloor}{\sigma_n^2}+\frac{(n-m)^2}{nm}\right).
\end{equation}

\end{theorem}

\medskip

\begin{remark}
We recall one of the Baum--Billingsley theorems from Section 1.2: if $m=m\in\{0,1,\ldots,n-1\}$ is an integer that depends on $n$ in such a way that
\begin{equation*}
m\to\infty\quad\textrm{ and }\quad\frac{n-m}{\sqrt{n}}\to\infty\quad\textrm{ as }\quad n\to\infty,
\end{equation*}
then $\overline{W}_{n,m}:=(W_{n,m}-\E W_{n,m})/\sqrt{\var W_{n,m}}$ has asymptotically standard normal distribution. This limit theorem was refined in Chapter 4 by showing that
\begin{equation*}
d_{\mathrm{K}}\left({\cal D}\left(\overline{W}_{n,m}\right),\mathrm{N}(0,1)\right)\leq C\frac{n}{m}\frac{1}{\sigma_n},
\end{equation*}
where $C=9.257$
. One can deduce that the same or better order of approximation is obtained in the discrete approximation given in our theorem, than with normal approximation, and now with the error measured with respect to the much stronger total variation distance.
\end{remark}

\medskip

\begin{remark}
Note that, with these parameters, $\pi_{\mu,a}$ has mean $\mu+a=\sigma_n^2-\langle \sigma_n^2-\mu_n\rangle=\mu_n+c$ and variance $\mu+2a=\sigma_n^2$.
\end{remark}

\medskip

\begin{remark}
We can express the bound of Theorem \ref{t_CP} more intuitively with the help of the asymptotic formulae given by Baum and Billingsley in \cite{BB} for the variance of the waiting time: if $n\to\infty$, then
\begin{equation*}\begin{array}{lll}
\frac{m}{n}\to 0, & \Rightarrow & \sigma_n^2\sim\frac{n^2}{m}\\
\frac{m}{n}\to c, c\in(0,1) & \Rightarrow & \sigma_n^2\sim\gamma n,\\
\frac{m}{n}\to 1, & \Rightarrow & \sigma_n^2\sim \frac{1}{2}\frac{(n-m)^2}{n},
\end{array}\end{equation*}
where $\gamma=(1-c+c\log c)/c$. In the latter case we distinguish two subcases: the case when $\liminf_{n,m\to\infty}\left\lfloor\sigma_n^2-\mu_n-(n-m)\right\rfloor>1$ and when $\limsup_{n,m\to\infty}\left\lfloor\sigma_n^2-\mu_n-(n-m)\right\rfloor<1$. We shall refer to the four categories above as "small", "medium", "large" and "very large" $m$. By these formulae, (\ref{discrete_approx1}) is equivalent to
\begin{equation}\label{discrete_approx2}
d_{\mathrm{TV}}\Big({\cal D}\Big(W_{n,m}+c\Big),\pi_{\mu,a}\Big)
=\left\{
      \begin{array}{ll}
        O\!\left(\frac{1}{\sqrt{m}}\right), & \hbox{in the "small" $m$ case;} \\
        O\!\left(\frac{1}{\sqrt{n}}\right), & \hbox{in the "medium" $m$ case;} \\
        O\!\left(\frac{1}{\sqrt{n}}\right), & \hbox{in the "large" $m$ case;}\\
        O\!\left(\frac{n-m}{n^{3/2}}\right), & \hbox{in the "very large" $m$ case.}
      \end{array}
    \right.
\end{equation}
\end{remark}

\bigskip

\noindent{\textbf{Proof of Theorem \ref{t_CP}.}} We apply Theorem 4.3 in \cite{BX}, which states that if $Z_j$, $j=1,\ldots,r$, are independent integer valued random variables with $\E|Z_j|^3<\infty$, $W=\sum_{j=1}^{r}Z_j$, and we define
$$e_c(W)=\left\{
           \begin{array}{ll}
             1, & \hbox{if $W+c\geq0$ almost surely;} \\
             0, & \hbox{otherwise;}
           \end{array}
         \right.$$
$$\psi_j:=\E|Z_j(Z_j-1)(Z_j-2)|+|\E Z_j|\E|Z_j(Z_j-1)|+2\E|Z_j||\var Z_j-\E Z_j|,$$
$$d_+:=\max_{1\leq i\leq r}\Big\{d_{\mathrm{TV}}({\cal D}(W_i),{\cal D}(W_i+1))\Big\},\textrm{ where }W_i:=W-Z_i,$$
then with $\mu=\var W-2\langle \var W-\E W\rangle$, $a=\langle \var W-\E W\rangle$ and $c=\lfloor \var W-\E W\rfloor$,
\begin{equation}\label{BX_theorem}
d_{\mathrm{TV}}\Big({\cal D}\left(W+c\right),\pi_{\mu,a}\Big)
\leq\frac{2e_c(W)+2\left(|\lfloor \var W-\E W\rfloor|+\sum_{j=1}^{r}\psi_j\right)d_+}{\var W}.
\end{equation}
We apply this theorem with $Z_j=X_j-1$, $j\in\{m+1,\ldots,n\}$, for the $X_j$ given in (\ref{hullamos}), in order to approximate the coupon collector's shifted waiting time $\widetilde{W}_{n,m}:=\sum_{j=m+1}^{n}[X_j-1]$, and then show that the upper bound in (\ref{BX_theorem}) is not greater than the right hand side of (\ref{discrete_approx1}). Then, since the two measures compared in (\ref{BX_theorem}) are the same for $W=\widetilde{W}_{n,m}$ and $W=W_{n,m}=\widetilde{W}_{n,m}+n-m$, the theorem for $W_{n,m}$ follows immediately.

To do so, for given $n\geq2$, $2\leq m\leq n-4$ and $j\in\{m+1,\ldots,n\}$, we bound $\psi_j$ and $d_+$ as defined above.

For $X$, a random variable that has geometric distribution with parameter $p$, we have
$$
\E X=\frac{1}{p},\quad
\E X^2=\frac{2-p}{p^2},\quad
\E X^3=\frac{p^2-6p+6}{p^3},\;\textrm{ and }\;
\var X=\frac{1-p}{p^2}.
$$
Thus for $Z=X-1$, one can easily calculate
\begin{align*}
&\psi:
=\E|Z(Z-1)(Z-2)|+|\E Z|\E|Z(Z-1)|+2\E|Z||\var Z-\E Z|\\
&=\E\{X^3-6X^2+11X-6\}+\E\{X-1\}\E\{X^2-3X+2\}+2\E\{X-1\}|\var X-\E X+1|\\
&=\frac{10(1-p)^3}{p^3},
\end{align*}
so $\psi_j=10\left(\frac{n}{j}-1\right)^3$. If we add the $\psi_j$ together, we obtain
\begin{equation}\label{piros_psik}
\sum_{j=m+1}^{n}\psi_j
=\sum_{j=m+1}^{n}\frac{10\left(1-\frac{j}{n}\right)^3}{\left(\frac{j}{n}\right)^3}
\leq 10\frac{(n-m)^2}{n m}\sum_{j=m+1}^{n}\frac{n(n-j)}{j^2}
=10\frac{(n-m)^2}{n m}\sigma_n^2.
\end{equation}

Next, we notice that $e_c(\widetilde{W}_{n,m})=0$ almost surely. Indeed, $0\leq \sigma_n^2-\E\widetilde{W}_{n,m}\leq \sigma_n^2$, because for each $X_j$ geometric random variable of parameter $j/n$ we have $\var X_j-\E(X_j-1)=\left(\frac{1-j/n}{j/n}\right)^2\leq\frac{1-j/n}{(j/n)^2}=\var X_j$.

Now combining the bound in (\ref{piros_psik}) with inequality (\ref{BX_theorem}) applied to the $\widetilde{W}_{n,m}$ waiting time, we obtain
\begin{equation}\label{proof}
d_{\mathrm{TV}}\Big({\cal D}\Big(\widetilde{W}_{n,m}+\tilde{c}\Big),\pi_{\tilde{\mu},\tilde{a}}\Big)
\leq20\left(\frac{\left\lfloor\sigma_n^2-\mu_n-(n-m)\right\rfloor}{\sigma_n^2}+\frac{(n-m)^2}{nm}\right)d_+,
\end{equation}
where $\tilde{c}$, $\tilde{\mu}$ and $\tilde{a}$ are defined by the formulae in (\ref{acmu}) with $W_{n,m}$ replaced with $\widetilde{W}_{n,m}$.

Before turning to the approximation of $d_+$, we bound $\sigma_n^2$. We see that
\begin{equation*}
\sigma_n^2=n\sum_{j=m+1}^{n}\frac{n-j}{j^2}\leq\frac{n}{(m+1)^2}\sum_{j=m+1}^{n}(n-j)=\frac{n(n-m)(n-m-1)}{2(m+1)^2},
\end{equation*}
also,
\begin{equation*}
\sigma_n^2=n\sum_{j=m+1}^{n}\frac{n-j}{j^2}\leq n(n-m-1)\int_{m}^{n}\frac{1}{x^2}d x\leq\frac{n(n-m-1)}{m},
\end{equation*}
thus
\begin{equation}\label{var}
\sigma_n^2\leq n(n-m-1)\min\left\{\frac{n-m}{2(m+1)^2},\frac{1}{m}\right\}.
\end{equation}

Now for $d_+$, by an inequality of Mattner and Roos \cite{MR} we have
\begin{equation*}
d_{\mathrm{TV}}({\cal D}(W_i),{\cal D}(W_i+1))
\leq\sqrt{\frac{2}{\pi}}\left(\sum_{j=m+1,j\neq i}^{n}\left[1-d_{\mathrm{TV}}({\cal D}(X_j),{\cal D}(X_j+1))\right]\right)^{-\frac{1}{2}},
\end{equation*}
and since $d_{\mathrm{TV}}({\cal D}(X_j),{\cal D}(X_j+1))$ is equal to
\begin{equation}\label{D_tv_geo}
\frac{1}{2}\sum_{k=1}^{\infty}|\p(X_j=k)-\p(X_j=k-1)|
=\frac{1}{2}\left(\frac{j}{n}+\left(\frac{j}{n}\right)^2\sum_{k=2}^{\infty}\left(1-\frac{j}{n}\right)^{k-2}\right)
=\frac{j}{n},
\end{equation}
we obtain
\begin{equation*}
d_+
\leq\sqrt{\frac{2}{\pi}}\left(\sum_{j=m+1}^{n}\left(1-\frac{j}{n}\right)-\max_{m+1\leq i\leq n}\left(1-\frac{i}{n}\right)\right)^{-\frac{1}{2}}
=\sqrt{\frac{2}{\pi}}\frac{\sqrt{n}}{\sqrt{(n-m-1)(n-m-2)}}.
\end{equation*}
It follows from this and (\ref{var}) that for any $K>0$
\begin{equation}\label{d_+}
d_+\leq\frac{2K}{\sqrt{\pi}}\frac{1}{\sigma_n},\quad\textrm{ if }\;\frac{n}{K}\leq m\leq n-4.
\end{equation}
Putting this bound into (\ref{proof}) gives a result which, when compared to (\ref{discrete_approx1}), has an extra factor $K \ge n/m$.  Thus it is of inferior order if $m \ll n$.  To prove the theorem for such values of~$m$, we need to use our proposition to bound~$d_+$.

Let us assume that $2\leq m\leq\frac{n}{2}$. If we apply the Proposition to the random variables $\{X_{j},j=m+1,\ldots,2m,j\neq i\}$, $i\in\{m+1,\ldots,2m\}$ fixed, with
$$l=\left\{
                                          \begin{array}{ll}
                                            \lfloor\frac{n}{m}\rfloor, & \hbox{if $\lfloor\frac{n}{m}\rfloor$ is even,} \\
                                            \lfloor\frac{n}{m}\rfloor-1, & \hbox{if $\lfloor\frac{n}{m}\rfloor$ is odd,}
                                          \end{array}
                                        \right.
\quad\textrm{and}\quad p=\left(1-\frac{2m}{n}\right)^{l}\frac{m}{n},
$$
we obtain
\begin{equation}\label{prop}
d_{\mathrm{TV}}\left\{{\cal D}\left(\sum_{j=m+1,j\neq i}^{2m}X_j\right)\!,{\cal D}\left(\sum_{j=m+1,j\neq i}^{2m}X_j+1\right)\right\}
\leq\frac{2}{l\sqrt{(m-1)lp}}+\frac{8d}{(m-1)lp},
\end{equation}
where
$$
d=\left\{
    \begin{array}{ll}
      d_{\mathrm{TV}}\left\{{\cal D}(X_{2m}),{\cal D}(X_{2m}+1)\right\}=\frac{2m}{n}, & \hbox{if $i\neq2m$,} \\
      d_{\mathrm{TV}}\left\{{\cal D}(X_{2m-1}),{\cal D}(X_{2m-1}+1)\right\}=\frac{2m-1}{n}, & \hbox{if $i=2m$}
    \end{array}
  \right.
$$
by (\ref{D_tv_geo}). For any $i\in\{m+1,\ldots,2m\}$ we have
$$d\leq\frac{2m}{n}\quad\textrm{and}\quad l\geq\frac{n}{2m},$$
since $\lfloor x\rfloor-1\geq\frac{x}{2}$, if $x\geq2$, and
$$lp\geq\frac{n}{2m}\left(1-\frac{2m}{n}\right)^{\frac{n}{m}}\frac{m}{n}\geq \frac{\e^{-2}}{2},$$
because $(1-x)^{\frac{2}{x}}$ decreases as $x$ increases in (0,1), and its limit at 0 is $\e^{-2}$. Now putting the bounds above together in (\ref{prop}) yields
\begin{align*}
d_+
&=\max_{i\in\{m+1,\ldots,2m\}}d_{\mathrm{TV}}\left\{{\cal D}\left(\sum_{j=m+1,j\neq i}^{2m}X_j\right),{\cal D}\left(\sum_{j=m+1,j\neq i}^{2m}X_j+1\right)\right\}\\
&\leq8\sqrt{2}\e\frac{m}{\sqrt{m-1}n}+32\e^2\frac{m}{(m-1)n}
\leq16\e\frac{\sqrt{m}}{n}+64\e^2\frac{1}{n}
\leq(16\e+64\e^2)\frac{\sqrt{m}}{n},
\end{align*}
where the last two inequalities hold for $m\geq2$. By (\ref{var}), $\frac{\sqrt{m}}{n}\leq\frac{1}{\sigma_n}$, thus we have
\begin{equation*}
d_+\leq(16\e+64\e^2)\frac{1}{\sigma_n},\quad\textrm{ if }2\leq m\leq\frac{n}{2}.
\end{equation*}
This and (\ref{d_+}) with $K=2$ substituted into (\ref{proof}) yield the theorem. $\blacksquare$


\chapter{Poisson--Charlier expansions}
\addtocounter{theorem}{-1}

In the final chapter of the thesis we approximate the coupon collector's shifted waiting time $\widetilde{W}_{n,m}=W_{n,m}-(n-m)$ with Poisson--Charlier signed measures in total variation distance. To do so, we shall apply a characteristic function technique proposed in \cite{BKN}. Throughout the chapter $C$ and $C_R$ denote positive constants, not necessarily the same ones at different occurrences, the first one is always a universal constant, while the latter one depends on $R$.

Let $\mu={\cal D}(\widetilde{W}_{n,m}+c)$, where
\begin{equation*}
c=\lfloor\var\widetilde{W}_{n,m}-\E\widetilde{W}_{n,m}\rfloor=\lfloor\sigma_n^2-[\mu_n-(n-m)]\rfloor
=\left\lfloor\sum_{k=m+1}^{n}\left(\frac{n-k}{k}\right)^2\right\rfloor.
\end{equation*}
Recalling the distributional equality (\ref{hullamos}), we see that since the characteristic function of the geometric distribution with success probability $p\in(0,1)$ is $\e^{\i t}/(1-\frac{1-p}{p}(\e^{\i t}-1))$, the Fourier--Stieltjes transform of $\mu$ is
\begin{align}\nonumber
\phi(t):&=\int_{-\infty}^{\infty}\e^{\i tx}\d\mu(x)
=\prod_{k=m+1}^{n}\left(\frac{\e^{\i t}}{1-\frac{n-k}{k}(\e^{\i t}-1)}\e^{-\i t}\right)\e^{\i tc}\\ \label{uj_phi}
&=\exp\left\{-\sum_{k=m+1}^{n}\log\left(1-\frac{n-k}{k}(\e^{\i t}-1)\right)\right\}
\exp\left\{\i t\left\lfloor\sum_{k=m+1}^{n}\left(\frac{n-k}{k}\right)^2\right\rfloor\right\}.
\end{align}
Introducing the new variable $w=w_t=\e^{\i t}-1$ and the sequences
\begin{equation}\label{uj_anj_def}
a_{n,j}:=\sum_{k=m+1}^{n}\left(\frac{n-k}{k}\right)^j,\quad j=1,2,\ldots,
\end{equation}
we write $\phi(t)$ in the form
\begin{align*}
\phi(t)
&=\exp\left\{-\sum_{k=m+1}^{n}\log\left(1-\frac{n-k}{k}(\e^{\i t}-1)\right)\right\}
\exp\left\{\i t\left\lfloor\sum_{k=m+1}^{n}\left(\frac{n-k}{k}\right)^2\right\rfloor\right\}\\
&=\exp\left\{-\sum_{k=m+1}^{n}\log\left(1-\frac{n-k}{k}w\right)+\lfloor a_{n,2}\rfloor\log(1+w)\right\}.
\end{align*}
Assuming $|t|\leq m/n$, for any $k=m+1,\ldots,n-1$ we have
$$
\frac{n-k}{k}|w|
=\frac{n-k}{k}|\e^{\i t}-1|
\leq\frac{n-m-1}{m+1}|\e^{\i t}-1|
\leq\frac{n-m-1}{m+1}|t|
\leq1,
$$
which together with $|w|=|\e^{\i t}-1|\leq |t|\leq1$ allows us to expand the logarithmic expressions in $\phi(t)$, therefore
\begin{equation*}
\phi(t)
=\exp\left\{\sum_{r=1}^{\infty}\Big(a_{n,r}+(-1)^{r+1}\lfloor a_{n,2}\rfloor\Big)\frac{w^r}{r}\right\},\quad |t|\leq\frac{m}{n}.
\end{equation*}
We note that $a_{n,1}=\mu_n-(n-m)$ and $a_{n,2}=\sigma_n^2-\left[\mu_n-(n-m)\right]$, hence the line above can be rewritten as
\begin{equation}\label{phifelbontas}
\phi(t)=\chi(t)\exp\{h(w)\},\quad |t|\leq\frac{m}{n},
\end{equation}
where
\begin{equation}\label{uj_chi}
\chi(t)=\exp\left\{\sigma_n^2(\e^{\i t}-1)\right\}
\end{equation}
is the characteristic function of the Poisson distribution with parameter $\sigma_n^2$, and
\begin{equation}\label{uj_h}
h(w)=-\Big(a_{n,2}-\lfloor a_{n,2}\rfloor\Big)w+\Big(a_{n,2}-\lfloor a_{n,2}\rfloor\Big)\frac{w^2}{2}
+\sum_{r=3}^{\infty}\Big(a_{n,r}+(-1)^{r+1}\lfloor a_{n,2}\rfloor\Big)\frac{w^r}{r}.
\end{equation}

Now we fix an integer $R\geq3$, and modify the function $\exp\{h(w)\}$ in (\ref{phifelbontas}) in two steps, each time replacing a certain expression of the previous function with the first $R$ terms of its series expansion around 0:
\begin{equation}\label{uj_approx}
\exp\{h(w)\}\approx\exp\{h_R(w)\}\approx H_{R}(w),
\end{equation}
where
\begin{equation}\label{uj_h_R}
h_R(w)=-\Big(a_{n,2}-\lfloor a_{n,2}\rfloor\Big)w+\Big(a_{n,2}-\lfloor a_{n,2}\rfloor\Big)\frac{w^2}{2}
+\sum_{r=3}^{R}\Big(a_{n,r}+(-1)^{r+1}\lfloor a_{n,2}\rfloor\Big)\frac{w^r}{r}
\end{equation}
and
\begin{equation}\label{uj_H_R}
H_{R}(w)=\sum_{l=0}^{R}\frac{h_R^l(w)}{l!}.
\end{equation}

We approximate the distribution $\mu={\cal D}(\widetilde{W}_{n,m}+c)$ with $\nu_R$, which we define to be the finite signed measure on the nonnegative integers whose characteristic function is
\begin{equation}\label{psifelbontas}
\psi(t)=\chi(t)H_{R}(\e^{\i t}-1),\quad t\in{\mathbb R}.
\end{equation}
Since $H_{R}(\e^{\i t}-1)$ is a polynomial of $\e^{\i t}-1$ of the form
$\sum_{r=0}^{R^2}\widetilde{a}_{n,m}^{(r)}(\e^{\i t}-1)^r$, where the $\widetilde{a}_{n,m}^{(r)}$ coefficients all depend on $n$ and $m$, it follows that $\nu_R=\nu_R(\sigma_n^2,\widetilde{a}_{n,m}^{(1)},\ldots,\widetilde{a}_{n,m}^{(R^2)})$ is a Poisson--Charlier signed measure, which, according to (\ref{Poisson--Charlier_measure}), is defined by
\begin{equation}\label{uj_nu_R}
\nu_R\{j\}=\mathrm{Po}(\sigma_n^2)\{j\}\left(1+\sum_{r=1}^{R^2}(-1)^{r+1}\widetilde{a}_{n,m}^{(r)}C_r(j,\sigma_n^2)\right),\quad j\in{\mathbb N},
\end{equation}
where $C_r(j,\sigma_n^2)$ denotes the $r$-th Charlier polynomial given in (\ref{Charlier_polynomial}).

\bigskip

\begin{theorem}\label{t_PCexp_1}
We assume $a_{n,2}>1$. For an arbitrary integer $R\geq3$ there exist threshold numbers $m_R$ and $n_R$ depending on $R$ such that if $m\geq m_R$ and $n\geq n_R$, then
\begin{equation*}
\sup_{k\in\Z}|\mu\{k\}-\nu_R\{k\}|\leq C_R\left(\frac{1}{\sqrt{m}}\right)^{R},\quad \textrm{if }m\leq\frac{n}{2}-1,
\end{equation*}
and
\begin{equation*}
\sup_{k\in\Z}|\mu\{k\}-\nu_R\{k\}|\leq C_R\frac{(\sqrt{n})^{R-2}}{(n-m)^{R-1}},\quad \textrm{if }m\geq\frac{n}{2}.
\end{equation*}
\end{theorem}


\bigskip

Before embarking on the proof of the theorem, we prove a sequence of propositions, which we shall need later.

\begin{proposition}\label{p_sigma_n}
Assume $1\leq m\leq n-1$ and $n-m\geq2$.
\begin{equation}\label{uj_s2mk}
\textrm{If } m\leq\frac{n}{2}-1, \quad\textrm{ then }\quad \frac{1}{20}\frac{n^2}{m}\leq\sigma_n^2\leq\frac{n^2}{m}.
\end{equation}
\begin{equation}\label{uj_s2mn}
\textrm{If } m\geq\frac{n}{2}, \quad\textrm{ then }\quad \frac{1}{24}\frac{(n-m)^2}{n}\leq\sigma_n^2\leq2\frac{(n-m)^2}{n}.
\end{equation}
\end{proposition}

\bigskip

\noindent{\bf Proof.}
First we prove the upper bounds. Since the terms in the sum $\sigma_n^2=\sum_{k=m+1}^{n}\frac{n(n-k)}{k^2}$ decrease as $k$ increases, we have the bound
\begin{equation*}
\sigma_n^2
\leq\int_{m}^{n}\frac{n(n-x)}{x^2}\d x
=n^2\int_{m}^{n}\frac{1}{x^2}\d x-n\int_{m}^{n}\frac{1}{x}\d x
=n\left(\frac{n}{m}-1-\log\frac{n}{m}\right).
\end{equation*}
Thus we see that $\frac{n^2}{m}$ is always an upper bound for $\sigma_n^2$. However, if $\frac{n}{m}\leq2$, one can apply the inequality $\log x\geq x-1-\frac{(x-1)^2}{2}$ with $x=\frac{n}{m}$ to obtain
$$
\sigma_n^2\leq n\frac{1}{2}\left(\frac{n-m}{m}\right)^2
\leq2\frac{(n-m)^2}{n},
$$
where at the last inequality we used the assumption $m\geq\frac{n}{2}$.

The proof of the lower bounds is similar. Again, we use the fact that the terms in the sum $\sigma_n^2=\sum_{k=m+1}^{n}\frac{n(n-k)}{k^2}$ decrease as $k$ increases to obtain
\begin{equation*}
\sigma_n^2
\geq\int_{m+1}^{n}\frac{n(n-x)}{x^2}\d x
=n^2\int_{m+1}^{n}\frac{1}{x^2}\d x-n\int_{m+1}^{n}\frac{1}{x}\d x
=n\left(\frac{n}{m+1}-1-\log\frac{n}{m+1}\right).
\end{equation*}
Now if $\frac{n}{m+1}\geq2$, then by the inequality $1+\log x\leq\frac{1+\log2}{2}x$, $x\geq2$, with $x=\frac{n}{m+1}$,
\begin{equation*}
\sigma_n^2
\geq n\left(\frac{n}{m+1}-1-\log\frac{n}{m+1}\right)
\geq\left(1-\frac{1+\log2}{2}\right)\frac{n^2}{m+1}
\geq\frac{1}{10}\frac{n^2}{m+1}
\geq\frac{1}{20}\frac{n^2}{m},
\end{equation*}
where we also used $m+1\leq2m$. If $\frac{n}{m+1}\leq2$, we can apply the inequality
$\log x\leq(x-1)-\frac{(x-1)^2}{2}+\frac{(x-1)^3}{3}$, $0\leq x\leq2$,
and $\frac{n}{2}-1\leq m\leq n-1$, which yield
\begin{align*}
\sigma_n^2
&\geq n\left(\frac{n}{m+1}-1-\log\frac{n}{m+1}\right)
\geq n\left(\frac{n}{m+1}-1\right)^2\left(\frac{5}{6}-\frac{1}{3}\frac{n}{m+1}\right)\\
&\geq\frac{1}{6}n\left(\frac{n-m-1}{m+1}\right)^2
\geq\frac{1}{6}\frac{(n-m-1)^2}{n}
=\frac{1}{6}\frac{(n-m)^2}{n}\left(\frac{n-m-1}{n-m}\right)^2
\geq\frac{1}{24}\frac{(n-m)^2}{n},
\end{align*}
where at the last inequality we used the fact that the function $(x-1)/x$ is increasing in $x$ with $x=n-m\geq2$. $\blacksquare$

\bigskip

\begin{proposition}\label{p_t_0}
We fix an arbitrary integer $R\geq3$, and define
\begin{equation}\label{uj_t_0}
t_0:=\frac{1}{\sigma_n}\sqrt{\pi R\log\sqrt{m}}.
\end{equation}
There exists a threshold number $m_R$ depending on $R$ such that if $m\geq m_R$, then
\begin{equation}\label{uj_t_0becs1}
t_0\leq\frac{1}{4}\frac{m}{n}\wedge \frac{m^{2/3}}{n},\quad\textrm{ if }m\leq\frac{n}{2}-1.
\end{equation}
There exists a threshold number $n_R$ depending on $R$ such that if $n\geq n_R$ and $m$ is such that
\begin{equation}\label{uj_prop2_felt}
\sqrt{24\pi R}\frac{\sqrt{n\log\sqrt{n}}}{n-m}\leq\frac{1}{8},
\end{equation}
then
\begin{equation}\label{uj_t_0becs2}
t_0\leq\frac{1}{4}\frac{m}{n}\wedge \frac{n^{2/3}}{n-m},\quad\textrm{ if }m\geq\frac{n}{2}.
\end{equation}
\end{proposition}

\bigskip

\noindent {\bf Proof.} First, if $m\leq\frac{n}{2}-1$, then by (\ref{uj_s2mk}) in Proposition \ref{p_sigma_n},
\begin{equation}\label{uj_t_0becskicsi}
t_0\leq \sqrt{20\pi R}\frac{\sqrt{m\log\sqrt{m}}}{n}.
\end{equation}
Thus we see that (\ref{uj_t_0becs1}) holds true if $m$ is greater than some threshold number depending on $R$.

Next, if $m\geq\frac{n}{2}$, then by (\ref{uj_s2mn}) of Proposition \ref{p_sigma_n}, we have
\begin{equation}\label{uj_t_0becsnagy}
t_0\leq \sqrt{24\pi R}\frac{\sqrt{n\log\sqrt{n}}}{n-m},
\end{equation}
where the bounding sequence is less than both $n^{2/3}/(n-m)$ for all large enough $n$, depending on $R$, and since $m\geq\frac{n}{2}$, assumption (\ref{uj_prop2_felt}) implies $t_0\leq\frac{1}{4}\frac{m}{n}$. Therefore we also have (\ref{uj_t_0becs2}).$\blacksquare$

\bigskip
\begin{proposition}\label{p_anj}
For $a_{n,2}$ defined in (\ref{uj_anj_def}) we have
\begin{equation}\label{uj_an2}
\frac{(n-m-1)^3}{3n^2}\leq a_{n,2}\leq\frac{(n-m)^3}{m^2}
\end{equation}
For $a_{n,j}$, with $j=2,3,\ldots$, also defined in (\ref{uj_anj_def}) we have
\begin{equation}\label{uj_a_nj}
a_{n,j}\leq 2^j\frac{n^j}{m^{j-1}},\textrm{ if }m\leq\frac{n}{2}-1,
\quad\textrm{and}\quad
a_{n,j}\leq 2^j\frac{(n-m)^{j+1}}{n^j},\textrm{ if }m\geq\frac{n}{2}-1.
\end{equation}
\end{proposition}

\bigskip

\noindent {\bf Proof.}
Since $a_{n,2}=\sum_{k=m+1}^{n}\frac{(n-k)^2}{k^2}$, the first assertion follows from
$$\frac{(n-m)\left(n-m-\frac{1}{2}\right)(n-m-1)}{3n^2}=\sum_{k=m+1}^{n}\frac{(n-k)^2}{n^2}
\leq\,a_{n,2}\,\leq
(n-m)\frac{(n-m-1)^2}{(m+1)^2}.$$

For any $j=2,3,\ldots$,
\begin{align*}
a_{n,j}&\leq\int_{m}^{n}\left(\frac{n}{k}-1\right)^j\d k
=\int_{m}^{n}\sum_{l=0}^{j}(-1)^l{j\choose l}\left(\frac{n}{k}\right)^{j-l}\d k
\leq\sum_{l=0}^{j}{j\choose l}\int_{m}^{n}\left(\frac{n}{k}\right)^{j}\d k\\
&=\sum_{l=0}^{j}{j\choose l}\frac{1}{j-1}\left(\frac{n^j}{m^{j-1}}-n\right)
\leq2^j\frac{n^j}{m^{j-1}}.
\end{align*}
The approximation applied at the second inequality in the display above is quite gross in the case when $m$ is close to $n$. In this case, with the help of Proposition \ref{p_sigma_n}, one is able to prove the following better upper bound for $a_{n,j}$: if $m\geq\frac{n}{2}-1$, then
\begin{align*}
a_{n,j}=\sum_{k=m+1}^{n}\left(\frac{n-k}{k}\right)^j
\leq\frac{(n-m-1)^{j-1}}{(m+1)^{j-2}n}\sigma_n^2
\leq2\frac{(n-m)^{j+1}}{(m+1)^{j-2}n^2}
\leq2^j\frac{(n-m)^{j+1}}{n^j}.\quad\blacksquare
\end{align*}

\noindent{\bf Proof of Theorem \ref{t_PCexp_1}.} First of all we note that by (\ref{uj_an2}), the assumption $a_{n2}>1$ implies (\ref{uj_prop2_felt}) for all $n$ that is greater than some threshold number depending on $R$. Let $m_R$ and $n_R$ be positive integers at least as big as the threshold numbers given by Proposition \ref{p_t_0}, and such that (\ref{uj_prop2_felt}) holds true for all $m\geq m_R$ and $n\geq n_R$. We fix an integer $R\geq3$, as well as integers $m\geq m_R$ and $n\geq n_R$.

We shall apply Theorem \ref{t_BKN} with the measures $\mu$ and $\nu_R$ and the constant $t_0$ given above. Recalling the decompositions of the characteristic functions corresponding to $\mu$ and $\nu_R$ in (\ref{phifelbontas}) and (\ref{psifelbontas}), we now give upper bounds for the differences $|\exp\{h(\e^{\i t}-1)\}-H_R(\e^{\i t}-1)|$, $|t|\leq t_0$, and $|\phi(t)-\psi(t)|$, $t_0<|t|\leq\pi$, that have the form required by the theorem.

We begin by bounding $|\exp\{h(\e^{\i t}-1)\}-H_R(\e^{\i t}-1)|$ when $|t|\leq t_0$. For an arbitrary such $t$,
\begin{equation}\label{uj_kozel0}
|\exp\{h(\e^{\i t}-1)\}-H_R(\e^{\i t}-1)|\leq\Delta_1+\Delta_2,
\end{equation}
where the $\Delta$s are the errors resulting from the approximations in (\ref{uj_approx}).

Starting with $\Delta_1$, the inequality
\begin{equation*}
|\e^{z_1}-\e^{z_2}|\leq\frac{1}{2}\left(\e^{|z_1|}+\e^{|z_2|}\right)|z_1-z_2|,\quad z_1,z_2\in\C,
\end{equation*}
yields
\begin{equation}\label{uj_delta1_1}
\Delta_1=\left|\exp\{h(w)\}-\exp\{h_R(w)\}\right|
\leq\frac{1}{2}\Big(\e^{|h(w)|}+\e^{|h_R(w)|}\Big)|h(w)-h_R(w)|.
\end{equation}
From the definitions of $h(w)$ in (\ref{uj_h}) and that of $h_R(w)$ in (\ref{uj_h_R}),
\begin{equation*}
|h_R(w)|\vee|h(w)|\,
\leq\,|w|+\frac{|w|^2}{2}+2\sum_{r=3}^{\infty}(a_{n,r}\vee a_{n,2})\frac{|w|^r}{r}
\end{equation*}
If $m\leq\frac{n}{2}-1$, then by (\ref{uj_a_nj}) in Proposition \ref{p_anj}, $(a_{n,r}\vee a_{n,2})\leq2^r\frac{n^r}{m^{r-1}}$, $r=2,3,\ldots$, which implies
\begin{equation}\label{uj_absz_mkicsi}
|h_R(w)|\vee|h(w)|\,\leq\,|w|+\frac{|w|^2}{2}+16\frac{n^3}{m^2}\frac{|w|^3}{3}\sum_{r=3}^{\infty}\left(2\frac{n}{m}|w|\right)^{r-3}
\leq\frac{9}{8}|w|+\frac{32}{3}\frac{n^3}{m^2}|w|^3.
\end{equation}
At the last inequality we used $|w|=|\e^{\i t}-1|\leq |t|\leq t_0\leq\frac{1}{4}\frac{m}{n}$ guaranteed by (\ref{uj_t_0becs1}). If however $m\geq\frac{n}{2}$, then $\frac{n-m-1}{m+1}<1$, thus each term in the defining sum of $a_{n,r}$ decreases as $r$ increases, which means that $a_{n,2}\geq a_{n,3}\geq\ldots$. Therefore by (\ref{uj_a_nj}), for any $r=2,3,\ldots$, $(a_{n,r}\vee a_{n,2})\leq a_{n,2}\leq4\frac{(n-m)^{3}}{n^2}$, and hence
\begin{equation}\label{uj_absz_mnagy}
|h_R(w)|\vee|h(w)|\,\leq\,|w|+\frac{|w|^2}{2}+8\frac{(n-m)^{3}}{n^2}\frac{|w|^3}{3}\sum_{r=3}^{\infty}|w|^{r-3}
\leq\frac{9}{8}|w|+\frac{32}{3}\frac{(n-m)^{3}}{n^2}|w|^3.
\end{equation}
At the last inequality we used $|w|=|\e^{\i t}-1|\leq |t|\leq t_0\leq\frac{1}{4}$, true because of (\ref{uj_t_0becs2}).
Note that by (\ref{uj_absz_mkicsi}) and (\ref{uj_absz_mnagy}), (\ref{uj_t_0becs1}) and (\ref{uj_t_0becs2}) also imply
\begin{equation}\label{uj_absz_kons}
|h_R(w)|\vee|h(w)|\,\leq\,C
\end{equation}
in both cases (that is for all $m$).

If we write (\ref{uj_absz_kons}) back into (\ref{uj_delta1_1}), we obtain
\begin{equation*}
\Delta_1
\leq C|h(w)-h_R(w)|
=C\left|\sum_{r=R+1}^{\infty}\Big(a_{n,r}+(-1)^{r+1}\lfloor a_{n,2}\rfloor\Big)\frac{w^r}{r}\right|
\leq C\sum_{r=R+1}^{\infty}2\Big(a_{n,r}\vee a_{n,2}\Big)\frac{|w|^r}{r}.
\end{equation*}
We proceed by bounding the sum in the last expression. Once again we distinguish two cases, according to the values of $m$.
If $m\leq\frac{n}{2}-1$, then we use $(a_{n,r}\vee a_{n,2})\leq2^r\frac{n^r}{m^{r-1}}$, $r=2,3,\ldots$, thus
$$
|h(w)-h_R(w)|
\leq\frac{2^{R+2}}{R+1}\frac{n^{R+1}}{m^R}|w|^{R+1}\sum_{r=R+1}^{\infty}\left(2\frac{n}{m}|w|\right)^{r-R-1}.
$$
While if $m\geq\frac{n}{2}$, then for any $r=2,3,\ldots$, $(a_{n,r}\vee a_{n,2})\leq a_{n,2}\leq4\frac{(n-m)^{3}}{n^2}$ as noticed before, thus
$$
|h(w)-h_R(w)|
\leq\frac{8}{R+1}\frac{(n-m)^3}{n^2}|w|^{R+1}\sum_{r=R+1}^{\infty}|w|^{r-R-1}.
$$
The last two sums can be bounded with the help of $|w|\leq t_0$ and (\ref{uj_t_0becs1}) or (\ref{uj_t_0becs2}) as seen before, therefore we conclude
\begin{equation}\label{uj_delta1}
\Delta_1\leq\left\{
                   \begin{array}{ll}
                     C_R\frac{n^{R+1}}{m^R}|w|^{R+1}, & \hbox{$m\leq\frac{n}{2}-1$;} \\
                     C_R\frac{(n-m)^3}{n^2}|w|^{R+1}, & \hbox{$m\geq\frac{n}{2}$.}
                   \end{array}
                 \right.
\end{equation}

\bigskip

We now deal with $\Delta_2$. Recalling the definitions in (\ref{uj_h_R}) and (\ref{uj_H_R}), by the series expansion of the exponential function, we have
\begin{align}\nonumber
\Delta_2&=|\exp\{h_R(w)\}-H_{R}(w)|
=\left|\sum_{l=R+1}^{\infty}\frac{h_R^l(w)}{l!}\right|\\ \label{uj_d2}
&\leq\frac{|h_R(w)|^{R+1}}{(R+1)!}\sum_{l=R+1}^{\infty}\frac{|h_R(w)|^{l-R-1}}{(l-R-1)!}
=\frac{|h_R(w)|^{R+1}}{(R+1)!}\exp\{|h_R(w)|\}.
\end{align}

If $m\leq\frac{n}{2}-1$, we use (\ref{uj_absz_kons}) and (\ref{uj_absz_mkicsi}), and then the inequality $(a+b)^k\leq2^k(a^k+b^k)$, $a,b\in{\mathbb R}^+$, $k\in{\mathbb N}$, which yield
\begin{equation*}
\Delta_2
\leq\frac{\e^C}{(R+1)!}\left[\frac{9}{8}|w|+\frac{32}{3}\frac{n^3}{m^2}|w|^3\right]^{R+1}
\leq C_R |w|^{R+1} + C_R\frac{n^{R+1}}{m^R}|w|^{R+1}\left[\frac{n^{2(R+1)}}{m^{R+2}}|w|^{2(R+1)}\right],
\end{equation*}
where
$$
\frac{n^{2(R+1)}}{m^{R+2}}|w|^{2(R+1)}
\leq\frac{n^{2(R+1)}}{m^{R+2}}|t_0|^{2(R+1)}
\leq C_R\frac{n^{2(R+1)}}{m^{R+2}}\!\left(\frac{\sqrt{m\log(\sqrt{m})}}{n}\right)^{2(R+1)}
\!\!\!\!\!=C_R\frac{(\log(\sqrt{m}))^{R+1}}{m}
$$
by (\ref{uj_t_0becskicsi}), and we see that there exists a constant depending on $R$, which bounds the last expression from above for all values of $m$.

If $m\geq\frac{n}{2}$, we use (\ref{uj_absz_kons}) and (\ref{uj_absz_mnagy}) to continue the approximation of $\Delta_2$ in (\ref{uj_d2}). Also applying the inequality $(a+b)^k\leq2^k(a^k+b^k)$, $a,b\in{\mathbb R}^+$, $k\in{\mathbb N}$, we obtain
\begin{align*}
\Delta_2
&\leq\frac{\e^C}{(R+1)!}\left[\frac{9}{8}|w|+\frac{32}{3}\frac{(n-m)^{3}}{n^2}|w|^3\right]^{R+1}\\
&\leq C_R |w|^{R+1} + C_R\frac{(n-m)^{3}}{n^2}|w|^{R+1}\left[\frac{(n-m)^{3R}}{n^{2R}}|w|^{2(R+1)}\right],
\end{align*}
where
\begin{align*}
\frac{(n-m)^{3R}}{n^{2R}}|w|^{2(R+1)}
&\leq \frac{(n-m)^{3R}}{n^{2R}}|t_0|^{2(R+1)}\\
&\leq C_R\frac{(n-m)^{3R}}{n^{2R}}\left(\frac{\sqrt{n\log\left(\sqrt{n}\right)}}{n-m}\right)^{2(R+1)}\\
&= C_R \frac{(n-m)^{R-2}}{n^{R-1}}\left(\log\left(\sqrt{n}\right)\right)^{(R+1)}\\
&\leq C_R \frac{\left(\log\left(\sqrt{n}\right)\right)^{(R+1)}}{n}
\end{align*}
by (\ref{uj_t_0becsnagy}), and again we see that there exists a constant depending on $R$, which bounds the last expression from above for all of $n$. We also note that in the latest bound for $\Delta_2$, the second term is the bigger one due to our assumption $a_{n2}>1$ and (\ref{uj_an2}).

These considerations lead to
\begin{equation*}
\Delta_2
\leq\left\{
      \begin{array}{ll}
        C_R\frac{n^{R+1}}{m^R}|w|^{R+1}, & \hbox{$m\leq\frac{n}{2}-1$;} \\
        C_R\frac{(n-m)^{3}}{n^2}|w|^{R+1}, & \hbox{$m\geq\frac{n}{2}$,}
      \end{array}
    \right.
\end{equation*}
which, up to a constant factor, is exactly the same upper bound we obtained for $\Delta_1$ in (\ref{uj_delta1}). Hence by (\ref{uj_kozel0}) and $|w|\leq|t|$,
\begin{equation}\label{uj_kozel}
|\exp\{h(\e^{\i t}-1)\}-H_R(\e^{\i t}-1)|\leq\left\{
      \begin{array}{ll}
        C_R\frac{n^{R+1}}{m^R}|t|^{R+1}, & \hbox{$m\leq\frac{n}{2}-1$;} \\
        C_R\frac{(n-m)^{3}}{n^2}|t|^{R+1}, & \hbox{$m\geq\frac{n}{2}$,}
      \end{array}
    \right.\quad|t|\leq t_0
\end{equation}


\bigskip

Next, from (\ref{uj_chi}), by the application of the inequality $1-\cos t\geq\frac{2}{\pi^2}t^2$, $0\leq t\leq\pi$, we see that
\begin{equation}\label{uj_chibecs}
|\chi(t)|=\exp\left\{-\sigma_n^2(1-\cos t)\right\}\leq\exp\left\{-\frac{2\sigma_n^2}{\pi^2}t^2\right\},\quad0\leq t\leq\pi.
\end{equation}

\bigskip

The next step in our proof is to bound $|\psi(t)|$. We recall the decomposition $\psi(t)=\chi(t)H_R(\e^{\i t}-1)$ in (\ref{psifelbontas}), and start by examining $|H_R(w)|$ (where $w=\e^{\i t}-1$) defined in (\ref{uj_H_R}) and (\ref{uj_h_R}). If we look at (\ref{uj_h_R}), we see that with the application of the triangle inequality and $|w|\leq|t|$, $|h_R(w)|$ can be bounded from above by a polynomial of $|t|$ of degree $R$. The coefficients in this polynomial are less than 1 for $|t|^0$ and $|t|$, and their order is given by (\ref{uj_a_nj}) for $|t|^r$, $r=3,\ldots,R$. If $m\leq\frac{n}{2}-1$, these orders are $n^r/m^{r-1}$, $r=3,\ldots,R$, respectively, and hence
$$|h_R(w)|\leq C+C_R\left(\frac{n^3}{m^2}|t|^3I_{\left[0,\frac{m}{n}\right]}(t)+\frac{n^R}{m^{R-1}}|t|^RI_{\left(\frac{m}{n},\pi\right]}(t)\right),$$
where for any $A\subset{\mathbb R}$, $I_A(t)$ is 1 if $|t|\in A$ and 0 otherwise. If $m\geq\frac{n}{2}$ the coefficient orders in the polynomial are $(n-m)^{3}/n^2$ for all $r=3,\ldots,R$, so
$$|h_R(w)|\leq C+C_R\left(\frac{(n-m)^3}{n^2}|t|^3I_{\left[0,1\right]}(t)+\frac{(n-m)^{3}}{n^{2}}|t|^RI_{\left(1,\pi\right]}(t)\right).$$

By (\ref{uj_H_R}), these bounds imply
\begin{align*}
|H_R(w)|
&\leq C_R+C_R\left(
\frac{n^3}{m^2}|t|^3I_{\left[0,\frac{m^{2/3}}{n}\right]}(t)
+\frac{n^{3R}}{m^{2R}}|t|^{3R}I_{\left(\frac{m^{2/3}}{n},\frac{m}{n}\right]}(t)
+\frac{n^{R^2}}{m^{R^2-R}}|t|^{R^2}I_{\left(\frac{m}{n},\pi\right]}(t)\right),
\end{align*}
if $m\leq\frac{n}{2}-1$, and
\begin{align*}
|H_R(w)|
&\leq C_R+C_R\left(
\frac{(n-m)^{3}}{n^{2}}|t|^3I_{\left[0,\frac{n^{2/3}}{n-m}\wedge1\right]}(t)
+\frac{(n-m)^{3R}}{n^{2R}}|t|^{3R}I_{\left(\frac{n^{2/3}}{n-m}\wedge1,1\right)}(t)+\right.\\ \nonumber
&\left.\quad\quad\quad\quad\quad\quad\quad
+\frac{(n-m)^{3}}{n^{2}}|t|^RI_{\left(1,\frac{n^{2/3}}{n-m}\vee 1\right)}(t)
+\frac{(n-m)^{3R}}{n^{2R}}|t|^{R^2}I_{\left[\frac{n^{2/3}}{n-m}\vee 1,\pi\right]}(t)
\right),
\end{align*}
if $m\geq\frac{n}{2}$.

We introduce a new variable $x_{nt}:=\sigma_n^2t^2$. With this, for $m\leq\frac{n}{2}-1$, by $\frac{1}{20}\frac{n^2}{m}t^2\leq\sigma_n^2t^2$ from (\ref{uj_s2mk}), we get
\begin{align}\nonumber
|H_R(w)|&\leq C_R+C_R\left(
x_{nt}^{\frac{3}{2}}I_{\left[0,\frac{m^{2/3}}{n}\right]}(t)
+x_{nt}^{\frac{3R}{2}}I_{\left(\frac{m^{2/3}}{n},\frac{m}{n}\right]}(t)
+x_{nt}^{\frac{R^2}{2}}I_{\left(\frac{m}{n},\pi\right]}(t)\right)\\ \label{uj_f_mkicsi}
&=:p(x_{nt}),\quad\textrm{ if }m\leq\frac{n}{2}-1;
\end{align}
and for $m\geq\frac{n}{2}$, with the help of $\frac{1}{24}\frac{(n-m)^2}{n}t^2\leq\sigma_n^2t^2$ from (\ref{uj_s2mk}), we obtain
\begin{align}\nonumber
&|H_R(w)|\\ \nonumber
&\leq C_R+C_R\left(
x_{nt}^{\frac{3}{2}}I_{\left[0,\frac{n^{2/3}}{n-m}\wedge1\right]}(t)
+x_{nt}^{\frac{3R}{2}}I_{\left(\frac{n^{2/3}}{n-m}\wedge1,1\right)}(t)
+x_{nt}^{\frac{3}{2}}I_{\left(1,\frac{n^{2/3}}{n-m}\vee 1\right)}(t)
+x_{nt}^{\frac{3R}{2}}I_{\left[\frac{n^{2/3}}{n-m}\vee 1,\pi\right]}(t)\right)\\ \label{uj_f_mnagy}
&=:p(x_{nt}),\quad\textrm{ if }m\geq\frac{n}{2}
\end{align}

Thus in either case, by (\ref{uj_chibecs}), we have
\begin{equation*}
|\psi(t)|\leq\exp\left\{-\frac{1}{\pi^2}x_{nt}\right\}\left(\exp\left\{-\frac{1}{\pi^2}x_{nt}\right\}p(x_{nt})\right).
\end{equation*}
The second term in the product above is clearly bounded from above by a constant depending on $R$ for all $x_{nt}$, and thus for all $0\leq|t|\leq\pi$ and all values of $\sigma_n$. The first term $\exp\left\{-\frac{1}{\pi^2}\sigma_n^2t^2\right\}$ is a decreasing function of $t$. Putting these together, and then using the definition of $t_0$ from (\ref{uj_t_0}) yields
\begin{equation}\label{kek_psi}
|\psi(t)|\leq C_R\exp\left\{-\frac{1}{\pi^2}\sigma_n^2t_0^2\right\}
=C_R\left(\frac{1}{\sqrt{m}}\right)^{R},\quad t_0\leq|t|\leq\pi.
\end{equation}

\bigskip

We proceed by bounding $|\phi(t)|$. It is easy to calculate
\begin{equation*}
|\phi(t)|=\left|\exp\left\{-\sum_{k=m+1}^{n}\log\left(1-\frac{n-k}{k}(\e^{\i t}-1)\right)\right\}\right|
=\prod_{k=m+1}^{n}\frac{1}{\sqrt{1+2\frac{n(n-k)}{k^2}(1-\cos t)}}.
\end{equation*}
We note here that $|\phi(t)|$ is obviously a decreasing function on the whole interval $[0,\pi]$, hence $|\phi(t)|\leq|\phi(t_0)|$ for all $t_0\leq|t|\leq\pi$.

Now for any $k=\{m+1,\ldots,n-1\}$, by the right hand side of the inequality
\begin{equation}\label{uj_cos}
\frac{t^2}{2}-\frac{t^4}{4!}\leq1-\cos t\leq\frac{t^2}{2}-\frac{t^4}{30},\quad 0\leq t\leq1,
\end{equation}
and (\ref{uj_t_0becs1}) or (\ref{uj_t_0becs2}), we get
$$x_{k,t_0}:=2\frac{n(n-k)}{k^2}(1-\cos t_0)\leq\frac{n(n-m-1)}{(m+1)^2}t_0^2<1.$$
Therefore we can apply the inequality
$$\frac{1}{1+x}\leq\exp\left\{-x+\frac{x^2}{2}\right\},\quad 0\leq x\leq1,$$
with $x=x_{k,t_0}$ to obtain
\begin{align*}
|\phi(t_0)|\leq\exp\left\{-\sigma_n^2(1-\cos t_0)+b_{n}(1-\cos t_0)^2\right\},
\end{align*}
where $b_n=n^2\sum_{k=m+1}^{n}\frac{(n-k)^2}{k^4}$.
Next, applying (\ref{uj_cos}) gives
\begin{equation}\label{kek_koztes}
|\phi(t_0)|\leq\exp\left\{-\sigma_n^2\left(\frac{t_0^2}{2}-\frac{t_0^4}{24}\right)+b_n\left(\frac{t_0^2}{2}-\frac{t_0^4}{30}\right)^2\right\}
\leq\exp\left\{-\sigma_n^2t_0^2\left(\frac{1}{2}-\frac{1}{24}-\frac{b_nt_0^2}{4\sigma_n^2}\right)\right\}.
\end{equation}

Now we shall prove $b_n t_0^2\leq\sigma_n^2$.
Since the terms in $b_n=n^2\sum_{k=m+1}^{n}\frac{(n-k)^2}{k^4}$ decrease as $k$ increases, with our usual technique we obtain
\begin{align*}
b_n&\leq n^2\int_{m}^{n}\frac{(n-x)^2}{x^4}\d x
=n^4\int_{m}^{n}\frac{1}{x^4}\d x-2n^3\int_{m}^{n}\frac{1}{x^3}\d x+n^2\int_{m}^{n}\frac{1}{x^2}\d x\\
&=\frac{1}{3}n^4\left(\frac{1}{m^3}-\frac{1}{n^3}\right)-n^3\left(\frac{1}{m^2}-\frac{1}{n^2}\right)+n^2\left(\frac{1}{m}-\frac{1}{n}\right)
=\frac{1}{3}n\left(\frac{n-m}{m}\right)^3.
\end{align*}
By (\ref{uj_t_0becs1}), (\ref{uj_t_0becs2}) and the bound above for $b_n$,
$$
b_n t_0^2
\leq\frac{1}{48}n\left(\frac{n-m}{m}\right)^3\left(\frac{m}{n}\right)^{2}
=\frac{1}{48}\frac{(n-m)^3}{mn},
$$
and the latter expression is at most $\sigma_n^2$ by (\ref{uj_s2mk}) and $n-m\leq n$ if $m\leq\frac{n}{2}-1$, and (\ref{uj_s2mn}) and $(n-m)/m\leq1$ if $m\geq\frac{n}{2}$. Therefore we have $b_n t_0^2\leq\sigma_n^2$, and substituting this into (\ref{kek_koztes}) yields
\begin{equation}\label{kek_phi}
|\phi(t)|\leq
C\exp\left\{-\frac{1}{\pi^2}\sigma_n^2t_0^2\right\}
=C\left(\frac{1}{\sqrt{m}}\right)^{R},\quad t_0\leq|t|\leq\pi.
\end{equation}

\bigskip

Now with the help of (\ref{kek_phi}) and (\ref{kek_psi}), it is easy to bound the difference $|\phi(t)-\psi(t)|$ at the points $t_0<|t|\leq\pi$, which is necessary for the application of the Theorem \ref{t_BKN}:
\begin{equation*}
|\phi(t)-\psi(t)|\leq|\phi(t)|+|\psi(t)|\leq C_R\left(\frac{1}{\sqrt{m}}\right)^{R},\quad t_0\leq|t|\leq\pi.
\end{equation*}
If $m\geq\frac{n}{2}$, then
$$
\left(\frac{1}{\sqrt{m}}\right)^{R}
\leq\left(\frac{\sqrt{2}}{\sqrt{n}}\right)^{R}
=(\sqrt{n})^{R-2}\left(\frac{\sqrt{2}}{n}\right)^{R-1}
\leq C_R\frac{(\sqrt{n})^{R-2}}{(n-m)^{R-1}}.
$$
Therefore we have
\begin{equation}\label{uj_messze}
|\phi(t)-\psi(t)|\leq
\left\{
  \begin{array}{ll}
    C_R\left(\frac{1}{\sqrt{m}}\right)^{R}, & \hbox{$m\leq\frac{n}{2}-1$;} \\
    C_R\frac{(\sqrt{n})^{R-2}}{(n-m)^{R-1}}, & \hbox{$m\geq\frac{n}{2}$,}
  \end{array}\quad t_0\leq|t|\leq\pi.
\right.
\end{equation}

\bigskip

Now we apply Theorem \ref{t_BKN}. Condition (\ref{uj_tetel1}) is given by (\ref{uj_kozel}) and (\ref{uj_chibecs}), thus
$S=1$, $\gamma_0=0$, $\gamma_1=C_R\frac{n^{R+1}}{m^R}$ if $m\leq\frac{n}{2}-1$ and $\gamma_1=C_R\frac{(n-m)^3}{n^2}$ if $m\geq\frac{n}{2}$, $\theta_1=R+1$, $\gamma=1$ and $\rho=\frac{2}{\pi^2}\sigma_n^2$. Condition (\ref{uj_tetel2}) is given by (\ref{uj_messze}), so $\eta=C_R\left(\frac{1}{\sqrt{m}}\right)^{R}$ if $m\leq\frac{n}{2}-1$ and $\eta=C_R\frac{(\sqrt{n})^{R-2}}{(n-m)^{R-1}}$ if $m\geq\frac{n}{2}$. We also use Proposition \ref{p_sigma_n} to bound $\sigma_n$. It follows that
\begin{equation}\label{uj_bizvege_loc}
\sup_{k\in\Z}|\mu\{k\}-\nu_R\{k\}|\leq
\left\{
  \begin{array}{ll}
    C_R\left(\frac{1}{\sqrt{m}}\right)^{R}, & \hbox{$m\leq\frac{n}{2}-1$;} \\
    C_R\frac{(\sqrt{n})^{R-2}}{(n-m)^{R-1}}, & \hbox{$m\geq\frac{n}{2}$,}
  \end{array}
\right.
\end{equation}
which is exactly what we wanted to prove. $\blacksquare$

\bigskip

We would like to prove an analogous result of Theorem \ref{t_PCexp_2} in the case when $a_{n,2}<1$, hence from now on we assume that this inequality holds true. In this case $\lfloor a_{n,2}\rfloor=0$, so the characteristic function in (\ref{uj_phi}) has the form
\begin{equation*}
\phi(t)
=\exp\left\{-\sum_{k=m+1}^{n}\log\left(1-\frac{n-k}{k}(\e^{\i t}-1)\right)\right\}.
\end{equation*}
Now we also assume $(n-m)/n\leq1/(4\pi)$, which implies
$$
\frac{n-k}{k}|w|
=\frac{n-k}{k}|\e^{\i t}-1|
\leq\frac{n-m-1}{m+1}|\e^{\i t}-1|
\leq\frac{n-m-1}{m+1}|t|
\leq1,
$$
for all $0\leq t\leq\pi$. This allows us to expand the logarithmic expression in $\phi(t)$, therefore with $w=\e^{\i t}-1$ and the notation introduced in (\ref{uj_anj_def}),
\begin{equation*}
\phi(t)
=\exp\left\{\sum_{r=1}^{\infty}a_{n,r}\frac{w^r}{r}\right\},\quad |t|\leq\pi.
\end{equation*}
We note that $a_{n,1}=\sigma_n^2-a_{n,2}$, hence the line above can be rewritten as
\begin{equation}\label{nn_phifelbontas}
\phi(t)=\chi(t)\exp\{h(w)\},\quad |t|\leq\pi,
\end{equation}
where $\chi(t)$ is the characteristic function of the Poisson distribution with parameter $\sigma_n^2$ given in (\ref{uj_chi}), and
\begin{equation}\label{nn_h}
h(w)=-a_{n,2}w+a_{n,2}\frac{w^2}{2}
+\sum_{r=3}^{\infty}a_{n,r}\frac{w^r}{r}.
\end{equation}

Now we fix an integer $R\geq3$, and following the argument in (\ref{uj_approx}), we modify the function $\exp\{h(w)\}$:
\begin{equation}\label{nn_approx}
\exp\{h(w)\}\approx\exp\{h_R(w)\}\approx H_{R}(w),
\end{equation}
where
\begin{equation}\label{nn_h_R}
h_R(w)=-a_{n,2}w+a_{n,2}\frac{w^2}{2}
+\sum_{r=3}^{R}a_{n,r}\frac{w^r}{r}
\end{equation}
and
\begin{equation}\label{nn_H_R}
H_{R}(w)=\sum_{l=0}^{3R-2}\frac{h_R^l(w)}{l!}.
\end{equation}

We approximate the distribution $\mu={\cal D}(\widetilde{W}_{n,m}+c)$ with $\nu_R$, which we define to be the finite signed measure on the nonnegative integers whose characteristic function is
\begin{equation}\label{nn_psifelbontas}
\psi(t)=\chi(t)H_{R}(\e^{\i t}-1),\quad t\in{\mathbb R},
\end{equation}
repeating (\ref{psifelbontas}). But now $H_{R}(\e^{\i t}-1)$ is a polynomial of $\e^{\i t}-1$ of degree $3R^2-R$, and the corresponding Poisson--Charlier signed measure $\nu_R=\nu_R(\sigma_n^2,\widetilde{a}_{n,m}^{(1)},\ldots,\widetilde{a}_{n,m}^{(3R^2-R)})$ is defined by
\begin{equation}\label{uj_nu_R2}
\nu_R\{j\}=\mathrm{Po}(\sigma_n^2)\{j\}\left(1+\sum_{r=1}^{3R^2-R}(-1)^{r+1}\widetilde{a}_{n,m}^{(r)}C_r(j,\sigma_n^2)\right),\quad j\in{\mathbb N},
\end{equation}
where $\widetilde{a}_{n,m}^{(r)}$ is the coefficient of $(\e^{\i t}-1)^r$ in $H_{R}(\e^{\i t}-1)$, and $C_r(j,\sigma_n^2)$ is the $r$-th Charlier polynomial defined in (\ref{Charlier_polynomial}).

\bigskip

\begin{theorem}\label{t_PCexp_2}
We assume $a_{n,2}<1$. For an arbitrary integer $R\geq3$ there exist a threshold number $n_R$ depending on $R$ such that if $n\geq n_R$, then
\begin{equation*}
\sup_{k\in\Z}|\mu\{k\}-\nu_R\{k\}|\leq C_R \frac{1}{(\sqrt{n})^{R}}.
\end{equation*}
\end{theorem}

\bigskip


\noindent{\bf Proof.} By (\ref{uj_an2}), condition $a_{n,2}<1$ implies $\frac{n-m}{n}<\frac{1}{4\pi}$ for all $n$ greater than some threshold number $n_R$ depending on $R$. We fix integers $R\geq3$ and $n\geq n_R$.

As it was the case for Theorem \ref{t_PCexp_1}, the key of the proof is Theorem \ref{t_BKN}. We apply it with the measures $\mu$ and $\nu_R$ defined above and $t_0=\pi$. Recalling the decompositions of the characteristic functions corresponding to $\mu$ and $\nu_R$ in (\ref{nn_phifelbontas}) and (\ref{nn_psifelbontas}), we now give an upper bound for the difference $|\exp\{h(\e^{\i t}-1)\}-H_R(\e^{\i t}-1)|$, $|t|\leq\pi$. For an arbitrary such $t$,
\begin{equation}\label{nn_kozel0}
|\exp\{h(\e^{\i t}-1)\}-H_R(\e^{\i t}-1)|\leq\Delta_1+\Delta_2,
\end{equation}
where the $\Delta$s are the errors resulting from the approximations in (\ref{nn_approx}).

Regarding $\Delta_1$, we apply the inequality
\begin{equation*}
|\e^{z_1}-\e^{z_2}|\leq\frac{1}{2}\left(\e^{|z_1|}+\e^{|z_2|}\right)|z_1-z_2|,\quad z_1,z_2\in\C,
\end{equation*}
and obtain
\begin{equation}\label{nn_delta1_1}
\Delta_1=\left|\exp\{h(w)\}-\exp\{h_R(w)\}\right|
\leq\frac{1}{2}\Big(\e^{|h(w)|}+\e^{|h_R(w)|}\Big)|h(w)-h_R(w)|.
\end{equation}
From the new definitions of $h(w)$ in (\ref{nn_h}) and that of $h_R(w)$ in (\ref{nn_h_R}),
\begin{equation*}
|h_R(w)|\vee|h(w)|\,
\leq\,a_{n,2}|w|+a_{n,2}\frac{|w|^2}{2}+\sum_{r=3}^{\infty}a_{n,r}\frac{|w|^r}{r}
\end{equation*}
By (\ref{uj_a_nj}),
\begin{align}\nonumber
|h_R(w)|\vee|h(w)|\,&\leq\,\frac{(n-m)^3}{n^2}|w|+\frac{(n-m)^3}{n^2}\frac{|w|^2}{2}+8\frac{(n-m)^{4}}{n^3}\frac{|w|^3}{3}\sum_{r=3}^{\infty}\left(2\frac{n-m}{n}|w|\right)^{r-3}\\ \nonumber
&\leq\left(1+\frac{\pi}{2}\right)\frac{(n-m)^3}{n^2}|w|+\frac{32}{3}\frac{(n-m)^{4}}{n^3}|w|^3\\ \label{nn_absz}
&\leq\left(1+\frac{\pi}{2}+\frac{32\pi^2}{3}\right)\frac{(n-m)^3}{n^2}|w|.
\end{align}
At the second inequality we used $(n-m)/n<1/(4\pi)$ and $|w|=|\e^{\i t}-1|\leq |t|\leq\pi$.
Note that the latter inequality and $a_{n,2}<1$ together with (\ref{uj_an2}) imply
$$
\frac{(n-m)^3}{n^2}|w|=\left(\frac{n-m}{n^{2/3}}\right)^3|w|\leq C,
$$
thus we also have
\begin{equation}\label{nn_absz_kons}
|h_R(w)|\vee|h(w)|\,\leq\,C.
\end{equation}

If we write (\ref{nn_absz_kons}) back into (\ref{nn_delta1_1}), we obtain
\begin{align*}
\Delta_1
&\leq C|h(w)-h_R(w)|
=C\left|\sum_{r=R+1}^{\infty}a_{n,r}\frac{w^r}{r}\right|
\leq C\sum_{r=R+1}^{\infty}a_{n,r}\frac{|w|^r}{r}\\
&\leq\frac{2^{R+1}}{R+1}\frac{(n-m)^{R+2}}{n^{R+1}}|w|^{R+1}\sum_{r=R+1}^{\infty}\left(2\frac{n-m}{n}|w|\right)^{r-R-1}.
\end{align*}
Again, we see from $|w|\leq\pi$ and $(n-m)/n<1/(4\pi)$ that the last sum above is finite for all $|t|\leq\pi$. We conclude
\begin{equation}\label{nn_delta1}
\Delta_1\leq C_R\frac{(n-m)^{R+2}}{n^{R+1}}|w|^{R+1}.
\end{equation}


We now deal with $\Delta_2$. Recalling the definitions in (\ref{nn_h_R}) and (\ref{nn_H_R}), by the series expansion of the exponential function, we have
\begin{align*}
\Delta_2&=|\exp\{h_R(w)\}-H_{3R-2}(w)|
=\left|\sum_{l=3R-1}^{\infty}\frac{h_R^l(w)}{l!}\right|\\
&\leq\frac{|h_R(w)|^{3R-1}}{(3R-1)!}\sum_{l=3R-1}^{\infty}\frac{|h_R(w)|^{l-3R+1}}{(l-3R+1)!}
=\frac{|h_R(w)|^{3R-1}}{(3R-1)!}\exp\{|h_R(w)|\}.
\end{align*}
Now $\exp\{|h_R(w)|\}\leq C$ by (\ref{nn_absz_kons}), and
$$
|h_R(w)|^{3R-1}
\leq C_R\left(\frac{(n-m)^{3}}{n^2}|w|\right)^{3R-1}
= C_R \frac{(n-m)^{3R}}{n^{2R}}|w|^{3R-1}\left(\frac{n-m}{n^{2/3}}\right)^{3(2R-1)},
$$
where the last fraction is less than some constant depending on $R$ by $a_{n,2}<1$ and (\ref{uj_an2}). Therefore we conclude
\begin{equation*}
\Delta_2\leq C_R \frac{(n-m)^{3R}}{n^{2R}}|w|^{3R-1}.
\end{equation*}

We substitute the inequality above and (\ref{nn_delta1}) into (\ref{nn_kozel0}), thus with $|w|\leq|t|$, we obtain
$$
|\exp\{h(\e^{\i t}-1)\}-H_R(\e^{\i t}-1)|\leq C_R\frac{(n-m)^{R+2}}{n^{R+1}}|t|^{R+1}+C_R \frac{(n-m)^{3R}}{n^{2R}}|t|^{3R-1}.
$$

Now we apply Theorem \ref{t_BKN}. Condition (\ref{uj_tetel1}) with $t_0=\pi$ is given by the last inequality and (\ref{uj_chibecs}), which is also true in this case, thus $S=2$, $\gamma_0=0$, $\gamma_1=C_R\frac{(n-m)^{R+2}}{n^{R+1}}$, $\gamma_2=C_R\frac{(n-m)^{3R}}{n^{2R}}$, $\theta_1=R+1$, $\theta_2=3R-1$, $\gamma=1$ and $\rho=\frac{2}{\pi^2}\sigma_n^2$. We also use Proposition \ref{p_sigma_n} to bound $\sigma_n$. It follows that
$$\sup_{k\in\Z}|\mu\{k\}-\nu_R\{k\}|\leq C_R \frac{1}{(\sqrt{n})^{R}},$$
and the proof of Theorem \ref{t_PCexp_2} is complete. $\blacksquare$

\bigskip

As outlined at the beginning of the chapter, our goal is to estimate the total variation distance of $\mu$ and $\nu_R$. Formula (\ref{uj_BKN_d_TV}) in Theorem \ref{t_BKN} also provides us a way to derive total variation error bounds form the local error bounds given by Theorems \ref{t_PCexp_1} and \ref{t_PCexp_2}.

But first, we would like to determine the difference in total variation between the approximating measures for successive values of $R$. By Lemma 6.~in \cite{BC}, we know that the total variation norm of $C_r(\cdot,\sigma_n^2)\mathrm{Po}(\sigma_n^2)$ is less than $\left((2r)/(\e\sigma_n^2)\right)^{r/2}$. With the help of this and Propositions \ref{p_anj} and \ref{p_sigma_n}, one can derive that the term that has the greatest total variation norm in the sum that defines $\nu_{R+1}-\nu_R$ by (\ref{uj_nu_R}), is the one belonging to $r=R+1$, thus
\begin{equation*}
||\nu_{R+1}-\nu_R||_{TV}\leq C_R\,\widetilde{a}^{R+1}_{n,m}\left(\frac{1}{\sigma_n}\right)^{R+1}.
\end{equation*}
From Proposition \ref{p_anj}, one can also deduce that the order of $\widetilde{a}^{R+1}_{n,m}$ is $n^{R+1}/m^R$, $(n-m)^3/n^2$ and $(n-m)^{R+2}/n^{R+1}$ respectively, in the three cases indicated in the display below. Hence, by Proposition \ref{p_sigma_n}, we have
\begin{equation}\label{uj_totnorm}
||\nu_{R+1}-\nu_R||_{TV}\leq
\left\{
  \begin{array}{ll}
    C_R\frac{1}{(\sqrt{m})^{R-1}}, & \hbox{in the case of Theorem \ref{t_PCexp_1} when $m\leq\frac{n}{2}-1$;} \\
    C_R\frac{(\sqrt{n})^{R-3}}{(n-m)^{R-2}}, & \hbox{in the case of Theorem \ref{t_PCexp_1} when $m\geq\frac{n}{2}$;}\\
    C_R\frac{n-m}{(\sqrt{n})^{R+1}}, & \hbox{in the case of Theorem \ref{t_PCexp_2}.}
  \end{array}
\right.
\end{equation}


\bigskip

\begin{corollary}\label{kov_CPexpansion}
For all $n$ and $m$ for which Theorem \ref{t_PCexp_1} is valid we have
\begin{equation*}
d_{\mathrm{TV}}(\mu,\nu_R)\leq C_R\sigma_n\log\sigma_n\left(\frac{1}{\sqrt{m}}\right)^{R},\quad\textrm{if }m\leq\frac{n}{2}-1,
\end{equation*}
and
\begin{equation*}
d_{\mathrm{TV}}(\mu,\nu_R)\leq C_R\frac{(\sqrt{n})^{R-3}}{(n-m)^{R-2}},\quad\textrm{if }m\geq\frac{n}{2},
\end{equation*}
with $\nu_R$ defined in (\ref{uj_nu_R}).
For all $n$ and $m$ for which Theorem \ref{t_PCexp_2} is valid we have
\begin{equation*}
d_{\mathrm{TV}}(\mu,\nu_R)\leq C_R\frac{n-m}{(\sqrt{n})^{R+1}},
\end{equation*}
with $\nu_R$ defined in (\ref{uj_nu_R2}).
\end{corollary}

\bigskip

\noindent{\bf Proof.}
First we assume that $n$ and $m$ satisfy the conditions of Theorem \ref{t_PCexp_1}. It follows from (\ref{uj_BKN_epsilon}) and the end of the proof of Theorem \ref{t_PCexp_1} that
\begin{equation}\label{uj_bizvege_epsilonab}
\varepsilon_{ab}\leq
\left\{
  \begin{array}{ll}
    C_R\left(\frac{1}{\sqrt{m}}\right)^{R-1}+C_R(b-a+2)\left(\frac{1}{\sqrt{m}}\right)^{R}, & \hbox{if $m\leq\frac{n}{2}-1$;} \\
    C_R\frac{(\sqrt{n})^{R-3}}{(n-m)^{R-2}}+C_R(b-a+2)\frac{(\sqrt{n})^{R-2}}{(n-m)^{R-1}}, & \hbox{if $m\geq\frac{n}{2}$.}
  \end{array}
\right.
\end{equation}

In order to obtain a total variation bound by (\ref{uj_BKN_d_TV}), we need to be able to control the tails
of the approximating measure $\nu_R$. As in \cite{BKN} p.~9, it can be deduced from the Chernoff inequalities for Po$(\sigma_n^2)$ that
\begin{equation}\label{chernoff_a}
|\nu_R|\{[0,a)\}\leq R^2\max_{1\leq r\leq R^2}\{\widetilde{a}_{n,m}^{(r)}\}\exp\left\{-\frac{(\sigma_n^2-a)^2}{3\sigma_n^2}\right\},\quad
0\leq a\leq\sigma_n^2
\end{equation}
and
\begin{equation}\label{chernoff_b}
|\nu_R|\{(b,\infty)\}\leq
R^2\max_{1\leq r\leq R^2}\{\widetilde{a}_{n,m}^{(r)}\}\exp\left\{-\frac{(b-R^2-\sigma_n^2)^2}{3\sigma_n^2}\right\},\quad
\sigma_n^2+R^2\leq b \leq2\sigma_n^2.
\end{equation}
If $m\leq\frac{n}{2}-1$, then by (\ref{uj_a_nj}),
$$\max_{1\leq r\leq R^2}\{\widetilde{a}_{n,m}^{(r)}\}\leq C_R\left(\frac{n^R}{m^{R-1}}\right)^R$$
and we choose
\begin{equation*}
a=\sigma_n^2-\sqrt{3\sigma_n^2\log\left(\left(\frac{n^R}{m^{R-1}}\right)^R(\sqrt{m})^{R-1}\right)}
\end{equation*}
and
\begin{equation*}
b=\sigma_n^2+R^2+\sqrt{3\sigma_n^2\log\left(\left(\frac{n^R}{m^{R-1}}\right)^R(\sqrt{m})^{R-1}\right)}.
\end{equation*}
By (\ref{uj_s2mk}), the argument of the logarithmic expressions in $a$ and $b$ are less than $\sigma_n^{2R^2}$, hence for all large enough $\sigma_n$, $0\leq a\leq\sigma_n^2$ and $\sigma_n^2+R^2\leq b \leq2\sigma_n^2$.
In the case when $m\geq\frac{n}{2}$, we have
$$\max_{1\leq r\leq R^2}\{\widetilde{a}_{n,m}^{(r)}\}\leq C_R\left(\frac{(n-m)^3}{m^2}\right)^R$$
by (\ref{uj_a_nj}), and we put
\begin{equation*}
a=\sigma_n^2-\sqrt{3\sigma_n^2\log\left(\left(\frac{(n-m)^3}{m^2}\right)^R\frac{(n-m)^{R-2}}{(\sqrt{n})^{R-3}}\right)}
\end{equation*}
and
\begin{equation*}
b=\sigma_n^2+R^2+\sqrt{3\sigma_n^2\log\left(\left(\frac{(n-m)^3}{m^2}\right)^R\frac{(n-m)^{R-2}}{(\sqrt{n})^{R-3}}\right)}.
\end{equation*}
Again, we see form (\ref{uj_s2mk}) that the argument of the logarithmic expressions in $a$ and $b$ are less than $\sigma_n^{2R-1}$, hence for all large enough $\sigma_n$, $0\leq a\leq\sigma_n^2$ and $\sigma_n^2+R^2\leq b \leq2\sigma_n^2$.
Therefore, in both cases one can apply the inequalities (\ref{chernoff_a}) and (\ref{chernoff_b}), which yield
\begin{equation*}
|\nu_R|\{[0,a)\cup(b,\infty)\}\leq
\left\{
  \begin{array}{ll}
    C_R\left(\frac{1}{\sqrt{m}}\right)^{R-1}, & \hbox{if $m\leq\frac{n}{2}-1$;} \\
    C_R\frac{(\sqrt{n})^{R-3}}{(n-m)^{R-2}}, & \hbox{if $m\geq\frac{n}{2}$.}
  \end{array}
\right.
\end{equation*}
Combining this with (\ref{uj_bizvege_loc}) and (\ref{uj_bizvege_epsilonab}) in (\ref{uj_BKN_d_TV}), with the choices of $a$ and $b$ given above, one can deduce
\begin{equation*}
d_{\mathrm{TV}}(\mu,\nu_R)\leq
\left\{
  \begin{array}{ll}
    C_R\left(\frac{1}{\sqrt{m}}\right)^{R-1}+C_R\sigma_n\log\sigma_n\left(\frac{1}{\sqrt{m}}\right)^{R}, & \hbox{if $m\leq\frac{n}{2}-1$;} \\
    C_R\frac{(\sqrt{n})^{R-3}}{(n-m)^{R-2}}+C_R\sigma_n\log\sigma_n\frac{(\sqrt{n})^{R-2}}{(n-m)^{R-1}}, & \hbox{if $m\geq\frac{n}{2}$.}
  \end{array}
\right.
\end{equation*}
Together with Proposition \ref{p_sigma_n} this gives
\begin{equation}\label{uj_bizvege_dTV0}
d_{\mathrm{TV}}(\mu,\nu_R)\leq
\left\{
  \begin{array}{ll}
    C_R\sigma_n\log\sigma_n\left(\frac{1}{\sqrt{m}}\right)^{R}, & \hbox{if $m\leq\frac{n}{2}-1$;} \\
    C_R\log\sigma_n\frac{(\sqrt{n})^{R-3}}{(n-m)^{R-2}}, & \hbox{if $m\geq\frac{n}{2}$.}
  \end{array}
\right.
\end{equation}

We see that in the case when $m\leq\frac{n}{2}-1$, we already have the inequality we aimed for. In the latter case, when $m\geq\frac{n}{2}$, one can omit the $\log\sigma_n$ factor in the bound above with the help of the following argument. Note that
$$
d_{\mathrm{TV}}(\mu,\nu_R)\leq||\nu_{R+1}-\nu_R||_{TV}+d_{\mathrm{TV}}(\mu,\nu_{R+1}).
$$
The first term on the right hand side of the inequality above can be bounded by (\ref{uj_totnorm}) and the second can be estimated by (\ref{uj_bizvege_dTV0}) with $R$ replaced by $R+1$. Then, by Proposition \ref{p_sigma_n}, the second inequality of Corollary \ref{kov_CPexpansion} follows.

Now we assume that $n$ and $m$ satisfy the conditions of Theorem \ref{t_PCexp_2}. In this case we can apply Theorem 3.2.~from \cite{BKN} to obtain
$$
d_{\mathrm{TV}}(\mu,\nu_R)\leq C_R\log\sigma_n\,\,\frac{n-m}{(\sqrt{n})^{R+1}}.
$$
We can apply the same argument as above to omit the $\log\sigma_n$ from the bound above. Indeed, we obtain
$$d_{\mathrm{TV}}(\mu,\nu_R)\leq C_R\frac{n-m}{(\sqrt{n})^{R+1}}+C_R\log\sigma_n\,\,\frac{n-m}{(\sqrt{n})^{R+2}}.$$
Since $\log\sigma_n\leq\sigma_n\leq\sqrt{2}(n-m)/\sqrt{n}$ by Proposition \ref{p_sigma_n}, and $(n-m)/n\leq1$, we finished the proof of the third inequality of Corollary \ref{kov_CPexpansion}. $\blacksquare$

\bigskip

Comparing the results of Corollary \ref{kov_CPexpansion} with (\ref{uj_totnorm}), we see that in the small $m$ case, when $m\leq n/2-1$, our results are not optimal in the sense that the error order of the approximation with $\nu_R$ does not coincide with the order of the total variation norm $||\nu_{R+1}-\nu_R||_{TV}$ of the $(R+1)$-th correction term. We see from the proof of Corollary \ref{kov_CPexpansion}, that the technique used in the other cases to omit the $\log\sigma_n$ from the total variation bounds does not work for small $m$, because for such $m$, $\sigma_n$ is of order $n/\sqrt{m}$, which is not comparable with the $1/(\sqrt{m})^{R+1}$ factor of the error order. It is an interesting open problem whether $d_{\mathrm{TV}}(\mu,\nu_R)\leq C_R/(\sqrt{m})^{R-1}$ can be achieved for $m\leq n/2-1$.

We finish the chapter by comparing its results with the ones obtained in Chapter 6 for compound Poisson approximation. First, we note that with some extra, but trivial considerations, the proofs of the chapter could be modified to hold true for $R=2$ also. Now, if the hypothesis formulated in the previous paragraph is true, than the bounds we obtain for $d_{\mathrm{TV}}(\mu,\nu_2)$ would exactly match the total variation bounds in (\ref{discrete_approx2}) for the compound Poisson approximation of $\mu$. We guess that if one defines $\chi$ in (\ref{uj_chi}) to be the characteristic function of the compound Poisson distribution given in Theorem \ref{t_CP}, then the technique used in this chapter with this more sensitive choice of approximating measure would lead to an improvement in terms of the error bounds.


\chapter*{Summary}
\addcontentsline{toc}{chapter}{Summary}

The coupon collector's problem is one of the classical problems of probability theory. In this thesis, we are interested in the version of the problem, when a collector samples with replacement a set of $n\ge 2$
distinct coupons so that at each time any one of the $n$ coupons is
drawn with the same probability $1/n$. For a fixed integer
$m\in\{0,1,\ldots,n-1\}$, this is repeated until $n-m$ distinct
coupons are collected for the first time. Let $W_{n,m}$ denote the
number of necessary repetitions to achieve this. Thus the random
variable $W_{n,m}$, called the coupon collector's waiting time, can
take on the values $n-m, n-m+1, n-m+2,\ldots$, and gives the number
of draws necessary to have a collection, for the first time, with
only $m$ coupons missing. In particular, $W_{n,0}$ is the waiting
time to acquire, for the first time, a complete collection.

Different limit theorems have been proved for the asymptotic distribution of $W_{n,m}$, depending on how $m=m(n)$ behaves as $n\to\infty$. (All asymptotic relations throughout are meant as $n\to\infty$.) The first result was proved by Erd\H{o}s and R\'{e}nyi for complete collections when $m=0$ for all $n\in\N$, obtaining a limiting shifted Gumbel extreme value distribution. This result was extended by Baum and Billingsley, who examined all relevant sequences of $m=m(n)$. They determined four different limiting distributions: the degenerate distribution at 0, the Poisson distribution, the normal distribution and a Gumbel-like distribution.

One of the aims of this thesis is to refine the limit theorems of Baum and Billingsley. Our basic goal is to approximate the distribution of the coupon collector's appropriately centered and normalized waiting time with well-known measures with high accuracy, and in many cases prove asymptotic expansions for the related probability distribution functions and mass functions. The approximating measures are chosen from five different measure families. Three of them -- the Poisson distributions, the normal distributions and the Gumbel-like distributions -- are probability measure families whose members occur as limiting laws in the limit theorems of Baum and Billingsley.

The fourth set of measures considered is a certain $\{\pi_{\mu,a}: \mu>0, a>0\}$ family of compound Poisson measures. For each $\mu>0$ and $a>0$, we define $\pi_{\mu,a}$ to be the probability distribution of $Z_1+2Z_2$, where $Z_1$ and $Z_2$ are independent random variables defined on a common probability space, $Z_1\sim\textrm{Po}(\mu)$ and $Z_2\sim\textrm{Po}(a/2)$.

The fifth set of approximating measures we consider is the family of Poisson--Charlier signed measures. For any positive real numbers $\lambda$, $\widetilde{a}^{(1)},\ldots,\widetilde{a}^{(S)}$ and $S\in\N$, the Poisson--Charlier signed measure $\nu=\nu(\lambda,\widetilde{a}^{(1)},\ldots,\widetilde{a}^{(S)})$ is a signed measure concentrated on the nonnegative integers defined by
$
\nu\{j\}=\mathrm{Po}(\lambda)\left(\sum_{r=1}^{S}(-1)^{r}\widetilde{a}^{(r)}C_r(j,\lambda)\right)$, $j\in\N,
$
where $C_r(j,\lambda)$ is the $r$-th Charlier polynomial.

\newpage
Our results are the following:

\begin{itemize}
  \item Chapter 3

In this chapter, we are interested in the the asymptotic behavior of the distribution function of the appropriately standardized waiting time $F_{n,m}$, if $m$ is a fixed constant for all $n$ and $n\to\infty$. With $F_m$ denoting the limiting distribution function, for every $m$, we give a one-term asymptotic expansion $F_m + G_{n,m}$ that approximates $F_{n,m}$ with the uniform order of $1/n$ such that the explicit sequence of
functions $G_{n,m}$ has the uniform order of $(\log n)/n$. We use characteristic functions in our proof.

We also give an argument, that not only proves that the error order of this approximation is sharp, but also that no longer asymptotic expansion of $F_{n,m}$ can improve the error order of $1/n$.
\medskip

  \item Chapter 4

In this chapter, with a classical characteristic function method, we prove that the error order for normal approximation to the coupon collector's standardized waiting time is at most $n/(m\sigma_n)$ in Kolmogorov distance, where $\sigma_n$ denotes the standard deviation of the waiting time. One can check that this bound is ideal in the sense that it tends to 0, iff $n$ and $m$ satisfy the conditions of the central limit theorem concerning the coupon collector's problem.
\medskip

  \item Chapter 5

In the first section of Chapter 5, we consider Poisson approximation to the distribution of sums of asymptotically negligible integer valued random variables in general. We complement a classical Poisson convergence theorem of Gnedenko and Kolmogorov. Considering an arbitrary triangular array $\{Y_{n1}, Y_{n2}, \ldots, Y_{nr_n}\}_{n\in\N}$ of row-wise independent nonnegative integer valued random variables, for each~$n$, we approximate the distribution of the $n$-th row sum with a Poisson distribution whose mean~$\lambda_n$ is defined only in terms of the distributions of the random variables in the~$n$-th row, but we do not assume the existence of moments. We give both lower and upper bounds, which have precisely the same form, up to a constant, provided that the means $\lambda_n$ are bounded away from infinity. We thus refine the obvious approximation of the $Y_n$-s that the Gnedenko-Kolmogorov limit theorem suggests.

In the next section we examine how the coupon collector's problem fits in the framework of the previous section. We show that the Poisson limit theorem concerning the coupon collector's waiting time is a special case of the Gnedeno--Kolmogorov theorem. Applying the general results of the previous section to the waiting time, we obtain a Poisson approximation of error order $1/\sqrt{n}$.

In the third section of the chapter we take advantage of the combinatorial structure of the coupon collector's problem. This combinatorial approach yields us a stronger result than the one of the previous section: we derive the first asymptotic correction of the $\p(W_{n,m}-(n-m)=k)$, $k=0,1,\ldots$, probabilities to the corresponding Poisson point probabilities.

In the final section of Chapter 5, we approximate the coupon collector's shifted waiting time $\widetilde{W}_{n,m}=W_{n,m}-(n-m)$ with another Poisson law, namely with the one that has the same mean as $\widetilde{W}_{n,m}$. One can easily calculate that in the range of parameters $n$ and $m$ for which the Poisson limit theorem holds true, the error order of this new approximation is $1/n$, which is clearly better than the error order $1/\sqrt{n}$ given in the preceding two sections for the same case.  The proof here is based on Stein's method, and heavily uses the fact that the means of the compared probability measures coincide.
\medskip

  \item Chapter 6

In the first section of Chapter 6 we consider translated compound Poisson approximation of sums of independent integer valued random variables in general. Using Stein's method, Barbour and others gave bounds for the errors of such approximations in total variation distance. Their upper bounds are expressed with the help of the first three moments of the summands $X_1,X_2,\ldots,X_n$ and the critical ingredient $d_{\mathrm{TV}}\left({\cal D}\left(W_n\right),{\cal D}\left(W_n+1\right)\right)$, where $W_n=\sum_{j=1}^{n}X_j$.

The expression $d_{\mathrm{TV}}\left({\cal D}\left(W_n\right),{\cal D}\left(W_n+1\right)\right)$
is usually bounded by the Mineka coupling, which typically yields a
bound of order $1/\sqrt{n}$. If the $X_j$'s are roughly similar in magnitude, this is
comparable with the order $O(1/\sqrt{\var W_n})$ expected for the error in the central
limit theorem.  However, if the distributions of the~$X_j$ become progressively more
spread out as~$j$ increases, then $\var W_n$ may grow faster than~$n$, and then
$1/\sqrt n$ is bigger than the ideal order $O(1/\sqrt{\var W_n})$. In fact, this is the situation in the case when we chose $W_n$ to be the coupon collector's waiting time. We introduce a new coupling which allows us to improve the Mineka bounds in such settings.

In the next section, we approximate the distribution of the appropriately centered coupon collector's waiting time with a compound Poisson measure $\pi_{\mu,a}$ defined above. We apply general results of translated compound Poisson approximation of sums of independent integer valued random variables and our new coupling. We prove that a translated compound Poisson approximation to the collector's waiting time, with ideal
error rate, can be obtained in all ranges of $n$ and~$m$ in which a central or Poisson limit theorem holds true.
Comparing this result with the ones we obtained for normal approximation, we see that the same or better order of approximation is obtained with this discrete approximation, and now with the error measured with respect to the much stronger total variation distance.
\medskip

  \item Chapter 7

In the final chapter of the thesis we approximate the coupon collector's shifted waiting time $\widetilde{W}_{n,m}=W_{n,m}-(n-m)$ with Poisson--Charlier signed measures in total variation distance. To do so, we apply a characteristic function technique. For an arbitrary $R\geq3$, we choose a Poisson--Charlier signed measure $\nu_R$ depending on $R$. Approximating $\widetilde{W}_{n,m}$ with $\nu_R$, we obtain error bounds of
order $\sigma_n\log\sigma_n(1/\sqrt{m})^{R}$, $(\sqrt{n})^{R-3}/(n-m)^{R-2}$ and $(n-m)/(\sqrt{n})^{R+1}$ depending on how the sequences $m/n$ and $\sigma_n^2-\mu_-(n-m)$ behave as $n$ tends to infinity, where $\mu_n$ is the mean and $\sigma_n^2$ is the variance of $W_{n,m}$.
\end{itemize}


\chapter*{Összefoglalás}
\addcontentsline{toc}{chapter}{Összefoglalás}

A kupongyűjtő probléma a valószínűségszámítás egyik klasszikus problémája. A dolgozatban a probléma következő változatával foglalkozunk: adva van $n\ge 2$ különböző kupon, melyekből egy gyűjtő véletlenszerű visszatevéses mintát vesz úgy, hogy minden egyes alkalommal az $n$ kupon bármelyikét azonos, tehát $1/n$ valószínűséggel szerzi meg. Valamely rögzített $m\in\{0,1,\ldots,n-1\}$ esetén a mintavételt addig folytatja, amíg előszörre pontosan $n-m$ különböző kupont nem gyűjtött. Jelölje $W_{n,m}$ az ehhez szükséges ismétlések számát. Tehát a $W_{n,m}$ véletlen változó, amelyet a kupongyűjtő várakozási idejének nevezünk, az $n-m, n-m+1, n-m+2,\ldots$ értékeket veheti fel, és megadja, hogy a gyűjtőnek hányszor kell húznia ahhoz, hogy $m$ darab kupon kivételével minden kupon a birtokában legyen. Speciálisan, $W_{n,0}$ a teljes gyűjtemény megszerzéséhez szükséges várakozási idő.

Különböző határeloszlás tételeket bizonyítottak $W_{n,m}$ aszimptotikus eloszlására attól függően, hogy az $m=m(n)$ sorozat hogyan viselkedik, amint $n\to\infty$. (A továbbiakban minden konvergencia és aszimptotikus reláció $n\to\infty$ mellett értendő.) Az első eredmény Erd\H{o}s és R\'{e}nyi nevéhez fűződik, akik a teljes gyűjtemény esetére, amikor minden $n\in\N$ esetén $m=0$, eltolt Gumbel-eloszlást kaptak határeloszlásként. Ezt az eredményt általánosította Baum és Billingsley, akik minden $m=m(n)$ sorozat típust vizsgáltak. Négy különböző határeloszlást határoztak meg: a 0-ra koncentrált eloszlást, a Poisson eloszlást, a normális eloszlást, és egy a Gumbel-eloszlásból származtatható eloszlást.

A dolgozatban  finomítjuk a fenti, Baum és Billingsley nevéhez fűződő határeloszlás tételeket. Célunk a megfelelően centralizált és normalizált várakozási idő eloszlásának jól ismert mértékekkel történő minél pontosabb közelítése, és sok esetben a kapcsolódó eloszlásfüggvények és valószínűségi pontfüggvények aszimptotikus sorfejtése. A közelítő mértékeket öt különböző mértékcsaládból választjuk. Ezek közül három -- a Poisson eloszlások, a normális eloszlások és a Gumbel-típusú eloszlások -- olyan mértékcsaládok, melyeknek tagjai határeloszlásaként szerepelnek Baum és Billingsley tételeiben.

A negyedik approximáló mértékcsalád az összetett Poisson eloszlásoknak bizonyos $\{\pi_{\mu,a}: \mu>0, a>0\}$ osztálya. Tetszőleges $\mu>0$ és $a>0$ esetén $\pi_{\mu,a}$ a $Z_1+2Z_2$ véletlen változó eloszlását jelöli, ahol $Z_1$ és $Z_2$ valamely közös valószínűségi mezőn definiált független véletlen változók, $Z_1\sim\textrm{Po}(\mu)$ és $Z_2\sim\textrm{Po}(a/2)$. Az ötödik közelítő mértékcsalád a Poisson--Charlier előjeles mértékek osztálya. Tetszőleges pozitív valós $\lambda$, $\widetilde{a}^{(1)},\ldots,\widetilde{a}^{(S)}$ és $S\in\N$ esetén a $\nu=\nu(\lambda,\widetilde{a}^{(1)},\ldots,\widetilde{a}^{(S)})$  Poisson--Charlier előjeles mérték az az előjeles mérték, amely a nemnegatív egészekre van koncentrálva, és $\nu\{j\}=\mathrm{Po}(\lambda)\left(\sum_{r=1}^{S}(-1)^{r}\widetilde{a}^{(r)}C_r(j,\lambda)\right)$, $j\in\N$, ahol $C_r(j,\lambda)$ az $r$-edik Charlier polinom.

\newpage
Eredményeink a következőek:

\begin{itemize}
  \item 3.~fejezet

A 3.~fejezetben a megfelelően standardizált várakozási idő $F_{n,m}$ eloszlásfüggvényének aszimptotikus viselkedését vizsgáljuk abban az esetben, amikor $m$ minden $n$-re rögzített konstans és $n\to\infty$. Jelölje $F_m$ a határeloszlás függvényt. Minden $m$ esetén olyan $F_m + G_{n,m}$ egytagú aszimptotikus sorfejtést adunk, amely egyenletesen $1/n$ rendben közelíti az $F_{n,m}$ eloszlásfüggvényt, továbbá az explicit módon megadott $G_{n,m}$ függvények sorozata egyenletesen $(\log n)/n$ rendű. A bizonyításban karakterisztikus függvényeket használunk.

Azt is belátjuk, hogy a közelítés hibarendja éles, és hogy $F_{n,m}$-nek semmilyen hosszabb aszimptotikus sorfejtése esetén nem kaphatunk az $1/n$ rendnél kisebb hibarendet.
\smallskip

  \item 4.~fejezet

A dolgozat ezen fejezetében klasszikus, karakterisztikus függvényeket használó módszerrel belátjuk, hogy a kupongyűjtő várakozási idejének normális eloszlással való közelítésének Kolmogorov távolságban mért hibája legfeljebb $n/(m\sigma_n)$ rendű, ahol $\sigma_n$ a várakozási idő szórását jelöli. Ellenőrizhető, hogy ez az approximáció jó abban az értelemben, hogy 0-hoz tarta, ha $n$ és $m$ teljesítik a kupongyűjtő problémára vonatkozó centrális határeloszlás tétel feltételeit.
\smallskip

  \item 5.~fejezet

Az 5.~fejezet első részében általános független nemnegatív egészértékű véletlen változók összegeinek eloszlását közelítjük Poisson eloszlással. A Gnedenko és Kolmogorov nevéhez fűződő klasszikus Poisson határeloszlás tételt finomítjuk.
Tetszőleges, soron-ként független, nemnegatív egészértékű véletlen változókból álló $\{Y_{n1}, Y_{n2}, \ldots, Y_{nr_n}\}_{n\in\N}$ szériasorozatot tekintve minden $n$ esetén az $n$-edik sorösszeg eloszlását olyan Poisson eloszlással közelítjük, melynek $\lambda_n$ várható értéke csak az adott sorban szereplő változók eloszlásától függ, de nem követeljük meg momentumok létezését.
A közelítés hibájára alsó és felső korlátot is adunk, melyek rendje konstans szorzótól eltekintve megegyezik, feltéve, hogy a $\lambda_n$ paraméterek korlátosak. Ezáltal jobb közelítését adjuk az $Y_n$ változóknak, mint amit a kézenfekvő, határeloszlással történő approximáció jelent.

A fejezet második alfejezetében megmutatjuk, hogy a várakozási időre vonatkozó Poisson határeloszlás tétel speciális esete a Gnedeno--Kolmogorov tételnek. Megmutatjuk, hogy ha az előző alfejezet eredményeit alkalmazzuk a várakozási időre, $1/\sqrt{n}$ hibarendű Poisson közelítését kapjuk.

A harmadik alfejezetben a kupongyűjtő probléma kombinatorikai struktúrájára építünk. A kombinatorikai megfontolásokra támaszkodó módszer segítségével erősebb eredményt tudunk igazolni, mint az előző alfejezetben: a $\p(W_{n,m}-(n-m)=k)$, $k=0,1,\ldots$, valószínűségek megfelelő Poisson valószínűségekkel történő közelítését pontosítjuk az első korrekciós tag meghatározása révén.

Az 5.~fejezet utolsó alfejezetében a kupongyűjtő $\widetilde{W}_{n,m}=W_{n,m}-(n-m)$ eltolt várakozási idejét egy újabb Poisson eloszlású véletlen változóval közelítjük, méghozzá olyannal, amelynek várható értéke megegyezik $\widetilde{W}_{n,m}$ várható értékével. Kiszámolható, hogy azon $n$ és $m$ paraméter értékek esetén, melyekre érvényes a Poisson határeloszlás tétel, az approximáció hibájának rendje $1/n$, ami világos, hogy kisebb, mint az előző két alfejezetben ugyanezen esetre bizonyított $1/\sqrt{n}$-es hibarend.  A bizonyítás a Stein-módszeren alapszik, és kihasználja azt a tényt, hogy az összehasonlított eloszlások várható értékei egyenlőek.
\smallskip

  \item 6.~fejezet

A hatodik fejezet első alfejezetében független egészértékű véletlen változók összegeinek összetett Poisson eloszlású változókkal történő közelítését vizsgáljuk általánosan. A Stein-módszer segítségével Barbour és mások felső korlátokat adtak az ilyen típusú approximációk teljes variációs távolságban mért hibáira. Ezek a korlátok az $X_1,X_2,..,X_n$ összeadandók első három momentumának és a $d_{\mathrm{TV}}\left({\cal D}\left(W_n\right),{\cal D}\left(W_n+1\right)\right)$ kifejezésnek függvényei, ahol $W_n=\sum_{j=1}^{n}X_j$.

A $d_{\mathrm{TV}}\left({\cal D}\left(W_n\right),{\cal D}\left(W_n+1\right)\right)$ kifejezést általában a Mineka-csatolás segítségével lehet becsülni, ami tipikusan $1/\sqrt{n}$ rendű eredményt ad. Ha az $X_j$ véletlen változók nagyjából azonos szórásúak, akkor ez az eredmény közel van a centrális határeloszlás tétel esetén elvárt $O(1/\sqrt{\var W_n})$ hibarendhez. Azonban ha az $X_j$ véletlen változók eloszlásai egyre laposabbak, amint $j$ nő, akkor $\var W_n$ nőhet gyorsabban, mint $n$, és ekkor $1/\sqrt n$ jóval nagyobb lesz, mint az elvárható $1/\sqrt{\var W_n}$-es rend. Pontosan ez a helyzet, ha $W_n$-nek a kupongyűjtő várakozási idejét választjuk. Bevezetünk egy új csatolást, mely segítségével ilyen esetekben jobb eredményeket tudunk bizonyítani, mint a Mineka-csatolás segítségével.

A következő alfejezetben a kupongyűjtő megfelelően centralizált várakozási idejének eloszlását a korábban definiált $\pi_{\mu,a}$ összetett Poisson eloszlással közelítjük. Független egészértékű véletlen változók összegeire vonatkozó általános összetett Poisson approximációs eredményeket és az új csatolásunkat alkalmazzuk. Belátjuk, hogy a $W_{n,m}$ várakozási idő jól közelíthető összetett Poisson eloszlással abban az esetben, amikor az $n$ és $m$ paraméterek teljesítik a kupongyűjtő problémára vonatkozó centrális vagy Poisson határeloszlás tétel feltételeit. Ezeket és a normális approximációra kapott eredményeket összehasonlítva látjuk, hogy az itt bevezetett diszkrét approximáció ugyanolyan, vagy jobb közelítését jelenti a várakozási időnek, mint a normális approximáció. Ráadásul a közelítés hibáját itt a Kolmogorov távolságnál sokkal erősebb teljes variációs távolságban mérjük.
\smallskip

  \item 7.~fejezet

Az utolsó fejezetben a kupongyűjtő $\widetilde{W}_{n,m}=W_{n,m}-(n-m)$ eltolt várakozási idejét Poisson--Charlier előjeles mértékekkel közelítjük teljes variációs távolságban. Ehhez egy karakterisztikus függvényeket használó módszert alkalmazzuk. Tetszőleges $R\geq3$ esetén definiálunk egy $R$-től függő $\nu_R$ Poisson--Charlier előjeles mértéket. Az eltolt várakozási időt ezzel közelítve $\sigma_n\log\sigma_n(1/\sqrt{m})^{R}$, $(\sqrt{n})^{R-3}/(n-m)^{R-2}$, illetve $(n-m)/(\sqrt{n})^{R+1}$ hibarendeket kapunk aszerint, hogy az $m/n$, illetve a $\sigma_n^2-\mu_n-(n-m)$ sorozatok hogyan viselkednek, amint $n\to\infty$, ahol $\mu_n$ a várakozási idő várható értékét, $\sigma_n^2$ pedig a szórásnégyzetét jelöli.
\end{itemize}


\chapter*{Acknowledgement}
\addcontentsline{toc}{chapter}{Acknowledgement}

I am heartily thankful to my supervisors Sándor Csörgő and Andrew Barbour.
Professor Csörgő was a remarkable teacher who evoked my affection for probability theory during his lectures. He introduced me to mathematical research and gave me a lot of encouragement. I am grateful for all I have learnt form him and for the personal guidance he gave me at the start of my academic career.
I would also like to thank Professor Barbour, who volunteered to be my advisor after the sudden and untimely death of Professor Csörgő. He introduced me to the area of probability theory I enjoy the most. I am grateful for his hospitality which I enjoyed many times at the University of Zurich, and for all the illuminating e-mails he sent. This thesis would not have been possible without his help.
I am also thankful to Vilmos Totik and Gyula Pap for their useful suggestions regarding my work and this thesis.



\begin{thebibliography}{100}
\addcontentsline{toc}{chapter}{Bibliography}

\bibitem{BD} Banderier, C.~and Dobrow, R.~P., A Generalized Cover Time for Random Walks on Graphs, Proceedings of FPSAC'00, 2000.

\bibitem{BC} Barbour, A.D.~and Cekanavicius, V, Total variation asymptotics for sums of independent integer random variables, {\it The Annals of Probability} {\bf 30} (2002), 509--545.

\bibitem{BChen} Barbour, A.~D.~and Chen, L.~H.~Y., {\it An Introduction to Stein's method}, Lecture Notes Series, Institute for Mathematical Sciences, National University of Singapore, Vol.~4, 2005.

\bibitem{BH} Barbour, A.D.~and Hall, P., On the rate of Poisson convergence, {\it Math. Proc. Cam. Phil. Soc.} {\bf 95} (1984), 473--480.

\bibitem{BHJ} Barbour, A.D., Holst, L. and Janson, S., {\it Poisson Approximation}, Clarendon Press, Oxford, 1992.

\bibitem{BKN} Barbour, A.~D., Kowalski, E., Nikeghbali, A., {\it Mod-discrete expansions}, {\it arXiv: 0912.1886v1 [math.PR]}, 2009.

\bibitem{BX} Barbour, A.D.~and Xia, A., Poisson Perturbations, {\it ESAIM Probab. and Statist.} {\bf 3} (1999), 131--150.

\bibitem{BB} Baum, L.E.~and Billingsley, P., Asymptotic distributions for the coupon collector's problem, {\it Ann. Math. Statist.} {\bf 36} (1965), 1835--1839.

\bibitem{Ch} Chihara, T.~S., {\it An introduction to orthogonal polynomials}, Gordon and Breach, New York, 1978.

\bibitem{C} Chen, L.~H.~Y., Poisson approximation for dependent trials, {\it Ann.~Probab.} {\bf 3} (1975), 534-545.

\bibitem{CS} Cs\"{o}rg\H{o}, S., A rate of convergence for coupon collectors, {\it Acta Sci. Math. (Szeged)} {\bf 57} (1993), 337--351.

\bibitem{ER} Erd\H{o}s, P. and R\'enyi, A., On a classical problem of probability theory, {\it Magyar Tud. Akad. Mat. Kutat\'o Int. K\"ozl.} {\bf 6} (1961), 215--220.

\bibitem{F} Feller, W., {\it An Introduction to Probability Theory and its Applications}, John Wiley $\&$ Sons, 1968.

\bibitem{GS} Gibbs, A.~L. and Su, F.~E., On Choosing and Bounding Probability Metrics, {\it International Statistical Review / Revue Internationale de Statistique}, Vol.~70, No.~3 (Dec., 2002), 419--435.

\bibitem{GK} Gnedenko, B. V. and Kolmogorov, A. N., {\it Limit Distributions for Sums of Independent Random Variables}, Addison-Wesley Publishing Company, Cambridge, Mass., 1954.

\bibitem{GH} Gut, A.~and Holst, L., On the waiting time in a generalized roulette game, {\it Statistics \& Probability Letters}, {\bf 2} (1984), 229--239.

\bibitem{Holst} Holst, L., The general birthday problem, Proceedings of the sixth international seminar on Random graphs and probabilistic methods in combinatorics and computer science, John Wiley \& Sons, 1995.

\bibitem{H} Huber, P.~J., {\it Robust Statistics}, John Wiley \& Sons, New York, 1981.

\bibitem{Cam} Le Cam, L., An approximation theorem for the Poisson binomial distribution, {\it Pacific J. Math.} {\bf 10} (1960), 1181--1197.

\bibitem{L} Lindvall, T., {\it Lectures on the Coupling Method}, Dover Publications, 1992.

\bibitem{MR} Mattner, L. and Roos, B., A shorter proof of Kanter's Bessel function concentration bound, {\it Probability Theory and Related Fields} {\bf 139} (2006), 191--205.

\bibitem{N} Nielsen, N., {\it Handbuch der Theorie der Gammafunktion}, Teubner, Leipzig, 1906. [Re\-printed as Band I of {\it Die Gammafunktion}, Chelsea, New York, 1965.]

\bibitem{P_limit} Petrov, V.~V., {\it Limit theorems of probability theory}, The Clarendon Press Oxford
    University Press, New York, 1995.

\bibitem{P_sums} Petrov, V.V., {\it Sums of Independent Random Variables}, Springer-Verlag, Berlin, 1975.

\bibitem{Polya} P\'olya, G., Eine Wahrscheinlichkeitsaufgabe zur Kundenwerbung, {\it Z. Angew. Math. Mech.}, {\bf 10} (1930), 96--97.

\bibitem{P_cp} P\'osfai A., An extension of Mineka's coupling inequality, {\it Electronic Communications in Probability}, {\bf 14} (2009), 464--473.

\bibitem{ParXiv} P\'osfai, A., A supplement to the paper Poisson approximation in a Poisson limit theorem inspired by coupon collecting, {\it arXiv:0904.4924 [math.PR]}, 2009.

\bibitem{P_p} P\'osfai A., Poisson approximation in a Poisson limit theorem inspired by coupon collecting, {\it Journal of Applied Probability}, {\bf 46} (2009), 585--592.

\bibitem{P_n} P\'osfai, A., Rates of convergence for normal approximation in incomplete coupon collection,
{\it Acta Scientiarum Mathematicarum (Szeged)} {\bf 73} (2007), 333--348.

\bibitem{PCS} P\'osfai, A. and Cs\"org\H{o}, S., Asymptotic approximations for coupon collectors, {\it Studia Scientiarum Mathematicarum Hungarica}, {\bf 46} (2009), 61--96.

\bibitem{R} Rachev, S.~T., {\it Probability Metrics and the Stability of Stochastic Models}, Wiley, 1991.

\bibitem{S} Stein, C., A bound for the error in the normal approximation to the distribution of a sum of dependent random variables, {\it Proceedings of the Sixth Berkeley Symposium on Mathematical Statististics and Probability}, Vol.~2, Univ.~of Calif.~Press, 1972.

\bibitem{T} Thorisson, H., {\it Coupling, Stationarity and Regeneration}, Springer, 2000.
\end{thebibliography}
\end{document}